\documentclass[a4paper, 12pt]{article}
\usepackage{amssymb}
\usepackage{amsmath}
\usepackage{amscd}
\usepackage{latexsym}
\usepackage{theorem}
\usepackage{color}
\theoremheaderfont{\scshape}
\newtheorem{thm}{\bfseries Theorem}[section]
\newtheorem{prop}[thm]{\bfseries Proposition}
\newtheorem{lemma}[thm]{\bfseries Lemma}
\newtheorem{cor}[thm]{\bfseries Corollary}

\theorembodyfont{\rmfamily}
\newtheorem{defn}[thm]{\bfseries Definition}
\newtheorem{ex}[thm]{Example}
\newtheorem{rem}[thm]{Remark}

\newtheorem{conj}[thm]{Idea}
\newtheorem{onepara}[thm]{}
\newtheorem{pf}{Proof.}

\def\N{{\mathbb N}}
\def\Z{{\mathbb Z}}

\def\m{\mathfrak m}
\def\n{\mathfrak n}
\def\ga{\mathfrak a}
\def\Hom{\mathrm{Hom}}

\def\Tor{\mathrm{Tor}}
\def\Ext{\mathrm{Ext}}
\def\length{\mathrm{length}}
\def\End{\mathrm{End}}

\def\Kdim{\mathrm{dim}}

\def\Spec{\mathrm{Spec}}
\def\kernel{\mathrm{Ker}}

\def\A{{{\mathcal A}_k}}
\def\Ahat{{\widehat{\mathcal A}_k}}
\def\C{{{\mathcal C}_k}}

\def\T{\mathcal T}

\def\F{\mathcal F}
\def\I{\mathcal I}

\def\L{\mathbb L}
\def\Lx{\mathbb L}
\def\Fx{\mathbb F}
\def\Gx{\mathbb G}
\def\Xx{\mathbb X}

\def\ctensor{{\ {\widehat{\otimes}_k}\ }}
\def\formal{k \langle\langle t_1, t_2 , \ldots , t_r \rangle\rangle}  

\def\qed{$\Box$}

\title{\bf 
 Universal lifts of chain complexes over non-commutative parameter algebras 
\\}
\author
{Yuji Yoshino \\
{\small Okayama University, } 
{\small Okayama 700-8530, Japan} \\
{\small \texttt{yoshino@math.okayama-u.ac.jp}}
}

\date{\empty}
\begin{document}
\maketitle


\begin{abstract}
We define the notion of universal lift of a projective complex based on non-commutative parameter algebras, and prove its existence and uniqueness.   
We investigate the properties of parameter algebras for universal lifts.

\bigskip

{\footnotesize{
\par 
\indent 
\hphantom{Keywords ::} non-commutative formal power series ring
}}

\end{abstract}

\bigskip

\begin{center}
{\bf Contents}
\end{center}
\begin{itemize}
\item[1] Introduction
\item[2] Non-commutative complete local algebras 
\begin{itemize}
\item[2.1] 
Definitions and properties 
\item[2.2] 
Small extensions
\item[2.3] 
Complete tensor products
\end{itemize}
\item[3] Universal lifts of chain complexes 
\begin{itemize}
\item[3.1]
Lifts to artinian local algebras 
\item[3.2]
Construction of maximal lifts
\item[3.3] 
Universal lifts
\item[3.4]
Every complete local algebra is a parameter algebra 
\item[3.5]
Deformation of modules 
\end{itemize}
\item[4] Properties of parameter algebras 
\begin{itemize}
\item[4.1] 
Obstruction maps 
\item[4.2] 
Universal lifts based on commutative algebras 
\item[4.3] 
Yoneda products
\item[4.4] 
Comparison of cohomology
\end{itemize}
\end{itemize}


\section{Introduction}

In this paper,  $k$ always denotes a field and $R$ is an arbitrary associative $k$-algebra. 
When we say an $R$-module, we always mean a left $R$-module unless otherwise stated.  

From the view point of representation theory, the final goal of the theory of $R$-modules should be to construct the moduli consisting of the isomorphism classes of  $R$-modules, by which we mean a geometric realization of the set of isomorphism classes. 
Generally speaking, it is however impossible to describe all the isomorphism classes of $R$-modules, even if we restrict ourselves to consider indecomposable ones. 
And one should say that the construction of moduli for $R$-modules is hopeless.

But there is a way to observe the moduli from the local view point. 
Fixing an $R$-module $M$, and assuming there is a modulus containing $M$  as a rational closed point, we can ask how it looks in the neighbourhood of the point,  which is nothing but to consider the universal deformation of  $M$. 
In such a context, the existence of formal local moduli is known (\cite{FGA}, \cite{H}, \cite{S}).

To explain this, let  $\C$  be the category of commutative artinian local $k$-algebras with residue field  $k$ and $k$-algebra homomorphisms. 
We consider the covariant functor 
$$
\F_M : \C \to (Sets), 
$$
which maps  $A \in \C$  to  the set of infinitesimal deformations of $M$ along  $A$, i.e. 
$$
\F_M (A) =  \left\{ 
\begin{aligned}
&\text{$(R,A)$-bimodules $X$ that are flat over $A$}\\  
&\text{　and $X \otimes^{L} _A k \cong M$  as left $R$-modules} 
\end{aligned}
\right\} \ / \cong,   
$$
where  $\cong$ means $(R, A)$-bimodule isomorphism. 
Under these circumstances the following theorem is known to hold.

\begin{thm}[Schlessinger's Theorem 1968]\label{schlessinger} 
Suppose $\Ext _R ^1(M,M)$ is of finite dimension as a $k$-vector space.
Then the functor  $\F_M$  is pro-representable. 
More precisely, there exist a commutative noetherian complete local $k$-alegbra  $Q$ with residue field  $k$ and an $(R, Q)$-bimodule  $U$  that is  flat over  $Q$  such that there is an isomorphism   
$$
\Hom _{k\text{-alg}} (Q, \ \ ) \cong  \F_M 
$$
as functors on  $\C$.
The isomorphism is given in such a way that  each  $f \in \Hom _{k\text{-alg}}(Q,  A)$  is mapped to  $[U \otimes _Q {}_fA] \in \F_M(A)$  for $A \in \C$, 
where  ${}_fA$ denotes the right $A$-module  $A$  regarded as a left $Q$-module  through  $f$. 
\end{thm}

In such a circumstance, we call $U$ the universal family of deformations of $M$, and call $Q$ (resp. $\Spec \ Q$) the commutative parameter algebra (resp. the parameter space) of  $U$. 

One of the easiest examples is the deformation of Jordan canonical forms. 

\begin{ex}\label{jordan}
Consider an $n \times n$ matrix which is of an irreducible Jordan canonical form:
$$
\begin{pmatrix}
0 & 1 & 0 &\cdots & 0 \\
0 & 0 & 1 &\cdots & 0 \\
: & : & : & : & : \\
: & :& : & : & 1 \\
0 & 0 & 0 &\cdots & 0 \\
\end{pmatrix}
$$
Setting  $R=k[x]$, we know that this is equivalent to consider the indecomposable  $R$-module $M = k[x]/(x^n)$. 
In this case, we can take  $Q = k[[ t_0,  \ldots , t_{n-1}]] $ as the commutative parameter algebra，and  $U = Q[x] / (x^n + t_{n-1}x^{n-1} + \cdots + t_0 )$ as the universal family of deformations of $M$.
If we consider this in a matrix form,  we obtain a so-called Sylvester family of matrices.
$$
\begin{pmatrix}
0 & 1 & 0 &\cdots & 0 \\
0 & 0 & 1 &\cdots & 0 \\
: & : & : & : & : \\
: & :& : & : & 1 \\
-t_0 & -t_1 & -t_2 &\cdots & -t_{n-1} \\
\end{pmatrix}
$$
\end{ex}

\bigskip

Under the setting of Theorem \ref{schlessinger}, since the bimodule $U$ is flat as a right $Q$-module,  the functor  
$
U \otimes^{L} _Q \ - \ \ : D(Q) \to D(R)
$
between derived categories is defined. 
Remark that  $U \otimes ^L_Q k = M$. 
Thus the functor induces  a map between Yoneda algebras. 
$$
\rho ^{\cdot} \ : \ \Ext ^{\cdot} _Q(k, k) \to 
 \Ext ^{\cdot} _R(M, M). 
$$
Of most interest is the mapping 
$$
\rho ^{2} \ : \ \Ext ^{2} _Q(k, k) \to 
 \Ext ^{2} _R(M, M), 
$$
which is often called the obstruction map. 
Our motivation of this paper starts with the observation that $\rho ^2$ does not work well as a comparison map between cohomology modules. 
We show this by the above example. 
In fact, we see in Example \ref{jordan} that 
$$
\begin{array}{rl} 
\Ext ^{2} _Q(k, k) &= \langle \text{Koszul relations of degree $2$ in the variables $t_i$'s} \rangle ^*  
\vspace{6pt} \\
&\downarrow \ \ \rho ^{2} 
\vspace{6pt} \\
\Ext ^{2} _R(M,&M) =(0). \\
\end{array} 
$$
Compared with that $\Ext ^{2} _R(M, M) = (0)$,  the  $k$-vector space  $\Ext ^{2} _Q(k, k)$ has dimension $n(n-1)/2$. 
 This is one of the examples that shows that $\rho ^2$ does not work well as  a comparison map of cohomology modules.
Here we should notice that the Koszul relations of degree $2$  are derived from the commutativity relations of the variables  $t_0,\ldots, t_{n-1}$. 

Thinking this phenomenon over, we get the idea that 
the parameters $t_0, \ldots , t_{n-1}$ should be regarded as non-commutative variables. 
Now we propose the following idea. 

\begin{conj}\label{1st conjecture} 
Parameter algebras should be non-commutative. 
\end{conj}

If we simply generalize the arguments in the commutative setting, 
we will have difficulty in showing the flatness of the universal family of deformations over the non-commutative parameter algebra. 
The reason for this is that the local criterion of flatness does not necessarily hold for modules over non-commutative rings. 
Therefore, to avoid the argument about flatness, we also propose the following idea.

\begin{conj}
We should consider the deformation of chain complexes instead of modules. 
\end{conj}

The deformation of chain complexes is nothing but the lifting of complexes, 
which we mainly discuss in this paper. 
In such a way, we necessarily come to think of \lq\lq the universal lifts of chain complexes over non-commutative parameter algebras\rq\rq.

\vspace{12pt}

Just to explain about the lifting of chain complexes, let us introduce several   notation concerning chain complexes. 
When we say  $\Fx = (F,  d)$ is a chain complex (or simply a complex) of $R$-modules, we mean that  $F = \oplus  _{i \in \Z} F_i$ is a graded $R$-module and 
$d : F \to F[-1]$ is a graded homomorphism satisfying $d^2 =0$. 
A projective complex  $\Fx = (F,d)$ is just a complex where the underlying graded module  $F$ is a projective $R$-module. 
If $\Fx = (F, d)$ is a projective complex, then we define 
$\Ext_R ^i(\Fx, \Fx)$ to be the set of homotopy equivalence classes of chain homomorphisms  on  $\Fx$  of degree $-i$.

We introduce the category $\A$,  
whose objects are artinian local $k$-algebras with  residue field $k$ with $k$-algebra homomorphisms as morphisms.
(Note that an object of $\A$ is not necessarily a commutative ring, but it is a finite dimensional $k$-algebra.)
Now let  $A \in \A$ and let  $\Fx =(F, d)$ be a projective complex of $R$-modules. 
Then, $(F \otimes _k A, \Delta)$ is said to be a lift of $\Fx$ to $A$ if 
it is a chain complex of $R \otimes _k A^{op}$-modules, and satisfies the equality  $\Delta \otimes _A k = d$.

The aim of this paper is to construct the universal lift of a given projective  complex  $\Fx = (F, d)$  which dominates all the lifts of  $\Fx$  to all non-commutative artinian $k$-algebras in  $\A$, and to investigate the properties of its parameter algebra. 

We should note that such a universal lift is no longer defined on an artinian algebra, but defined on a \lq pro-artinian\rq \ local $k$-algebra. 
We call such a pro-artinian algebra a complete local $k$-algebra by an abuse of  the terminology for commutative rings. 
The non-commutative formal power series ring  $k\langle\langle t_1, \ldots , t_r\rangle\rangle$ with non-commutative variables  $t_1, \ldots , t_r$ is an example of complete local $k$-algebra. 
This is actually complete and separated in the $(t_1, \ldots , t_r)$-adic topology. 
And a complete local $k$-algebra is defined to be a residue ring of the  non-commutative formal power series ring by a closed ideal. 
(See Definition \ref{def complete local algebra} and Proposition \ref{image of T}.) 
In particular all artinian algebras in  $\A$  are complete local $k$-algebras. 
But the difficulty here is that complete local $k$-algebras are not necessarily  noetherian rings.

We can extend the notion of lifting to the lifting to complete local $k$-algebras. 
In fact,  $(F \widehat\otimes _k A , \Delta _A)$ is said to be a lift of  $\Fx$ to  a complete local $k$-algebra $A$  if it is  a chain complex of $R \widehat\otimes _k A^{op}$-modules and the equality $\Delta _A \otimes _A k = d$ holds. 
(See Section 2.3 for the complete tensor product  $\widehat\otimes$.) 

To give a precise definition of universal lifts, 
let $\Fx = (F, d)$ be a projective complex of $R$-modules which we fix. 
Then we define a covariant functor $\F : \A  \to (Sets)$  by setting as $\F (A)$  the set of chain-isomorphism classes of lifts of  $\Fx$ to $A$  for any  $A \in \A$. 
If we have a complete local $k$-algebra  $P$ and 
 a lift $\Lx = (F \widehat\otimes _k P  , \Delta _P)$ of $\Fx$ to  $P$, 
then we can define a natural transformation  
$ \phi _{\Lx} : \Hom _{k\text{-alg}} (P, - ) \to \F $
of functors by setting 
$
\phi _{\Lx} (f) = (F\otimes _k A, \Delta _P \otimes _P {}_fA) 
$
for $A \in \A$ and $f \in \Hom _{k\text{-alg}} (P, A)$, where 
  ${}_fA$ denotes the right $A$-module  $A$  regarded as a left $P$-module  through  $f$.

A chain complex  $\Lx = (F \widehat\otimes  _k P , \Delta _P)$ is said to be a universal lift of  $\Fx$, if $\phi _{\Lx}$ is an isomorphism of functors. 
In this case, we say that $P$ is a parameter algebra.

The first main result of this paper is about the existence and the uniqueness of universal lifts, which we summarize as follows. 
(See Theorem \ref{uniqueness of max lifts} and Theorem \ref{univ=max}.)

\begin{thm}\label{summarize}
Let  $\Fx = (F, d)$ be a projective complex of $R$-modules.  
We assume that it satisfies  $r = \dim _k \Ext _R^1 (\Fx, \Fx) < \infty$. 
Then the following statements hold true. 

\begin{itemize}
\item[$(1)$]
There exists a universal lift $\Lx _0= ( F \widehat\otimes _k P_0 , \Delta _0)$  of   $\Fx$. 

\item[$(2)$]
A parameter algebra  $P_0$ is unique up to $k$-algebra isomorphisms.

\item[$(3)$]
Fixing a parameter algebra  $P_0$,  a universal lift $\Lx_0$ is unique up to chain isomorphisms of complexes of $R \widehat\otimes _k P_0 ^{op}$-modules. 

\item[$(4)$]
The parameter algebra has a description  $P_0 \cong T/I$, 
where  $T = k \langle\langle t_1, \ldots , t_r \rangle\rangle$ is a non-commutative formal power series ring of  $r$ variables  and  $I$ is a closed ideal which is contained in the square of the unique maximal ideal of  $T$. 
\end{itemize}
\end{thm}

\vspace{12pt}

We shall give a proof of this theorem in Section 3, where we need several new ideas to do so, because complete local $k$-algebras are not necessarily noetherian. 
We should remark that every complete local $k$-algebra can be a parameter algebra. 
In fact, for any complete local $k$-algebra $P$  with maximal ideal  $\m _P$,  
 $P$ itself is the parameter algebra for the universal lift of a free resolution of the left $P$-module $k = P /\m _P$. 
(See Theorem \ref{resol of k}). 

This theorem is essentially used in the proofs in Section 4,
 where we investigate the properties of parameter algebras by considering the  comparison of cohomology modules. 
As one of the main results there, 
we can give a certain structure theorem for parameter algebras. 
In fact, assuming  that  $r = \dim _k \Ext _R^1(\Fx, \Fx) < \infty$ and $\ell = \dim _k \Ext _R^2(\Fx, \Fx) < \infty$  for a projective complex  $\Fx$  of  $R$-modules, 
we have a description of the parameter algebra $P_0$  as 
$
P_0  \cong  
k \langle\langle t_1,\ldots t_r\rangle\rangle/ 
\overline{ ( f_1, \ldots , f_{\ell} )}.  
$
(See Theorem \ref{analytic gen}.)
In particular, if $\dim _k \Ext _R^2(\Fx, \Fx) =0$, then the parameter algebra equals a non-commutative formal power series ring. 

\vspace{12pt}

Let  $P_0$ be the parameter algebra of the universal lifts of  $\Fx$ which is 
 described as  $P_0 = T/I_0$, where $T$ is a non-commutative formal power series ring and  $I_0$  is a closed ideal of  $T$  with $I_0 \subseteq \m_T^2$. 
Then we prove in Theorem \ref{square iso} that there is an isomorphism of $k$-vector spaces  
$$
\Ext _R^1(\Fx, \Fx) ^2 \cong \Hom _k (I_0 /I_0 \cap \m_T^3, k), 
$$
where the left hand side means the $k$-subspace of   $\Ext _R^2(\Fx, \Fx)$  generated by all the products of two elements in  $\Ext _R^1(\Fx, \Fx)$. 
This isomorphism shows that 
$I_0 \subseteq \m_T^3$  if and only if  $\Ext _R^1(\Fx, \Fx) ^2 =0$. 
(See Corollary \ref{cor to square iso}.)

\vspace{12pt}

We can also regard such all observations as results of comparison of cohomology modules. 
For this, we assume that  $\Fx = (F, d)$  is a right bounded projective complex of  $R$-modules, and let $\Lx _0 = (F \widehat\otimes _k P_0,  \Delta _0)$ be the universal lift of $\Fx$.
For any integer  $n$, we have a projective complex of $R \otimes _k (P_0/\m _{P_0}^n)^{op}$-modules;
$$
\Lx _0 ^{(n)} = (F \otimes _k P_0/\m_{P_0}^n , \  \Delta _0 \otimes _{P_0}  P_0/\m_{P_0}^n), 
$$ 
which is a lift of  $\Fx$ to  $P_0/\m_{P_0}^n$. 
Therefore we have a morphism of Yoneda algebras as before; 
$$
\Ext _{P_0/\m_{P_0}^n}^{\cdot} (k , k) \to \Ext _R^{\cdot} (\Fx , \Fx).
$$
Taking the direct limit, we finally get the $k$-algebra homomorphism 
$$
\rho ^{\cdot} \ : \ \varinjlim \Ext _{P_0/\m_{P_0}^n}^{\cdot} (k , k) \to \Ext _R^{\cdot} (\Fx , \Fx).
$$
Our main problem is to see how the mapping $\rho ^i$  behaves for  $i \geq 0$. 
One can easily observe that $\rho ^0 \ : \ \varinjlim \Hom _{P_0/\m_{P_0}^n}(k , k) = k \to \End _R (\Fx)$  is a natural embedding and hence it is always an injection. 
Furthermore, by our construction of  $\Lx _0$ in Theorem \ref{summarize}, we see that 
$\rho ^1 \ : \ \varinjlim \Ext _{P_0/\m_{P_0}^n}^1 (k , k) = (\m_{P_0}/\m_{P_0}^2)^*  \to \Ext _R^1 (\Fx , \Fx)$  is a bijection

One of the main theorems of this paper is  Theorem \ref{main comparison}, in which we prove that 
$\rho ^2 \ : \ \varinjlim \Ext _{P_0/\m_{P_0}^n}^2 (k , k) \to \Ext _R^2 (\Fx , \Fx)$ is always an injection.
This actually realizes Idea \ref{1st conjecture}. 
We should notice that this holds because we had extended the notion of parameter algebras to non-commutative rings.

\newpage

\section{Non-commutative complete local algebras}


\subsection{Definitions and properties}
Throughout this paper, $k$  always denotes a field. 
Let  $A$  be an associative  $k$-algebra. 
By an ideal of  $A$  we always mean a two-sided ideal. 
When  $S$  is a subset of  $A$,  we denote by  $(S)$  the minimum ideal of $A$ that contains  $S$.

\begin{defn}\label{def complete local algebra}
Let  $A$  be an associative local $k$-algebra with Jacobson radical  $\m _A$. 
We say that  $A$  is a complete local $k$-algebra if the following three conditions are satisfied. 

\begin{itemize}
\item[$(a)$]
The natural inclusion  $k \subset A$  induces an isomorphism  $k \cong A / \m _A$. 
\item[$(b)$]
The $k$-vector space  $\m _A /\m _A^2$  is of finite dimension. 
\item[$(c)$]
$A$  is complete and separated in the $\m _A$-adic topology, i.e. 
the natural projections $A  \to  A/\m ^n \ (n \in \N)$  induce an isomorphism 
$ A \cong  \varprojlim  A / \m_A^n$.    
\end{itemize}

For  a complete local $k$-algebra  $A$,  we always denote by  $\m _A$  the Jacobson radical of  $A$, and we regard  $A$  as a topological ring with  $\m _A$-adic topology. 
\end{defn}

Note that any artinian local $k$-algebra  $A$  with  $A/\m _A \cong k$  is a complete local $k$-algebra in our sense.  

\begin{ex}
Let  $S = k \langle t_1, t_2 , \ldots , t_r \rangle$  be a free $k$-algebra over variables  $t_1, t_2 , \ldots , t_r$, and let  $J = ( t_1, t_2 , \ldots , t_r )$. 
We denote by  $T$  the $J$-adic completion of  $S$, i.e. 
$$
T  =  \varprojlim  S/ J^n,  
$$ 
and we call  $T$  the {\bf non-commutative formal power series ring}, which is denoted by  $\formal$. 
Clearly from the definition,  $T$  is a complete local $k$-algebra with maximal  ideal  $\m _T = (t_1, t_2, \ldots , t_r)$. 

Note that each element of  $T$  has a unique expression as a formal infinite sum $\sum_{\lambda} c_{\lambda} m _{\lambda}$, where  $c_{\lambda} \in k$ and the $m_{\lambda}$'s are distinct monomials on $t_1, t_2 , \ldots , t_r$. 
\end{ex}

\begin{rem}
Let  $f : A \to B$  be a $k$-algebra homomorphism of complete local $k$-algebras. Then it is easy to see that  $f$  is a local homomorphism, i.e.  $f(\m _A) \subseteq \m _B$. 
In particular, $f$  is a continuous map. 
\end{rem}

\begin{defn}
Let  $A$  be a complete local $k$-algebra and let  $I$  be an ideal (resp. a left or right ideal). 
Then we denote the closure of  $I$  by  $\overline{I}$, i.e.  $\overline{I} =  \bigcap _{n =0}^{\infty} \left( I + \m _A ^n \right)$. 
It is easy to see that  $\overline{I}$  is also an ideal (resp. a left or right ideal). 
We say that $I$ is a closed ideal (resp. a closed left or right ideal) if  $I = \overline {I}$. 
\end{defn}

\begin{rem}
If  $A$  is a commutative complete local $k$-algebra, then it is well-known that  $A$  is noetherian and every ideal of  $A$  is closed (cf. \cite{bourbaki}). 
But, in general, a non-commutative complete local $k$-algebra is not necessarily noetherian, and an ideal may not be closed.

For example, let  $T = k \langle\langle x, y \rangle\rangle$ and let  $I = (x)$. Since any element of  $I$ is a finite sum of elements of the form  $a x b$  with  $a, b \in T$, 
one can easily see that  $\sum _{n=1}^{\infty} y^nxy^n$  belongs to $\overline{I}$, but not to  $I$. 
\end{rem}

\begin{rem}
If  $I$  is a closed ideal of a complete local $k$-algebra. 
Then, $I$ is complete and separated in the relative topology  on  $I$, i.e. 
$$
I = \varprojlim  I/I \cap \m_A ^n. 
$$
\end{rem}

\begin{lemma}\label{quotient by closed ideal}
Let  $A$  be a complete local $k$-algebra and let  $I$ be an ideal of  $A$.
Then,  $A/I$  is a complete local $k$-algebra if and only if  $I$  is a closed ideal. 
\end{lemma}

\begin{pf}
Note that the residue ring  $A/I$ is complete (but may not be separated) in $\m _A$-adic topology. 
If $I$ is a closed ideal, then  $A/I$  is separated, hence  $A/I$  is a complete local $k$-algebra. 
Conversely, if  $A/I$  is a complete local $k$-algebra,  then the natural projection  $f : A \to A/I$  is continuous and  $\{0\} \subseteq A/I$ is closed. 
Therefore $I = f^{-1}(\{0\})$  is closed. 
\qed\end{pf}

\begin{lemma}\label{surjective}
Let  $f : A \to B$  be a $k$-algebra homomorphism of complete local $k$-algebras. 
Suppose that the induced mapping  $\overline{f} : \m_A / \m _A^2 \to \m _B/\m_B^2$  is surjective. 
Then  $f$  is a surjective homomorphism. 
\end{lemma}

\begin{pf}
It is easy to see by induction on  $n$  that the induced mappings 
$\overline{f} : \m_A^n /\m _A^{n+1} \to \m_B^n /\m _B^{n+1}$  are surjective for all $n \geq 1$.
Then, for a given  $b \in B$, we can find  $a_i \in \m _A^i \ (0 \leq i \leq n)$  such that  $f(a_0 + a_1 + \cdots + a_n ) - b \in \m _B^{n+1}$  for   $n \geq 0$. 
Thus, putting  $a = \sum _{n=0}^{\infty} a_n$, we have  $a \in A$  and  $f (a)=b$, since  $f$  is continuous. 
\qed\end{pf}

\begin{prop}\label{image of T}
Let  $A$  be a complete local $k$-algebra. 
Then, there are a non-commutative formal power series ring  $T = \formal$  and a $k$-algebra homomorphism  $f : T \to A$  such that the induced mapping  $\overline{f} :  \m_T /\m _T^{2} \to \m_A /\m _A^{2}$ is bijective.

In particular, $A$  can be described as  $A \cong T /I$,   where  $I$  is a closed ideal of  $T$  and  $I \subseteq  \m _T ^2$.
\end{prop}

\begin{pf}
Take  $x_1, x_2, \ldots , x_r \in \m _A \backslash \m_A^2$  which give rise to a basis of the $k$-vector space  $\m _A /\m_A^2$. 
Now define a $k$-algebra homomorphism $f : T  =  \formal \to A$  by  $f(t_i ) = x_i \ (1 \leq i \leq r)$. 
Then it is obvious that $f$ satisfies the desired conditions.
\qed\end{pf}

\begin{defn}
We denote by $\Ahat$  the category of  complete local $k$-algebras and $k$-algebra homomorphisms.
We also denote by $\A$  the category of  artinian local $k$-algebras $A$  with  $A /\m_A \cong k$  and $k$-algebra homomorphisms.
Obviously,  $\A$  is a full subcategory of  $\Ahat$.   
\end{defn}

\begin{rem}\label{projlim}
Let  $A$  be a complete local $k$-algebra. 
Then  $A/\m_A^n \in  \A$  for any $n \geq 1$ and by definition  $A = \varprojlim A/\m_A^n$. 
Conversely, let 
$$
\begin{CD}
\cdots @>>> A _{n+1} @>{f_{n+1}}>> A_n @>{f_n}>> A_{n-1} @>>> \cdots @>{f_2}>> A_1  
\end{CD}
$$
be a projective system in $\A$  such that  each $f_n$  induces an isomorphism 
$\m_{A_n}/\m_{A_n}^2 \cong \m_{A_{n-1}}/\m_{A_{n-1}}^2$. 
Then we have that  $\varprojlim A_n  \in \Ahat$. 

In fact,  we see from Lemma \ref{image of T} that each  $A_n$  is isomorphic to $T/I_n$  for any  $n\geq1$, where  $I_1 \supseteq I_2 \supseteq \cdots \supseteq I_n \supseteq I_{n+1} \supseteq \cdots$ are closed ideals of the non-commutative formal power series ring  $T$ . 
Then we have $\varprojlim A_n \cong T/\bigcap _{n=1}^{\infty}I_n$ and $\bigcap _{n=1}^{\infty}I_n$ is a closed ideal of  $T$. 
Thus the claim follows from Lemma \ref{quotient by closed ideal}.
\end{rem}

We remark here on the closedness of certain ideals in the non-commutative formal power series ring.
First we note the following lemma.

\begin{lemma}\label{f g ideals are free} 
Let   $I$  be a left ideal of  $T = \formal$, and suppose that $I$  is finitely generated as a left ideal. 
Then  $I$ is a free module as a left $T$-module. 
\end{lemma}

\begin{pf}
We note that  $I/\m_T I$  is  a finite dimensional $k$-vector space.
Hence we can take a finite number of elements  $f_1,  \ldots , f_n \in I$  which yield  a base of the $k$-vector space  $I/\m_T I$. 
First we claim that  $I$  is generated by  $f_1,  \ldots , f_n$ as a left ideal. 

To show this, let  $x$  be any element of  $I$. 
Since  $I = T\{ f_1,  \ldots , f_n \} + \m_T I$,  there are elements 
 $a_{01}, \ldots , a_{0n} \in T$  such that 
$x - \sum _{i=1}^n a_{0i}f_i \in \m_T I$. 
Then, apply the same argument to this element, we can find 
 $a_{11}, \ldots , a_{1n} \in \m_T$  such that 
$x - \sum _{i=1}^n a_{0i}f_i - \sum_{i=1}^n a_{1i}f_i \in \m_T ^2I$. 
Inductively, one can show that there are 
 $a_{\ell 1}, \ldots , a_{\ell n} \in \m_T^{\ell}$  with 
$x - \sum _{i=1}^n (a_{0i} + a_{1i}+ \cdots + a_{\ell i}) f_i \in \m_T ^{\ell +1} I$ for any  $\ell \geq 1$. 
Now put $\alpha _{i} = \sum _{\ell = 0} ^{\infty} a_{\ell i}$ which are  well-defined elements in  $T$,   and we have  $x = \sum _{i=1}^n \alpha _i f_i$. 
Thus the set   $\{ f_1,  \ldots , f_n \}$  generates  $I$  as a left ideal.

Now we prove that  $\{ f_1, \ldots , f_n \}$  is a free basis of  $I$ as a left $T$-module. 
To show this, let  $\sum _{i=1}^n a_i f_i =0$, where  $a_i \in T \ (1 \leq i \leq n)$. 
We have to show  $a_i =0$  for each  $i$. 

For this, we only have to prove,  by induction on  $\ell \geq 1$,  that  $a_i \ (1 \leq i \leq n)$  belong to $\m _T ^{\ell}$  for  all   $a _i \in T \ (1 \leq i \leq n)$  which satisfy the equality  $\sum _{i=1}^n a_i f_i =0$.  

Since   $\{ f_1, \ldots , f_n \}$  is a $k$-base of  $I/ \m_T I$ , 
it is trivial that  $a _i \in \m_T \ (1 \leq i \leq n)$. 
Hence the claim holds for $\ell = 1$. 
Now assume  $a_i \in \m _T ^{\ell} \ (1 \leq i \leq n)$ for  $\ell \geq 1$. 
Then we may write  $a_i = \sum _{j=1}^r t_jb_{ji}$  for some $b_{ji} \in \m_T ^{\ell-1}$. 
Thus we have 
$$
\sum_{j=1}^r  t_j \left( \sum _{i=1}^n b_{ji} f_i \right) = 0,  
$$
in  $T$. 
Since an element of  $T$  has a unique expression as a formal infinite sum of monomials with coefficients in $k$, 
it follows that  $\sum _{i=1}^n b_{ji}f_i = 0$  for any $j$. 
Then, by the induction hypothesis, we have  $b_{ji} \in \m_T ^{\ell}$, and hence   $a_i = \sum_{j=1}^r t_j b_{ji} \in \m _T ^{\ell +1}$  as desired.      
\qed
\end{pf}

The following lemma is known as Nagata's theorem for commutative formal power series ring, which is easily generalized to non-commutative ones. 

\begin{lemma}\label{nagata}
Let  $T = \formal$  be a non-commutative formal power series ring. 
Suppose a descending sequence $\ga _1 \supset \ga_2 \supset \ga _3 \supset \cdots$  of left ideals of  $T$  satisfies the equality $\bigcap _{i=1}^{\infty} \ga _i = (0)$. 
Then the linear topology on  $T$  defined by  $\{ \ga _i\ | \ i=1,2,\ldots \}$ is stronger than the $\m_T$-adic topology. 
\end{lemma}

\begin{pf}
The proof given in \cite[(30.1)]{nagata} is valid for non-commutative case. 
\qed\end{pf}

\begin{prop}\label{condition for closed ideal}
Let  $I$  be a left ideal of a complete local $k$-algebra  $A$. 
If one of the following conditions holds, then  $I$  is a closed left ideal in $A$.
\begin{itemize}
\item[$(a)$]
$A$ is a non-commutative formal power series ring  $T = \formal$ and $I$ is finitely generated as a left ideal. 
\item[$(b)$]
$I$  is of finite length as a left $A$-module, i.e. $\Kdim _k I < \infty$. 
\end{itemize}
\end{prop}

\begin{pf}
(a) 
By Lemma \ref{f g ideals are free} we may write  $I = T\{ f_1, \ldots , f_n\}= Tf_1 \oplus \cdots \oplus Tf_n$.  
We prove the lemm by induction on  $n$. 
If  $n =0$, then it is trivially true. 

Suppose  $n >0$  and set  $J = T\{ f_2, \ldots , f_n \}$. 
(We understand  $J =(0)$ if $n=1$.) 
Note that we have a direct decomposition  $I = Tf_1 \oplus J$  as a left $T$-module.
Now we set 
$$
\ga _{\ell} = \{ c \in T \ | \ cf_1 + g \in \m_T^{\ell} \ \ \text{for some} \ \ g \in J\},
$$  
for each  $\ell >0$. 
Note that  $\ga _{\ell}$  is a left ideal, $\ga _{\ell} \supseteq \ga _{\ell +1}$  and  $\m _T ^{\ell} \subseteq \ga _{\ell}$  for all  $\ell$. 
First of all, we claim that the following equality holds.  
$$
 \bigcap _{\ell =1}^{\infty} \ga _{\ell} = (0) \qquad\quad (*)
$$
In fact, for any element  $c \in  \bigcap _{\ell =1}^{\infty} \ga _{\ell}$, there is an element  $g_{\ell} \in J$ with  $cf_1 +g_{\ell} \in \m_T^{\ell}$  for each  $\ell$. 
Since  $g_{\ell} -g_{\ell +1} \in \m_T^{\ell}$,  we see that  $\{ g_{\ell} \}$  forms  a Cauchy sequence in the $\m _T$-adic topology. 
Therefore we see that  $cf_1 + \lim _{\ell \to \infty} g_{\ell} = 0$. 
Since  $J$  is a closed ideal by the induction hypothesis, we have  $\lim _{\ell \to \infty } g_{\ell} \in J$, and thus we have  $cf_1 \in J$. 
Then the direct decomposition  $I = Tf_1 \oplus J$  forces $cf_1 =0$, hence  $c =0$. 
This proves the equality  $(*)$. 
Note from Lemma \ref{nagata}  that the ideals  $\ga_{\ell}$  define the topology equivalent to the $\m_T$-adic topology. 

Now, to prove that  $I$  is closed, take an element  $x \in \overline{I}$. 
We want to show  $x \in I$. 
Take a sequence  $\{ a_{\ell} \ |\ \ell =1,2,\ldots \}$  in $I$  which converges to $x$  in  the $\m_T$-adic topology.
We may assume that  $a_{\ell} - a_{\ell +1} \in \m _T ^{\ell}$  for each $\ell$.  Each $a_{\ell}$  has a unique description $a_{\ell} = \sum _{i=1}^{n} b_{\ell, i}f_i$  for some  $b_{\ell, i}\in T$. 
Thus  $\sum _{i=1}^{n} (b_{\ell, i} - b_{\ell +1, i} ) f_i \in \m_T ^{\ell}$. 
Therefore  $b_{\ell, 1} - b_{\ell+1, 1} \in \ga _{\ell}$ for any $\ell$. 
Then, by the fact we have shown above, we see that  $\{ b_{\ell, 1}\ | \ i = 1,2, \ldots \}$  is a Cauchy sequence in the $\m_T$-adic topology. 
This is true for the sequences   $\{ b_{\ell, i} \ | \ \ell = 1,2,\ldots \}$  
 for all  $i \ (1 \leq i \leq n)$. 
Since  $T$  is complete in the $\m_T$-adic topology, 
 the sequence  $\{ b_{\ell, i}\ | \ \ell \geq 1\}$  converges to an element  $c_i \in T$  for each  $i$. 
Then, 
$x = \lim _{\ell \to \infty} a_{\ell} 
= \sum _{i=1} ^n  \lim _{\ell \to \infty} b_{\ell i} f_i 
=\sum _{i=1} ^n  c_i f_i \in I$  as desired. 

\vspace{6pt}

(b)
Since  $\bigcap _{n=1}^{\infty} I \cap \m_A^n \subseteq \bigcap _{n=1}^{\infty} \m_A^n = (0)$, and since  $\Kdim _k I < \infty$, there is an integer  $n_0$ such that  $I \cap \m_A ^n =(0)$  for  $n \geq n_0$.
Thus  $I + \m_A^n = I \oplus \m_A^n $  for  $n \geq n_0$. 
Therefore, 
$$
\overline{I } = \bigcap _{n=n_0}^{\infty} I + \m_A^n  
 = \bigcap _{n=n_0}^{\infty} I \oplus \m_A^n      
 = I  \oplus \bigcap _{n=n_0}^{\infty}  \m_A^n      
= I.  \qquad \text{\qed}
$$
\end{pf}

\begin{cor}\label{condition for Ahat}
Let  $A \in \Ahat$  and let  $I$  be an ideal of  $A$. 
Suppose one of the conditions in the previous proposition holds. 
Then we have   $A/I \in \Ahat$.  
\end{cor}

The Artin-Rees lemma  for non-commutative formal power series ring  holds in the following form.

\begin{cor}\label{relative top}
Let  $I $  be a finitely generated left ideal of the non-commutative formal power series ring  $T = \formal$. 
Then, the relative topology on  $I$  induced from  $T$  is equivalent to the $\m_T$-adic topology on $I$. 
That is, for any  $m \geq 1$, there is an integer  $\ell \geq 1$  such that 
$\m_T ^{\ell } \cap I \subseteq \m_T ^{m} I$.  
\end{cor}

\begin{pf}
In the proof of (a) in Proposition \ref{condition for closed ideal}, 
we have shown that, for a given $m \geq 1$, there is an integer  $\ell _1 \geq 1$  such that  $\ga _{\ell _1} \subseteq \m _T ^{m}$. 
This shows that  $c_1 f_1 + \cdots + c_n f_n \in \m_T^{\ell_1}$  implies  $c_1 \in \m _T ^{m}$. 
This is true for any  $i \ (1 \leq i \leq n)$, that is,  there is an integer  $\ell _i >0$  such that  $c_1 f_1 + \cdots + c_n f_n \in \m_T^{\ell_i}$  implies  $c_i \in \m _T ^{m}$. 
Now  take $\ell$  so that  $\ell > \ell _i$  for all $i\ (1 \leq i \leq n)$. 
Then we have that  $c_1 f_1 + \cdots + c_n f_n \in \m_T^{\ell}$  implies  $c_i \in \m _T ^{m}$  for all $i$. 
Hence  $I \cap \m_T ^{\ell} \subseteq  \m^{m}I$.  
\qed\end{pf}

\begin{defn}
Let  $A$  be a complete local $k$-algebra and let  $S$ be a subset of  $A$. 
Then we say that an ideal  $I$  {\bf is analytically generated by}  $S$  if  $I = \overline{ (S)}$. 
\end{defn}

\begin{prop}\label{nakayama}
Let  $A$  be a complete local $k$-algebra and let  $S$  be a subset of a closed  ideal $I$ of  $A$. 
Then,   $I$ is analytically generated by  $S$   if the image of  $S$  generates  $I / \overline{\m_A I + I \m_A}$ as a $k$-vector space. 

Furthermore, if  $S$  is a finite subset, then the converse is also true. 
\end{prop}

\begin{pf}
Suppose that the image of  $S$  generates  $I / \overline{\m_A I + I \m_A}$, and let $n$ be an arbitrary  natural number. 
Since  $\overline{\m_A I+I\m_A} \subseteq  \m_A I+I\m_A + (\m_A^n \cap I) \subseteq I$, 
the set  $S$  generates  $I/\m_A I+I\m_A + (\m_A ^n \cap I)$  as  a $k$-vector space, hence 
$$
\begin{aligned}
I &= (S) + \m_A I + I \m_A + (\m_A^n \cap I) \\
&= (S) + \m_A ((S)+ \m_A I + I \m_A ) + ((S) +  \m_A I + I \m_A) \m_A + (\m_A ^n \cap I) \\
&= (S) + (\m_A ^2 I + \m_A I \m_A + I \m_A ^2) + (\m_A^n \cap I) \\ 
& \cdots \\
&= (S) + \sum_{i+j=s} \m_A ^i I \m_A ^j + (\m_A ^n \cap I), \\  
\end{aligned}
$$
for $1 \leq s \leq n$. 
Finally, putting  $s = n$,  we have that 
$I = (S) + (\m _A^n \cap I)$. 
Since this equality holds for all  $n \geq 1$, we have  $I = \overline{(S)}$. 
  
\vspace{6pt}
To prove the converse, we assume that  $I$  is analytically generated by a  finite subset $S$.
Then the equality  $I = \bigcap _{n=1}^{\infty} ((S) + \m_A^n )$  holds. 
Since  $I \subseteq (S) + \m_A^n$, we have  $I = (S) + (\m_A^n \cap I)$ for all  $n \geq 1$. 
Thus the image of  $S$  generates  $I/ \m_A I + I\m_A + (\m_A ^n \cap I)$ as a $k$-vector space,  for  all  $n \geq 1$. 
In particular, $\dim _k I/\m_A I + I \m_A + (\m _A^n \cap I) \leq |S|$. 
Since  $|S|$  is finite,  there is an integer  $n_0 >0$  such that  
$I/\m_A I + I \m_A + (\m _A^n \cap I) = I/\m_A I + I \m_A + (\m _A^{n+1} \cap I)$ for all  $n \geq n_0$. 
Thus, we have the equality  $I/\overline{\m_A I + I \m_A} = I/\m_A I + I \m_A + (\m _A^{n_0} \cap I)$, which is generated by $S$  as a $k$-vector space.
\qed\end{pf}

\begin{cor}\label{cor nakayama}
Let  $I$  be a closed ideal in a complete local $k$-algebra $A$. 
Then, the equality  $I = \overline{ \m_A I + I \m_A}$ implies  $I =(0)$ 
\end{cor}

\begin{cor}
Let  $I$  be a closed ideal in a complete local $k$-algebra $A$. 
Then, $I$ is analytically generated by a finite number of elements of $I$ if and only if  $\Kdim _k \left( I / \overline{\m_A I + I \m_A} \right)< \infty$.
\end{cor}

\begin{cor}\label{artin-rees}
Let  $I$  be a closed ideal in a complete local $k$-algebra $A$ that is analytically generated by a finite number of elements. 
Then, the equality
$$
\overline{\m_A I + I \m_A} = \m_A I + I \m_A + (\m_A^n \cap I) 
$$ 
holds for any large integer  $n \gg 1$. 
\end{cor}

\begin{pf}
See the proof of Proposition \ref{nakayama}.
\qed\end{pf}

It is well-known that the category  $\A$  admits the fiber products.

\begin{lemma}\label{fiber}
The category  $\Ahat$  admits the fiber products, that is, any diagram in  $\Ahat$  
$$
\begin{CD}
  @.   B \\
@. @VgVV \\
A @>f>> C  \\
\end{CD}
$$
can be embedded into a pull-back diagram  
$$
\begin{CD}
Q  @>>>  B \\
@VVV @VgVV \\
A @>f>> C.  \\
\end{CD}
$$
\end{lemma}

\begin{pf}
For any integer $n$, we have a diagram in  $\A$ 
$$
\begin{CD}
  @.   B/\m _B^n  \\
@. @Vg_nVV \\
A/ \m_A^n  @>f_n>> C/\m _C ^n   \\
\end{CD}, 
$$
from which we  have a fiber product  $Q_n := A/ \m_A^n  \times _{C/\m _C^n}  B/\m _B^n$  in the category  $\A$. 
It is clear that  $\{ Q_n | \ n \geq 1\}$ forms a projective system in  $\A$. 
Put  $Q = \varprojlim \ Q_n$, and we have  $Q \in  \Ahat$  by Remark \ref{projlim}. 
It is  routine to show that  $Q$  is a fiber product in  $\Ahat$. 
\qed
\end{pf}

\begin{rem}
In the setting of Lemma \ref{fiber}, the fiber product $Q$  and its Jacobson radical  $\m _Q$  can be described in the following way : 
$$
Q = \{ (a, b) \in A \times B \ | \ f(a) = g(b) \}, \quad 
\m_Q = \{ (a, b) \in A \times B \ |\ f(a) = f(b) \in \m _C \}.  
$$
We denote the fiber product $Q$  by  $A \times_C B$.  
\end{rem}


\subsection{Small extensions}

Let  $A$  be a complete local $k$-algebra. 
We say that an element  $\epsilon \not= 0$ in  $A$  is a socle element of  $A$  if  $\m _A \epsilon = \epsilon \m_A = 0$. 
Note that an element  $\epsilon$ of  $A$  is a socle element if and only if the ideal  $(\epsilon)$  is a one-dimensional $k$-vector space. 
Note that if  $A$  is an artinian local $k$-algebra then there exists at least one socle element. 

One should remark from Corollary \ref{condition for Ahat} that, 
if   $\epsilon$  is a socle element in a complete local $k$-algebra $A$, 
then  $\overline{A} = A /(\epsilon)$  is also a complete local $k$-algebra.

\begin{defn}
A pair  $(A' , \epsilon)$  is called a small extension of a complete local $k$-algebra  $A$  if $\epsilon$  is a socle element of a complete local $k$-algebra  $A'$  and  $A'/(\epsilon) \cong A$  as a $k$-algebra. 
To describe the small extension  $(A', \epsilon)$ of $A$, we often write it as a short exact sequence  
$$
\begin{CD} 
0 @>>> k @>\epsilon>> A' @>\pi>> A @>>> 0, 
\end{CD}
$$ 
where  $\pi$  is the natural projection.
\end{defn}

\begin{lemma}\label{trivial or nontrivial small extensions}
Let  $(A' , \epsilon)$  be a small extension of a complete local $k$-algebra $A$. 

\begin{itemize}
\item[$(a)$]
If  $\epsilon \not\in \m _{A'}^2$, then there is a $k$-algebra homomorphism $\iota : A \to A'$  that is a right inverse of  $\pi : A' \to A$. 
In this case, $A'$  is isomorphic to  $A[x] / (x^2, \m_A x, x \m_A)$  as a $k$-algebra, which we call a {\bf trivial small extension} of  $A$.

\item[$(b)$]
If  $\epsilon \in \m_{A'}^2$, and if  $A = T/I$  where  $I$  is a closed ideal of  $T = k \langle\langle t_1, t_2. \ldots , t_r\rangle\rangle$  and  $I \subseteq  \m_T^2$,  then there is a closed ideal  $J \subseteq I$  of  $T$  such that 
$A' \cong T/J$  and  the length  $\ell _T (I/J) = 1$. 
In this case, we say that  $(A', \epsilon)$  is a {\bf nontrivial small extension}. 
\end{itemize}
\end{lemma}

\begin{pf}
(a)  Suppose  $\epsilon \not\in \m _{A'}^2$. 
Then, since  $(\epsilon) \cong k$, we have  $(\epsilon ) \cap \m_{A'}^2 = (0)$. Thus we can take a $k$-subspace  $\n$ of  $\m_{A'}$  such that $\m_{A'}^2 \subseteq \n$ and $\m_{A'} = (\epsilon) \oplus \n$ as a $k$-vector space. 
Noting that  $\n ^2 = \m_{A'}^2$, we see that the $k$-subspace  $k \oplus \n \subseteq A'$  is actually a  $k$-subalgebra and the restriction to $k \oplus \n$ 
of $\pi : A' \to A$  yields an isomorphism 
$k \oplus \n \cong  A$. 

(b)   Suppose  $\epsilon \in \m _{A'}^2$. 
Then we have  $\m_{A'}/\m_{A'}^2 \cong \m_{A}/\m_{A}^2$. 
It follows from Lemma \ref{image of T} that 
there is a commutative diagram in $\Ahat$  
$$
\begin{CD}
T @= T \\
@Vf'VV @VfVV \\
A' @>>> A,  
\end{CD} 
$$
where  $f$  and  $f'$  are surjective and  $I = \kernel (f)$. 
It is easy to see that  $J = \kernel (f')$  satisfies the desired conditions.
\qed\end{pf}

\begin{defn}
Let  $A \in \Ahat$. 
For small extensions  $(A_1 , \epsilon _1)$  and  $(A_2, \epsilon _2)$  of  $A$,  we say that   $(A_1 , \epsilon _1)$  and  $(A _2, \epsilon _2)$ are equivalent, denoted by  $(A_1 , \epsilon _1) \sim (A _2, \epsilon _2)$, if there is a $k$-algebra isomorphism $f : A _1 \to A_2$  with  $f (\epsilon _1) = \epsilon _2$.
We denote by  $\T (A) $  the set of equivalence classes of small extensions of  $A$ : 
$$
\T (A) = \{ (A', \epsilon ) \ | \ \epsilon \in A' \in \Ahat, \ (\epsilon ) \cong k, \  A'/(\epsilon) \cong A \} / \sim.
$$  
For a small extension  $(A', \epsilon )$  we denote its equivalence class by 
$[A', \epsilon]$. 
\end{defn}

Note from Lemma \ref{trivial or nontrivial small extensions} that trivial small extensions defines a unique element of  $\T (A)$.

\begin{lemma}\label{abelian group}
Let $A \in \Ahat$. 
Then  $\T(A)$  is an abelian group in which 
the zero element is the class of a trivial small extension. 
\end{lemma}

\begin{pf}
Let  $[A_1, \epsilon _1]$  and  $[A_2, \epsilon _2]$  be elements in  $\T(A)$. 
Then we have the following commutative diagram by taking the fiber product.
$$
\begin{CD}
   @. @. 0 @. 0 \\
@. @. @VVV  @VVV \\
  @.   @.  k @= k \\
  @.   @.  @V(0, \epsilon_2)VV  @V\epsilon_2VV  \\
0 @>>> k @>(\epsilon_1, 0)>> A_1 \times _{A} A_2 @>>> A_2 @>>> 0 \\
@. @|  @VVV @VVV \\
0 @>>> k @>{\epsilon_1}>> A_1 @>>> A @>>> 0 \\
@. @. @VVV  @VVV \\
@. @. 0 @. 0 \\
\end{CD}
$$
Put  $B = A_1 \times _{A} A_2 /(\epsilon_1, - \epsilon _2)$ and it follows from the exact sequence of the middle row in the diagram that there is an exact sequence 
$$
\begin{CD}
0 @>>> k @>{(\epsilon _1, 0)}>> B @>>> A @>>> 0. 
\end{CD}
$$
Note that, since  $A_1 \times _{A} A_2$ is a complete local $k$-algebra by  Lemma \ref{fiber} and  $(\epsilon_1, - \epsilon _2)$ is its socle element, it follows from Corollary  \ref{condition for Ahat} that $B$ is a complete local $k$-algebra. 
Hence  $(B, (\epsilon _1, 0))$  is a small extension of  $A$. 
Note that  $(\epsilon _1, 0) = (0, \epsilon _2)$  in  $B$. 
Now we define the sum by 
$$
[A_1, \epsilon_1] + [A_2,  \epsilon _2] = [B, (\epsilon _1, 0)].
$$ 
Then it is routine to verify that  $\T(A)$  is an abelian group by this definition of addition.
Actually, the commutativity of sum is given by the isomorphism 
$$
A_1 \times _{A} A_2 / (\epsilon _1, -\epsilon _2) \cong A_2 \times _{A} A_1 / (\epsilon _2, -\epsilon _1), \quad   (\epsilon _1 , 0) \leftrightarrow (\epsilon _2, 0).  
$$
The associativity is induced by 
$$
\{A_1 \times _{A} A_2 / (\epsilon _1, -\epsilon _2)\}\times_{A}A_3 /((\epsilon _1, 0), -\epsilon _3) 
\cong A_1 \times _{A} \{A_2\times _{A} A_3 / (\epsilon _2, -\epsilon _3)\} / 
(\epsilon _1, -(\epsilon _2, 0)).
$$
Let   $(A_0, \epsilon _0)$  be a trivial small extension of  $A$. 
Then we can show  $A_1 \times _A A_0 /(\epsilon _1, -\epsilon _0) \cong A_1$, which implies that $[A_0 , \epsilon _0]$  is the zero element in  $\T (A)$. 
Note that  the inverse element is given in the following. 
$$
-[A_1, \epsilon _1] = [A_1, - \epsilon _1]
$$ 
In fact, using the following lemma \ref{dual number},  one can show the isomorphism 
$$
A_1 \times _A A_1 /(\epsilon _1 , \epsilon _1 ) 
\cong 
A_1 \times _k D / (\epsilon _1, 0) \cong 
A_0.  \qquad  \text{\qed}
$$
\end{pf}

\begin{lemma}\label{dual number}
Let  $A \in \Ahat$  and let  $D = k[\epsilon _0]/(\epsilon _0^2)$, which we call the ring of  dual numbers over $k$. 
Then we have the following isomorphism of complete local $k$-algebras for any  $(A_1, \epsilon _1) \in \T (A)$. 
$$
A_1 \times_A A_1 \cong  A _1 \times _k D
$$
\end{lemma}

\begin{pf}
Define  $f :  A _1 \times _k D \to A_1 \times_A A_1$  by 
$f ((a_1, \overline{a_1} + c \epsilon _0)) = (a_1, a_1 + c \epsilon _1)$, 
where  $a_1 \in A_1$  and  $c \in k$, and  $\overline{a_1} \in k = A_1 / \m _{A_1}$  is the natural image of  $a_1 \in  A_1$. 
Then it is easy to see that $f$  is an isomorphism of $k$-algebras.   
\qed
\end{pf}

Let   $[A_1, \epsilon _1] \in \T(A)$ for  $A \in \Ahat$. 
We define the scalar product by an element  $c \in k$  as follows : 
$$
c \cdot [A_1, \epsilon ] = 
\begin{cases}
[A_1, c^{-1} \epsilon ]  & ( c \not= 0) \\
\text{the class of a trivial small extension}  & (c=0). \\ 
\end{cases}
$$

\begin{lemma}
Let  $A \in \Ahat$. 
Then  $\T (A)$  is a $k$-vector space by the above action of $k$. 
\end{lemma}

\begin{pf}
Let $c_1, c_2 \in k$  and  $[A_1, \epsilon _1], [A_2, \epsilon _2] \in \T (A)$. It is obvious from the definition that  $(c_1c_2)\cdot [A_1, \epsilon _1] = c_1(c_2\cdot [A_1, \epsilon _1])$. 
When $c_1 \not= 0$, the identity  $c_1\cdot ([A_1, \epsilon _1]+[A_2, \epsilon _2]) = c_1\cdot [A_1, \epsilon _1]+c_1 \cdot [A_2, \epsilon _2]$  follows from the isomorphism 
$A_1 \times _A A_2/(\epsilon _1, -\epsilon _2) \cong 
A_1 \times _A A_2/(c_1 ^{-1} \epsilon _1, -c_1^{-1} \epsilon _2)$.
We have to verify the equality  $(c_1+c_2)\cdot [A_1, \epsilon _1] = c_1 \cdot [A_1, \epsilon _1] + c_2 \cdot [A_1, \epsilon _1]$. 
If  one of  $c_1, c_2$ and $c_1+c_2$ is equal to zero, then it is  easy to see the equality holds. 
We assume that  $c_1 \not=0, c_2 \not=0$ and  $c_1+c_2\not=0$. 
In this case, we have from Lemma \ref{dual number} the isomorphism 
$$
A_1 \times_A  A_1 /(c_1 ^{-1} \epsilon _1, \  -c_2^{-1} \epsilon _1)
\cong A_1 \times _k D /(c_1 ^{-1} \epsilon _1,  -(c_1^{-1}+c_2^{-1}) \epsilon _0) \cong A_1, 
$$
and by this isomorphism $(c_1 ^{-1} \epsilon _1, 0)$ corresponds to 
$(c_1+c_2) ^{-1} \epsilon _1 \in A_1$. 
Hence, $[ A_1 \times_A  A_1 /(c_1 ^{-1} \epsilon _1, \  -c_2^{-1} \epsilon _1)
, \ (c_1 ^{-1} \epsilon _1, 0)] = [A_1, \ (c_1+c_2) ^{-1} \epsilon _1 ]$. 
\qed\end{pf}

\begin{lemma}\label{functoriality of T}
Let  $A_1, A_2 \in \Ahat$  and let  $f : A_1 \to A_2$  be a $k$-algebra homomorphism. 
Then  $f$  induces a $k$-linear map  $f^* : \T (A_2) \to \T(A_1)$. 
Therefore, $\T$  is a contravariant functor from  $\Ahat$  to the category of $k$-vector spaces.
\end{lemma}

\begin{pf}
For a given  $[A_2', \epsilon '_2] \in \T(A_2)$, take a fiber product 
$$
\begin{CD}
0 @>>> k @>{\epsilon' _1}>> A_1' @>>> A_1 @>>> 0 \\
@.  @| @VVV @VfVV  @. \\ 
0 @>>> k @>{\epsilon' _2}>> A_2' @>>> A_2 @>>> 0, \\
\end{CD}
$$
and we get a small extension  $(A'_1, \epsilon' _1)$ of  $A_1$. 
Now define  $f^* ([A'_2, \epsilon' _2]) = [A'_1, \epsilon' _1]$. 
It is not difficult to verify that $f^*$  is a $k$-linear mapping.
\qed\end{pf}

\begin{defn}\label{def of con}
Let  $I$  be a closed ideal of  $T = k \langle \langle t_1, t_2, \ldots, t_r \rangle \rangle$. 
We always regard $I$ as a topological $T$-bimodule by the relative topology induced from  $T$. 
Therefore, the set  $\{ I \cap \m_T^n \ | \ n = 1,2,\ldots \}$  gives the fundamental open neighbourhoods of  $0$  in $I$.  
We also consider the unique simple $T$-bimodule $k$ with discrete topology. 
We set 
$$
\Hom_{con}(I, k) = \{ f : I \to k \ | \ f \ \text{is a continuous } T\text{-bimodule homomorphism}\}. 
$$
Note that  $f \in \Hom _{T\text{-bimod}}(I, k)$  belongs to  $\Hom_{con}(I, k)$ if and only if  $f (I \cap \m_T^n) = 0$ for a large integer  $n$. 
It is clear that  $\Hom _{con}(I, k)$  is naturally a $k$-vector space.

Since  $f(\m_T I + I \m_T) = 0$  for  $f \in \Hom _{con}(I, k)$,  such an  $f$ induces the continuous map  $\overline{f} : I/\overline{\m_T I + I \m_T} \to k$. 
Hence,   
$$
\Hom _{con} (I, k) \cong \{ \overline{f} : I/\overline{\m_T I + I \m_T} \to k \ | \ \overline{f} \ \text{is a continuous } \ k\text{-linear map}\}.
$$  
Note, however,  that the induced topology on $I/\overline{\m_T I + I \m_T}$ may not be discrete. 

Let  $A \in \Ahat$, which we describe  as  $A = T/I$  where  $T$ is a non-commutative formal power series ring and  $I \subseteq  \m_T^2$. 
Under such a circumstance, we define the mapping 
$$
\tau : \Hom _{con} (I, k)  \to  \T(A), 
$$
as follows: 
For $f \in \Hom _{con} (I, k)$, if  $f = 0$, define  $\tau (f)$  to be the class of a trivial small extension. 
If  $f\not=0$, then $I_f := \kernel (f)=f^{-1}(0)$  is a closed ideal of  $T$ and hence  $A_f := T/I_f$  is a complete local $k$-algebra and we can take a unique element $\epsilon _f \in I/I_f \subseteq A_f$  with  $f(\epsilon _f) =1$. 
Since  $I = I_{f} + (\epsilon _f)$, $(A_f, \epsilon _f)$  is a small extension of  $A$. 
We define   
$\tau (f) = [ A_f , \epsilon _f]$. 
\end{defn}

\begin{prop}\label{con iso T}
The mapping  $\tau : \Hom _{con} (I, k) \to \T(A)$ is an isomorphism of $k$-vector spaces. 
\end{prop}

\begin{pf}
First we show that $\tau$  is a $k$-linear mapping. 
To show that  $\tau (cf) = c\cdot \tau(f)$  for  $c \in k$ and $f \in \Hom _{con}(I, k)$, we may assume that  $c \not=0$. 
Then it is trivial that  $I_f = I_{cf}$, hence  $A_f = A_{cf}$ and  $\epsilon _{cf} = c^{-1} \epsilon _f$. 
Thus it follows that $\tau (cf) = c \cdot \tau (f)$.

To show  $\tau (f + g) = \tau(f) + \tau (g)$  for $f, g \in \Hom _{con}(I, k)$,  we assume that  $f\not=0, g\not=0$ and  $f+g \not=0$. 
(Otherwise, the equality is proved easily.)
Suppose  $f$  and  $g$  are linearly dependent over $k$, hence  $f = cg$  for some $c \ (\not= 0, -1) \in k$. 
In this case, we have  $I_f = I_g = I_{f+g}$.  
 Since  $(f+g) (\epsilon _g ) = c+1$, we see   $\epsilon _{f+g} = (c+1)^{-1}\epsilon_g$. 
Therefore, 
$\tau (f+g) = [A_{f+g}, \epsilon _{f+g}] = 
[A_{g}, (c+1)^{-1}\epsilon _{g}]=(c+1)\cdot[A_g, \epsilon_g] =
c \cdot [A_g, \epsilon _g] + [A_g, \epsilon _g] = c\cdot \tau (g) + \tau (g)
= \tau (f) + \tau (g)$. 

Now suppose  $f$  and  $g$ are linearly independent over $k$. 
In this case  $I_f \not= I_g$, and hence  $I_f + I_g =I$. 
It then follows from the obvious exact sequence
$$
\begin{CD}
0 @>>> T/I_f \cap I_g @>\phi>> A_f \times A_g @>>> A @>>> 0
\end{CD}
$$
that we can take  $e_f, e_g \in I$  whose images in $T/ I_f \cap I_g$ are mapped respectively to  $(\epsilon _f,0), (0, \epsilon _g)$  by $\phi$. 
Note that  $I_f \cap I_g \subset I_{f+g}$. 
And note also that  $f(e_f)=1$, $f(e_g)=0$, $g(e_f)=0$ and  $g(e_g)=1$, hence  $(f+g)(e_f-e_g) =0$  and  $(f+g)(e_f)=1$. 
It is then easy to see that  $I_{f+g} = (I_f \cap I_g) + (e_f-e_g)$. 
Since  $A_f \times_A  A_g \cong T/I_f \cap I_g$, we have from the definition that $\tau(f) + \tau (g) = [A_f , \epsilon _f]+[A_g , \epsilon _g] =[T/I_{f+g}, e_f] = \tau (f+g)$.

Now we have proved that  $\tau$ is $k$-linear. 
Assume  $f \not= 0$. 
Then, since  $\epsilon _f \in I/I_f$ and  $I \subseteq \m_T^2$, we have  $\epsilon _f \in \m_{A_f}^2$. 
This implies  $\tau (f) \not= 0$  by  Lemma \ref{trivial or nontrivial small extensions}.   
Thus  $\tau$  is injective. 
The surjectivity of  $\tau$  is obvious from the definition of  $\tau$ and  $\T(A)$. 
\qed\end{pf}

\subsection{Complete tensor products}

In this section, let  $R$  be an associative algebra over a field $k$.  

\begin{defn}
For a complete local $k$-algebra  $A \in \Ahat$, 
we define the complete tensor product  $R \ctensor A$  as follows :
$$
R \ctensor A = \ \varprojlim R \otimes _k A/\m_A^n.  
$$
\end{defn}

Note that   $R \ctensor A$  is an associative $k$-algebra, since each mapping   $R \otimes _k A/\m _A ^{n+1} \to R \otimes _k A/ \m_A^n$  is a $k$-algebra homomorphism for  $n \geq 1$. 
Also note that, if  $A \in \A$,  then  $R \ctensor A = R \otimes _k A$  is an ordinary tensor product of $k$-algebras.

\begin{rem}
In general, $R \otimes _k A$  is a subalgebra of  $R \ctensor A$. 
However, they are distinct  in general. 

For example, let  $R = k[x]$  and $T = k \langle\langle t \rangle\rangle = k[[t]]$ (with one variable).  
Then, we have 
$
R \otimes _k T = k[[t]] \ [x] \subset R \ctensor T = k[x]\ [[t]],       
$
which are actually distinct. 
\end{rem}

\begin{defn}
Let  $M$  be a left $R$-module and let  $X$  be a right (resp. left)  $A$-module,  where  $A \in \Ahat$. 
Then note that,   for each  $n \geq 1$,   
$M \otimes _k X/ X \m_A ^n$ (resp.  $M \otimes _k X/ \m_A ^nX$)   is a left $R \otimes _k (A/\m_A^n)^{op}$-module (resp.  a left $R \otimes _k (A/\m_A^n)$-module), i.e. a left module over  $R$  and a right (resp. left) module over  $A/\m_A^n$.    
We define the complete tensor product by 
$$
M \ctensor X = \varprojlim M \otimes _k X/X \m_A ^n 
\quad (\text{resp.} \ \  M \ctensor X = \varprojlim M \otimes _k X/ \m_A ^nX), 
$$
which is a left  $R \ctensor A^{op}$-module (resp. a left $R \ctensor A$-module) by the reason above. 
\end{defn}

We always consider  $M \otimes _k X$  and  $M \ctensor X$  with  $\m_A$-adic topology. 
In general,  there is a natural mapping  $M \otimes _k A \to M \ctensor A$, which is the completion map in $\m_A$-adic topology.

\begin{rem}
(a)  If  $M$  is of finite dimension as a $k$-vector space with a $k$-basis  $\{ e_1 , \ldots, e_{\ell}\}$, then we have
$$
M \ctensor A  = \varprojlim \left( \bigoplus _{i=1}^{\ell} e_i k \otimes _k A/\m_A ^n \right)  = \bigoplus _{i=1}^n e_i \left( \varprojlim A/\m_A ^n \right) = \bigoplus _{i=1}^n e_i A
$$
for any  $A \in \Ahat$.
Thus  $M \ctensor A$  is a free module as a right $A$-module if  $\Kdim _k M < \infty$.  

(b) 
Suppose  $M$  is of infinite dimension as a $k$-vector space with basis $\{ e_{\lambda} \ | \ \lambda \in \Lambda \}$. 
In this case, we have 
$$
M \ctensor A = \varprojlim \left[ \bigoplus _{\lambda} e_{\lambda}(A/\m_A ^n) \right]. 
$$
Therefore, an element of  $M \ctensor A$  is described to be a formal sum 
$\sum _{\lambda} e_{\lambda} x_{\lambda}\ (x_{\lambda} \in A)$  as an element of  $\prod _{\lambda} e_{\lambda} A$. 
Note that  
$\sum _{\lambda} e_{\lambda} x_{\lambda}  \in  \prod _{\lambda} e_{\lambda} A$ 
 belongs to  $M \ctensor A$  if and only if 
$$
\sharp \{ \lambda \in \Lambda \ | \ x_{\lambda} \not\in \m_A ^n \} < \infty 
$$
for all  $n \geq 1$. 
\end{rem}

\begin{lemma}\label{exact seq}
Let  $A \in \Ahat$. 
\begin{itemize}
\item[$(a)$]
Let   $X$  be a right $A$-module, and let 
$$
\begin{CD}
0 @>>> L @>>> M @>>> N @>>> 0
\end{CD}
$$
be a short exact sequence of left $R$-modules. 
Then the complete tensor product by  $X$  induces the exact sequence of left 
 $R \ctensor A^{op}$-modules 
$$
\begin{CD}
0 @>>> L \ctensor X @>>> M  \ctensor X @>>> N  \ctensor X @>>> 0.
\end{CD}
$$
\item[$(b)$]
Let  $M$  be a left  $R$-module, and let 
$$
\begin{CD}
0 @>>> X @>>> Y @>>> Z @>>> 0, 
\end{CD}
$$
be a short exact sequence of left $A$-modules. 
Then we have an exact sequence of left  $R \ctensor A$-modules
$$
\begin{CD}
0 @>>> \varprojlim \left( M \otimes _k  (X/X\cap \m_A^nY)\right)  @>>> M  \ctensor Y @>>> M  \ctensor Z @>>> 0.
\end{CD}
$$
In particular, if the relative topology on  $X$  induced from  the $\m_A$-adic topology on  $Y$ is equivalent to the $\m_A$-adic topology on  $X$, then we have  an exact sequence of  $R \ctensor A$-modules
$$
\begin{CD}
0 @>>> M \ctensor X @>>> M  \ctensor Y @>>> M  \ctensor Z @>>> 0.
\end{CD}
$$
\end{itemize}
\end{lemma}

\begin{pf}
The proof is similar to the commutative complete case in \cite{bourbaki} \cite{matsumura} or \cite{nagata}. 
\qed\end{pf}

Let  $M$  be a left  $R$-module and let  $X$  be a left $A$-module where  $A \in \Ahat$. 
Then, from the definition of complete tensor products, we see that there is a natural mapping
$$
\gamma _{M,X} \ : \ (M \ctensor A) \otimes _A X \to M \ctensor X. 
$$
 In fact, $\gamma_{M,X}$  is induced from the natural mappings 
$$
(M \ctensor A) \otimes _A X \to  ( M \otimes _k A/\m_A^n) \otimes _A X 
=    M \otimes _k X/\m_A^n X.
$$

\begin{lemma}
Under the circumstances above, suppose that the left  $A$-module $X$ is finitely generated. 
Then,  $\gamma _{M, X} : \ (M \ctensor A) \otimes _A X \to M \ctensor X$  is surjective for any left  $R$-module $M$. 
\end{lemma}

\begin{pf}
By the assumption, there is a surjective homomorphism of left  $A$-modules
$ f : F = \oplus _{i=1}^{\ell} A e_i \to  X$, where  $F$  is a free left $A$-module of finite rank $\ell$. 
Remark that  $\gamma _{M, F}$  is an isomorphism. 
Naturally we have a commutative diagram
$$
\begin{CD}
(M \ctensor A) \otimes _A F @>{1 \otimes f}>> (M \ctensor A) \otimes _A X \\
@V{\gamma_{M, F}}VV  @V{\gamma _{M, X}}VV  \\
M \ctensor F @>{1 \widehat{\otimes} f}>> M \ctensor X.  \\
\end{CD}
$$ 
Since  the horizontal mappings in the diagram are surjective (see Lemma \ref{exact seq}), and since  $\gamma _{M, F}$  is an isomorphism, we see that  $\gamma _{M, X}$  is surjective. 
\qed
\end{pf}

\begin{prop}\label{flat for formal}
Let  $T = \formal$  be a non-commutative formal power series ring, 
and let  $M$  be an arbitrary left $R$-module. 
Then  $M \ctensor T$  is flat as a right $T$-module. 
\end{prop}

\begin{pf}
To prove the flatness, it is enough to show the following.

\begin{itemize}
\item[(*)]
 For any finitely generated left  ideal  $\ga$, the mapping 
$1 \otimes j : (M \ctensor T)\otimes _T \ga \to (M \ctensor T) \otimes _T T$ 
induced from the inclusion  $j : \ga \to T$  is injective. 
 \par\noindent
\end{itemize}

Note that there is a commutative diagram 
$$
\begin{CD}
(M \ctensor T) \otimes _T \ga @>{1 \otimes j}>> (M \ctensor T) \otimes _T T \\
@V{\gamma_{M, \ga}}VV  @V{\gamma _{M, T}}VV  \\
M \ctensor \ga  @>{1 \widehat{\otimes} j}>> M \ctensor T,  \\
\end{CD}
$$
where we should note that  $\gamma _{M, \ga}$  is an isomorphism, since  $\ga$  is a free module of finite rank by Lemma \ref{f g ideals are free}.
Thus, to prove the proposition, it is sufficient to show that 
$1 \widehat{\otimes} j : M \ctensor \ga \to M \ctensor T$  is injective. 
 By virtue of Lemma \ref{exact seq}, we only have to show that  the relative topology on  $\ga$  from  $T$  is equal to the $\m_A$-adic topology. 
But this has been proved in Corollary \ref{relative top}.  
\qed\end{pf}

\begin{prop}\label{flat}
Let  $A \in \Ahat$, and let  $M$  be an arbitrary left $R$-module. 
Suppose  $A$  is of the form  $A = T/I$  where $T = \formal$  is a non-commutative formal power series ring and $I$  is an ideal of  $T$  that is finitely generated as a left ideal. 
Then  $M \ctensor A$  is flat as a right $A$-module. 
\end{prop}

\begin{pf}
From the short exact sequence of  left  $T$-modules 
$$
\begin{CD}
0 @>>> I @ >>> T @>>> A @>>> 0, 
\end{CD}
$$
we have the commutative diagram 
$$
\begin{CD}
 @. (M \ctensor T) \otimes _T I @ >>> (M \ctensor T) \otimes _T T @>>> (M \ctensor T) \otimes _T A @>>> 0, \\
 @.   @V{\gamma _{M, I}}VV  @V{\gamma_{M, T}}VV  @V{\gamma _{M, A}}VV \\
 0  @>>> M \ctensor I @ >>> M \ctensor T @>>> M \ctensor A @>>> 0, \\
\end{CD}
$$
where  the both rows are exact sequences by Corollary \ref{relative top} and Lemma \ref{exact seq}. 
We already know that  $\gamma _{M, I}$  and  $\gamma _{M, T}$  are isomorphisms,  since  $I$  is a free module of finite rank by Lemma \ref{f g ideals are free}. 
It follows that  $\gamma _{M, A} : (M \ctensor T) \otimes _T A \to M \ctensor A$  is bijective, which is actually an isomorphism of  left  $R \ctensor A^{op}$-modules. 
Since  $(M \ctensor T) \otimes _T A$  is flat as a right $A$-module by Proposition \ref{flat for formal}, we can conclude that  $M \ctensor A$  is also flat over  $A$. 
\qed
\end{pf}

Note,  in general,  $M \ctensor A$  is not flat  as a right $A$-module.

\begin{ex}\label{not flat}
Let  $T = k \langle\langle x, y \rangle\rangle$ be the non-commutative formal power series ring of two variables and let us consider the closed ideals in  $T$,  
$$
I = \overline{(xy^nx \ | \ n =1,2,\ldots)} \quad  \subseteq  \quad   
J = \overline{(x)}.  
$$ 
Note that an element of  $I$  (resp. $J$)  is a formal infinite sum 
$\sum _{\lambda} c_{\lambda} m_{\lambda}$  with  $c_{\lambda} \in k$  and with   monomials  $m_{\lambda}$  involving  $x$  at least twice (resp. once). 

Consider the mapping  $\varphi : T/J \to T/I$  defined by right multiplication by  $x$, i.e.  $\varphi (a\mod J) = ax\mod I$. 
Note that  $\varphi$  is a well-defined homomorphism of left $T$-modules, and 
 it is injective. 

Now let  $A$  be the residue ring $T/I$ and consider  $\varphi$  to be an injective homomorphism of left  $A$-modules.
Let  $M = \bigoplus _{i=1}^{\infty} e_i k$  be a $k$-vector space of countably infinite dimension. 
Then we can show that the mapping  
$$
(M \ctensor A)\otimes \varphi \ : \ (M \ctensor A) \otimes _A T/J \longrightarrow  (M \ctensor A) \otimes _A T/I = M \ctensor A
$$
 is not injective. 
In fact, an element  $z = \sum _{i=1}^{\infty} e_i \otimes y^ixy^i  \in  M \ctensor A$  is mapped to   $\varphi (z) =  \sum _{i=1}^{\infty} e_i \otimes y^ixy^ix$  by  $\varphi$,  which is zero in  $M \ctensor A$. 
However,  $z$  never belongs to  $(M \ctensor A)J$, 
because  any element of  $(M \ctensor A)J$  is a finite sum of the form 
$\sum _i z_ij_i \ (z_i \in M \ctensor A, j_i \in J)$  and 
$z$ is never of this form. 

We can conclude from this observation that $M \ctensor  A$ is not flat as a right  $A$-module. 
\end{ex}

\newpage
\section{Universal lifts of chain complexes}

\subsection{Lifts to artinian local algebras}

In this section $k$  is a field and  $R$  is an associative $k$-algebra. 

By a graded left $R$-module $F$, we just mean a direct sum  $F = \bigoplus _{i \in \Z} F_i$  where  each  $F_i$  is a left  $R$-module.  
If  $F$  is  a  graded left $R$-module and  if  $j$ is an integer, then the shifted graded left  $R$-module  $F[j]$  is defined to be  $F[j]_i = F_{i+j}$  for any $i \in \Z$.  
A graded homomorphism $f : F \to G$ of graded left $R$-module  is an $R$-homomorphism with  $f(F_i) \subseteq  G_i$  for any  $i \in \Z$. 
If  $f : F \to G$  is a graded homomorphism, we denote by $f_i$  the restriction of $f$ on  $F_i$  for each $i$.  
We refer to a graded homomorphism $F \to G[j]$ as a graded homomorphism of degree $j$.

By a chain complex of left $R$-modules or simply a complex over  $R$, we mean a pair  $\Fx = (F, d)$   where  $F$ is a graded left $R$-module  and $d$  is a graded homomorphism of degree $-1$ such that  $d^2 =0$. 
A  complex  $\Fx = (F, d)$  over $R$ is described as 
$$
\begin{CD}
\cdots @>>> F_{i+1} @>{d_{i+1}}>> F _i @>{d_i}>> F_{i-1} @>>> \cdots.
\end{CD}
$$
We say that a complex  $\Fx = (F, d)$ is a projective complex over  $R$  if the underlying graded left $R$-module  $F$  is projective.

Let  $\Fx = (F, d)$  and  $\Gx = (G, d')$  be  chain complexes over  $R$. 
A chain homomorphism  $f : \Fx \to \Gx$  of degree $j$  is a graded homomorphism $f : F \to G[j]$  satisfying  $f\cdot~d +(-1)^{j+1} d'\cdot f = 0$. 
A chain isomorphism  $f : \Fx \to \Gx$  of complexes is a chain homomorphism of degree $0$  that is bijective. 
If there is a chain isomorphism between  $\Fx$  and  $\Gx$, then we say that they are isomorphic as chain complexes over  $R$  and we denote it by $\Fx \cong \Gx$.

Now let  $\Fx$  and  $\Gx$ be projective complexes over $R$, and let  $f, g : \Fx \to \Gx$  be chain homomorphisms of degree $j$. 
We say that  $f$  and  $g$  are homotopically equivalent, denoted by  $f \sim g$, if  there is a graded homomorphism  $h : F \to G [j+1]$  such that 
$f =g + (h \cdot d + (-1)^j d' \cdot h )$. 
We denote  the set of all the homotopy equivalence classes of chain homomorphisms of degree $j$  by  $\Ext _R ^{-j} (\Fx , \Gx)$, which is clearly equipped with structure of $k$-vector space.

For graded homomorphisms  $f : F \to F[j]$  and  $g : F \to F[{\ell}]$, we define a graded homomorphism $[f, g] : F \to F[j+\ell]$ by 
$$
[f, g] = f\cdot g + (-1)^{j+\ell} g \cdot f.
$$  
Note that  $f: F \to F[j]$  is a chain homomorphism if and only if  $[d, f]=0$. Also note that  $f \sim 0$  if and only if there is a graded homomorphism $g : F \to F[j+1]$  with  $f = [d, g]$.

Let  $\varphi : R \to S$  be a $k$-algebra homomorphism and let  $\Fx = (F,d)$  be a projective complex over $R$. 
In this case, we denote by  $S_{\varphi}$  (resp.  ${}_{\varphi}S$)  the left (resp. right)  $S$-module $S$ with right (resp. left)  $R$-module structure through  $\varphi$. 
Then the chain complex  $S_{\varphi} \otimes _R \Fx$  (resp.  $\Fx \otimes _R {}_{\varphi}S$)  of projective left (resp. right)  $S$-modules is defined to be  $(S_{\varphi} \otimes _R F, \ S_{\varphi} \otimes _R d)$ (resp. $(F \otimes _R {}_{\varphi}S, \ d \otimes _R {}_{\varphi}S)$).

Recall that we denote by  $\A$  the category of artinian local $k$-algebras $A$  with  $A/\m_A \cong k$  and $k$-algebra homomorphisms. 
If  $F$  is a graded projective (resp. free) left $R$-module and if  $A \in \A$, then  $F \otimes _k A$  is a graded projective (resp. free)  left $R\otimes _k A^{op}$-module.

\begin{defn}
Let  $\Fx = (F, d)$  be a  projective complex over $R$ and let  $A \in \A$. 
We say that a projective complex $(F \otimes _k A, \Delta)$ over  $R \otimes _k A^{op}$ is {\bf a lifting chain complex of}  $\Fx$  {\bf to} $A$  (or simply a {\bf lift of} $\Fx$  {\bf to} $A$)  if  it satisfies  the equality   $(F \otimes _k A, \Delta) \otimes _A k  = \Fx$.

To be more general, let  $\varphi : A \to B$ be a morphism in  $\A$. 
A projective complex  $(F \otimes _k A, \Delta _A)$ over $R \otimes _k A^{op}$    is said to be a {\bf lift} of a projective complex  $(F \otimes _k B, \Delta _B)$  over  $R \otimes _k B^{op}$  if  it satisfies the equality 
$(F \otimes _k A, \Delta _A) \otimes _A (_{\varphi}B) =  (F \otimes _k B, \Delta _B)$. 
And a projective complex  $(F \otimes _k B, \Delta _B)$ over $R \otimes_k B^{op}$ is said to be {\bf liftable} to $A$ if there is a lift of   $(F \otimes _k B, \Delta _B)$ to  $A$. 
\end{defn}

The aim of this section is to construct a universal one among those lifts of a given projective complex $\Fx = (F, d)$  over $R$.
For this, in the rest of this paper,   $\Fx = (F, d)$ always denotes a fixed projective complex over $R$.

\begin{lemma}
Let  $A \in \A$. 
Then, since  $A$  is of finite dimension as a $k$-vector space, we may take a $k$-basis $\{ 1\} \cup \{ x_1, \ldots , x_r\} \cup  \{ y_j | \ 1\leq j \leq s \}$  of  $A$ so that 
 $\{ x_1 \ldots , x_r\}$  yields  a $k$-basis of  $\m_A/\m_A^2$  and
 $\{ y_j | \ 1\leq j \leq s\}$  is a $k$-basis of  $\m_A^2$.

\begin{itemize}\label{description}
\item[$(a)$]
For any graded homomorphism  $\Delta : F \otimes _k A  \to F \otimes _k A[n]$  of left  $R \otimes_k A^{op}$-modules,  there uniquely exist graded homomorphisms  $f , g_i , h_j : F \to F [n] \ (1 \leq i \leq r, \ 1 \leq j \leq s)$  of left  $R$-modules such that 
$$
\Delta  = f \otimes 1 + \sum_{i=1}^r g _i \otimes x_i + \sum_{j=1}^s h_j \otimes y_j. 
$$ 

\item[$(b)$]
Let $(F \otimes _k A, \Delta)$  be  a lift  of  $\Fx = (F, d)$ to  $A$. 
Then $\Delta$  has a description as in  $(a)$, where   $f = d$  and  
each  $g_i : F \to F[-1] \ (1 \leq i \leq r)$  is  a chain homomorphism. 
\end{itemize}
\end{lemma}

\begin{pf}
(a) 
For any $z  \in F_{\ell}$, we can uniquely write 
$\Delta (z) = z_0 \otimes 1 + \sum _{i=1}^r z_i \otimes x_i + \sum_{j=1}^s w_j \otimes y_j$  for some  $z_0, z_i, w_j \in F_{\ell+n}\ (1 \leq i \leq r, \ 1 \leq j \leq s)$. 
Then, define  $f , g_i , h_j : F \to F [n]$   by  $f(z) = z_0$, $g_i (z) = z_i$  and  $h_j(z) = w_j$  and it is easy to see that they are graded homomorphisms of left  $R$-modules. 
Note that  
$$
\Delta (z \otimes a ) 
= f (z) \otimes a + \sum_{i=1}^r g _i(z) \otimes x_ia + \sum_{j=1}^s h_j(z) \otimes y_ja, 
$$   
for any  $z \otimes a \in F \otimes _k A$.

(b)
Since  $d = \Delta  \otimes _A k =  \Delta  \otimes _A A/\m_A$, 
we have  $f =d$. 
Similarly, we have  $\Delta \otimes _A A/ \m_A ^2 = d \otimes 1 + \sum _{i=1}^r g_i\otimes x_i$  as a graded homomorphism  $F \otimes _k A/ \m_A ^2 \to F \otimes _k A/ \m_A ^2 [-1]$. 
Since  $(\Delta \otimes _A A/ \m_A ^2 )^2 =0$ and  since  $(d \otimes 1)^2=0$, it follows that 
$\sum_{i=1}^r dg_i\otimes x_i + g_id \otimes x_i = \sum_{i=1}^r [d, g_i] \otimes x_i = 0$  as a graded homomorphism on  $F \otimes _k A/ \m_A ^2$.  
Hence we have  $[d, g_i]=0$ for all $i$. 	
\qed\end{pf}

\begin{cor}\label{lift to m^2=0}
Let  $A \in \A$  and suppose   $\m_A ^2 =0$. 
Then, for any lift $(F \otimes _k A, \Delta)$  of  $\Fx = (F, d)$  to $A$, 
the differentiation  $\Delta$  is given by 
$$
\Delta  = d \otimes 1 + \sum_{i=1}^r g _i \otimes x_i, 
$$ 
where  $\{ x_1 , \ldots , x_r\}$ is a $k$-basis of  $\m_A$  and each $g_i : F \to F[-1]$  is a chain homomorphism $(1 \leq i \leq r)$. 
\end{cor}

\begin{lemma}\label{lifting homo}
Let  $\varphi : A \to B$  be a surjective morphism in $\A$  and 
let  $(F \otimes _k B, \Delta)$ be a lifting chain complex of  $\Fx = (F, d)$ to $B$. 

\begin{itemize}
\item[$(a)$]
Any graded homomorphism  $\alpha : F \otimes _k B \to F \otimes _k B$ of graded  $R \otimes _k B^{op}$-modules is liftable to a graded homomorphism  $F \otimes _k A \to F \otimes _k A$ of graded $R \otimes _k A^{op}$-modules. 
That is, there is a graded homomorphism  $\beta : F \otimes _k A \to F \otimes _k A$ with  $\beta \otimes _A \ {}_{\varphi}B = \alpha$. 

\item[$(b)$]
If  $\alpha$  is an isomorphism in $(a)$, then  $\beta$  is also an isomorphism.
\end{itemize}
\end{lemma}

\begin{pf}
(a) 
Since  $F \otimes _k A$  is a left projective $R \otimes _k A^{op}$-module, and since  $1 \otimes \varphi : F \otimes _k A \to F \otimes _kB$  is a surjective homomorphism,  one can find a left  $R\otimes _k A^{op}$-homomorphism $\beta$  which makes the following diagram commutative. 
$$
\begin{CD}
F\otimes _k A  @>\beta>>  F\otimes_k A \\
@V{1\otimes \varphi}VV   @V{1\otimes \varphi}VV \\
F\otimes _k B  @>\alpha>>  F\otimes_k B \\
\end{CD}
$$

(b)
To prove that  $\beta$  is an isomorphism, we may assume that  $\varphi :A \to B$  is a small extension, because any surjective morphism in  $\A$  is a composition of a finite sequence of small extensions. 
So we may have a short exact sequence 
$$
\begin{CD}
0 @>>> k @>\epsilon>> A  @>>>  B @>>> 0.  \\
\end{CD}
$$
Hence, we have a commutative diagram of  left  $R \otimes _k A^{op}$-modules 
$$
\begin{CD}
0 @>>> F  @>\epsilon>> F \otimes _k A  @>>>  F \otimes _k B @>>> 0  \\
@. @V{\alpha \otimes k}VV @V{\beta}VV  @V{\alpha}VV \\
0 @>>> F  @>\epsilon>> F \otimes _k A  @>>>  F \otimes _k B @>>> 0.  \\
\end{CD}
$$
Since $\alpha$  is an isomorphism, it is clear that so is   $\beta$. 
 \qed\end{pf}

\begin{cor}\label{chain iso}
Let  $\varphi : A \to B$  be a surjective morphism in $\A$  as in the lemma. 
Suppose we have two chain complexes  $(F \otimes _k B, \Delta _1)$ and  $(F \otimes _k B, \Delta_2)$  which are lifts of  $\Fx$  to  $B$  and are isomorphic to each other as chain complexes over $R\otimes _k B^{op}$.  
If  $(F \otimes _k B, \Delta _1)$ is liftable to  $A$, then so is   $(F \otimes _k B, \Delta_2)$.  
\end{cor}

\begin{pf}
By the assumption, there is a graded isomorphism  $\alpha : F \otimes _k B \to F \otimes _k B$  such that  $\Delta _2 = \alpha \Delta _1 \alpha ^{-1}$. 
Let  $(F\otimes _k A, \Delta'_1)$  be a lift of  $(F \otimes _k B, \Delta _1)$ to  $A$. 
By Lemma \ref{lifting homo}, $\alpha$  is lifted to an isomorphism  $\beta : F \otimes _k A \to F \otimes _k A$. 
Then it is easy to see that 
$(F\otimes _k A, \beta \Delta'_1 \beta^{-1})$  is a lift of  
$(F \otimes _k B, \Delta _2)$ to $A$.   
\qed\end{pf}

\begin{lemma}\label{difference}
Let  $(A', \epsilon )$  be a small extension of  $A \in \A$, and let  
$(F \otimes _k A, \Delta)$  be  a lift of  $\Fx$ to  $A$. 
Suppose  that chain complexes  $(F \otimes _k A', \Delta _1)$ and  $(F \otimes _k A', \Delta _2)$  are lifts of  $(F \otimes _k A, \Delta)$ to  $A'$. 

\begin{itemize}
\item[$(a)$] 
Then there is a chain homomorphism  $h : F \to F[-1]$  such that  $\Delta  _2 = \Delta _1 + h \otimes \epsilon$. 

\item[$(b)$] 
The following two conditions are equivalent. 

\begin{itemize}
\item[$(1)$]
The equivalence class $[h] \in \Ext _R^1(\Fx, \Fx)$  is zero. 
\item[$(2)$]
There is an isomorphism  $\varphi : (F \otimes _k A', \Delta _1) \to (F \otimes _k A', \Delta _2)$  of chain complexes over $R \otimes _k A'^{op}$  such that  $\varphi \otimes _{A'} A$  is the identity mapping on  $F \otimes _k A$. 
\end{itemize}
\end{itemize}
\end{lemma}

\begin{pf}

(a)
We can take a $k$-basis of  $\m_{A'}$  containing $\epsilon$  as a member. 
Then, both  $\Delta _1$  and  $\Delta _2$  have the descriptions as in Lemma \ref{description}. 
Since  $\Delta _1 \otimes _{A'} A = \Delta = \Delta _2 \otimes _{A'} A$, the difference  $\Delta _2 - \Delta _1$  has a description  $h \otimes \epsilon$. 
We have to show that  $h$  is a chain map. 
Since  $\Delta _1 ^2 = \Delta _2 ^2 =0$ and  $\epsilon ^2 =0$, we have 
$$
0 = \Delta _2 ^2 = (\Delta _1 + h \otimes \epsilon)^2
= \Delta _1 \cdot (h\otimes \epsilon ) + (h\otimes \epsilon) \cdot \Delta _1. 
$$
 Note that  $\Delta _1 \cdot (1\otimes \epsilon ) = d \otimes \epsilon = (1\otimes \epsilon) \cdot \Delta _1$, since  
$\Delta _1$  has a description 
$$
\Delta _1 = d \otimes 1 + \sum _i g_i \otimes x_i + \sum _j h_j \otimes y_j, 
$$ 
as in Lemma \ref{description} and  $\epsilon x_i = x_i \epsilon = 0$  so on. 
Therefore, we have  $(dh + hd) \otimes \epsilon = 0$, hence  $[d, h]=0$. 

(b)
$[(1) \Rightarrow (2)]$: 
If $[h]=0$, then there is a graded homomorphism  $g : F \to F$  of degree $0$  such that  $h = [d, g]$. 
Define a mapping  $\varphi : F \otimes _k A' \to F \otimes _k A'$  by 
$
\varphi = 1 \otimes 1 + g \otimes \epsilon, 
$
which maps  $x \otimes a \in F \otimes _k A'$  to  $x \otimes a + g(x) \otimes \epsilon a $.
Then it is easy to see that  $\varphi$  is an automorphism of a graded left  $R \otimes _k A'^{op}$-module, and the inverse is given by 
$\varphi ^{-1} = 1 \otimes 1 - g \otimes \epsilon$. 
Then, we have the following equalities. 
$$
\begin{array}{rl}
\varphi ^{-1} \Delta _1 \varphi &= 
(1 \otimes 1 - g \otimes \epsilon) \Delta _1 (1 \otimes 1 +  g \otimes \epsilon) \\
&= 
\Delta _1 + \Delta _1 (g \otimes \epsilon)- (g \otimes \epsilon) \Delta _1 \\
&= 
\Delta _1 + (dg  \otimes \epsilon - gd \otimes \epsilon) \\
&= 
\Delta _1 + [d, g] \otimes \epsilon = \Delta _2 \\
\end{array}
$$
Therefore, $\varphi$ satisfies the conditions in  (2).

$[(2) \Rightarrow (1)]$: 
By Lemma \ref{description}, we have a description 
 $\varphi = 1 \otimes 1 + g\otimes \epsilon$  and  $\varphi ^{-1} = 1 \otimes 1 - g \otimes \epsilon$  for some graded homomorphism  $g : F \to F$  of degree $0$. 
Hence, by the same computation as above, we have 
$$
\Delta _2 = \varphi ^{-1} \Delta _1 \varphi = \Delta _1 + [d, g] \otimes \epsilon.  
$$
Therefore, $h = [d, g] \sim 0$.  
\qed\end{pf}

\begin{prop}\label{fiber prod}
Let 
$$
\begin{CD}
B @>>> A_1 \\
@VVV  @V{a_1}VV  \\
A_2 @>{a_2}>> A 
\end{CD}
$$
be a diagram of a fiber product in  $\A$ with  $a_2$  being a surjective map. 

\begin{itemize}
\item[$(a)$]
Let 
$\varphi _1 : F \otimes _k A_1 \to F \otimes _k A_1 [j]$  and 
$\varphi _2 : F \otimes _k A_2 \to F \otimes _k A_2 [j]$  be graded homomorphisms of degree $j$  such that  
$\varphi \otimes _{A_1} A = \varphi _2 \otimes _{A_2} A (= \varphi)$. 
Then there is a graded homomorphism  $\Phi : F \otimes B  \to F \otimes _k B [j]$  of degree $j$  with  $\Phi \otimes _{B} A_i = \varphi _i$  for $i =1,2$. 

\item[$(b)$]
Let 
$(F \otimes _k A_1, \Delta _1)$  and $(F \otimes _k A_2, \Delta _2)$  
be lifts of a chain complex  $(F \otimes _k A, \Delta )$. 
Then there is a chain complex  $(F \otimes _k B, \Delta _B)$  which is a lift of both of  $(F \otimes _k A_1, \Delta _1)$  and $(F \otimes _k A_2, \Delta _2)$. 

\item[$(c)$]
Let 
$(F \otimes _k A_1, \Delta _1)$  and $(F \otimes _k A_2, \Delta _2)$  
be chain complexes such that there is an isomorphism 
 $$
(F \otimes _k A_1, \Delta _1 ) \otimes _{A_1} A \cong (F \otimes _k A_2,   \Delta _2 ) \otimes _{A_2} A 
$$ 
of chain comlexes over  $R \otimes _k A^{op}$. 
Then there is a chain complex  $(F \otimes _k B, \Delta _B)$  which satisfies  
  $(F \otimes _k B , \Delta _B) \otimes _{B} A_i \cong (F \otimes _k A_i, \Delta _i)$  for  $i =1,2$. 
\end{itemize}
\end{prop}

\begin{pf}
(a) Since there is commutative diagram with exact rows 
$$
\begin{CD}
0 @>>> F \otimes _k B @>>> (F \otimes _k A_1) \times (F \otimes _k A_2) @>>> F \otimes _k A @ >>> k \\ 
@.  @.  @V{\varphi _1 \times \varphi _2}VV  @V{\varphi}VV  \\
0 @>>> F \otimes _k B [j] @>>> (F \otimes _k A_1) [j] \times (F \otimes _k A_2) [j]  @>>> F \otimes _k A [j] @ >>> k,  \\ 
\end{CD}
$$
it induces a mapping  $\Phi : F \otimes _k B \to F \otimes _k B [j]$. 

(b) Just apply (a) to  $\Delta _1$  and  $\Delta _2$, and we get a graded homomorphism   $\Delta _B : F \otimes _k B \to F \otimes _k B[-1]$. 
It is clear that  $\Delta _B ^2 = 0$. 

(c) By definition of isomorphisms of chain complexes, there is a graded isomorphism  $\alpha : F \otimes _k A \to F \otimes _k A$  of graded $R \otimes _k A^{op}$-modules, such that 
$\Delta _1 \otimes _{A_1} A = \alpha \cdot (\Delta _2 \otimes _{A_2} A) \cdot \alpha ^{-1}$. 
Since  $a_2$  is surjective, 
there is a graded isomorphism $\beta : F \otimes _k {A_2} \to F \otimes _k A_2$  which lifts  $\alpha$, by Lemma \ref{lifting homo}. 
Put  $\Delta _2' = \beta \cdot \Delta _2 \cdot \beta ^{-1}$ and apply (b) to   the chain complexes  $(F \otimes _k A_1,  \Delta _1)$  and  
$(F \otimes _k A_2, \Delta _2')$, and we obtain a chain complex 
$(F \otimes _k B, \Delta_B)$  which is a lift of the both of them.     
\qed\end{pf}

\subsection{Construction of maximal lifts}

As in the previous section, let  $R$  be an associative algebra over a field $k$  and let  $\Fx = (F, d)$  be a projective complex over $R$. 
In the rest of the paper we always assume that 

\begin{equation}\label{assumption}
r = \Kdim _k \ \Ext _R^1(\Fx, \Fx ) < \infty.  
\end{equation}

Under this assumption, we take chain homomorphism 
$t_i^* : F \to F[-1] \ ( 1 \leq i \leq r)$  whose equivalence classes 
$\{ [t_1^*],\ldots , [t_r^*]\}$  is a $k$-basis of  $\Ext _R^1(\Fx, \Fx )$. 
We take variables $t_1, \ldots , t_r$ corresponding to this basis, and consider  the non-commutative formal power series ring  $T = \formal$. 
Now define  $\delta : F \otimes _k T/\m_T^2 \to F \otimes _k T/\m_T^2$  by 

\begin{equation}\label{delta}
\delta = d \otimes 1 + \sum _{i=1}^r t_i^* \otimes t_i. 
\end{equation}

It follows from Corollary \ref{lift to m^2=0} that 
$(F \otimes _k T/\m_T^2, \delta )$  is a lift of  $\Fx$  to  $T/\m_T^2$.

\begin{defn}
Let  $I$  be a closed ideal of  $T$. 
We define the complete tensor product of a graded projective left $R$-module with  $T/I$  as follows:
$$
\begin{CD}
F \ctensor T/I := \bigoplus _{i\in \Z} (F_i \ctensor T/I). 
\end{CD}
$$
\end{defn}

Now let  $I$  be a closed ideal of  $T$  and  let  $(F \ctensor T/I, \ \Delta)$  be a chain complex. 
If  $I' \supseteq  I$  be another closed ideal of  $T$, then there is a natural projection  $T/I \to T/I'$, which induces, by Lemma \ref{exact seq},  a surjective homomorphism 
$$
F \ctensor  T/ I \to F \ctensor T/I'. 
$$
Thus we have a surjective homomorphism of left  $R \ctensor (T/I)^{op}$-modules 
$$
\pi : (F \ctensor  T/ I)  \otimes_{T/I}  T/I'  \to F \ctensor T/I'.
$$
Note that  $\pi$  may not be an isomorphism,   

For each  $n \geq 1$,  $\Delta$  induces a graded homomorphism  
$\Delta _n  : F \otimes _k (T/I+\m_T^n)  \to   F \otimes _k (T/I+\m_T^n)  [-1]$  and it holds that 
$\Delta = \varprojlim \Delta _n$. 
Since  $I \subseteq I'$,  each  $\Delta _n$ induces a graded homomorphism  
$\Delta' _n  : F \otimes _k T/(I' +\m_T^n) \to   F \otimes _k (T/I'+\m_T^n) [-1]$, 
  and we obtain a graded homomorphism 
$\Delta ' = \varprojlim  \Delta '_n  : F \ctensor T/I' \to   F \ctensor T/I' [-1]$. 
By an abuse of notation, we denote this mapping  $\Delta '$  by  $\Delta \otimes _{T/I} T/I'$.

\begin{defn}
Let  $I_1 \subseteq I_2 \subseteq \m_T ^2$  be closed ideals of  $T$. 
A chain complex  $(F \ctensor T/I_1, \Delta _1)$  is called a lift of a chain complex  $(F \ctensor T/I_2,  \Delta _2)$  if  $\Delta _2 = \Delta _1 \otimes _{T/I_1} T/I_2$. 
\end{defn}

\begin{defn}
Let  $T/I \in \Ahat$ where  $T = \formal$  is  a non-commutative formal power series ring  and  $I \subseteq \m_T^2$. 
And  let  $(F \ctensor T/I, \Delta )$  be  any chain complex which is a lift of  $\Fx$. 
Now we consider the following set of lifting chain complexes of  $(F \ctensor T/I, \ \Delta )$:   
$$
\begin{array}{rl}
\I (I, \Delta)  = 
&\{ (T/I', \ \Delta ' )\ | \   I' \ \text{is a closed ideal of} \ T \ \text{with } \  I' \subseteq I  \  \text{and} \\ 
&(F \ctensor T/I', \ \Delta') \ \text{is a lifting chain complex of } \ (F \ctensor T/I, \ \Delta) \} 
\end{array}
$$
We define an order relation on the set  $\I (I, \Delta) $ as follows: 
$$
(T/I_1, \ \Delta _1) > (T/I_2, \ \Delta _2)
\Longleftrightarrow
I_1 \subseteq I_2 \ \ \text{and } \ \ \Delta _1 \otimes _{T/I_1} T/I_2 = \Delta _2
$$
\end{defn}

\begin{lemma}\label{Zorn}
The ordered  set  $\I (I, \Delta) $  is an inductively ordered set. 
In particular, there exists a maximal element $\I (I, \Delta)$.    
\end{lemma}

\begin{defn}
If  $(T/I_0, \ \Delta _0)$ is a maximal element in  $\I (I, \Delta) $  as in the lemma, 
then we say that  the chain complex  $(F \ctensor (T/I_0), \Delta _0)$  is {\bf  a maximal lift of}  $(F \ctensor T/I, \ \Delta)$.
\end{defn}

\begin{pf}
Let  $\{ (T/I_{\lambda}, \Delta _{\lambda}) \ | \ \lambda \in \Lambda \}$  be a totally ordered subset of  $\I  (I, \Delta)$. 
Note that  $J = \bigcap _{\lambda \in \Lambda} I_{\lambda}$  is a closed ideal of  $T$  and that  $\varprojlim  T/I_{\lambda}= T/J$.   
Hence  $(F \ctensor T/J, \Delta _J) = \varprojlim _{\lambda} (F \ctensor T/I_{\lambda}, \Delta _{\lambda})$  is a lifting chain complex of  $(F \ctensor T/I, \Delta)$. 
Therefore, $(T/J, \Delta _J) \in \I (I, \Delta)$  and  $(T/J, \Delta _J) > (T/I_{\lambda}, \Delta _{\lambda})$  for any $\lambda \in \Lambda$. 
Thus  $\I (I, \Delta)$  is an inductively ordered set. 
The existence of maximal element of $\I (I, \Delta)$ follows from Zorn's lemma.
\qed\end{pf}

We should remark the following

\begin{lemma}\label{not maximal}
If  $(T/I_1, \ \Delta _1)$ is not a  maximal element in  $\I (I, \Delta) $, then  there is a nontrivial small extension   $T/I_2$   of  $T/I_1$  such that 
 $(T/I_2, \Delta _2) \in  \I (I, \Delta)$ is strictly bigger than  $(T/I_1, \ \Delta _1)$,  for some  $\Delta _2$
\end{lemma}

\begin{pf}
Take a strictly bigger element $(T/I_1', \ \Delta _1') > (T/I_1, \ \Delta _1)$  in  $\I (I, \Delta)$.  
Since  $I_1' \subsetneq I_1$ are closed ideals,  there is an integer $n$ with   $\overline{I_1' + (\m_T ^n \cap I)} \not= I_1$. 
In fact,  if not, we will have   $I_1 \subseteq I_1' + \m _T ^n$ for any $n$, 
because the right hand side is the closed ideal containing  $I_1' + (\m _T ^n \cap I)$. 
Then we shall have  $I_1 \subseteq \overline{I_1'} = I_1'$, a contradiction. 
 
Now, since  $A = T/\overline{I_1' + (\m_T ^n \cap I)}$  is a complete local $k$-algebra and  since  the image  $J_1$  of  $I_1$ in  $A$   is an ideal of $A$  of finite length, we can find a closed ideal $J_2$ of  $A$  contained in  $J_1$  with  $\length (J_1/J_2)=1$. 
See Proposition \ref{condition for closed ideal}. 
Taking the inverse image  of  $J_2$ in  $T$,  we have a closed ideal  $I_2$  of  $T$  contained in  $I_1$  and  $\length (I_1/I_2) =1$. 
Finally set  $\Delta _2 = \Delta _1' \otimes _{T/I_1'} T/I_2$, and we easily see that  $(T/I_2, \Delta _2)$  meets the requirements. 
\qed\end{pf}

The following is an easy consequence of  Lemma \ref{lifting homo}.

\begin{lemma}\label{lifting homo complete1} 
Let  $\varphi : A \to B$  be a surjective morphism in $\Ahat$  where  $\kernel (\varphi)$  is of finite length, and let  $(F \ctensor  B, \Delta)$ be a lifting chain complex of  $\Fx$. 

\begin{itemize}
\item[$(a)$]
Any graded homomorphism  $\alpha : F \ctensor B \to F \ctensor B$ is liftable to a graded homomorphism  $F \ctensor  A \to F \ctensor A$. 
That is, there is a graded homomorphism  $\beta : F \ctensor  A \to F \ctensor  B$ with  $\beta \otimes _A \ {}_{\varphi}B = \alpha$. 

\item[$(b)$]
If  $\alpha$  is an isomorphism in $(a)$, then  $\beta$  is also an isomorphism. \end{itemize}
\end{lemma}

\begin{pf}
(a) 
By induction on the length of $\kernel (\varphi)$, we may assume that $A \to B$  is a small extension. 
In this case, it is easily seen that the following diagram is a pull-back diagram of right $A$-modules  for any integer $n$  which satisfies  $\kernel (\varphi) \cap \m _A ^n = (0)$. 
$$
\begin{CD}
A/\m_A ^{n+1} @>{\varphi _{n+1}} >> B/\m_B^{n+1} \\
@VVV    @VVV  \\
A/\m_A ^{n} @>{\varphi _n}>> B/\m_B^{n},  \\
\end{CD}
$$
where  $\varphi _n$  is the induced mapping by $\varphi$  and the vertical arrows are natural projections. 
Thus the diagram 
$$
\begin{CD}
F \otimes _k A/\m_A ^{n+1} @>{1 \otimes \varphi _{n+1}} >> F \otimes _kB/\m_B^{n+1} \\
@V{p_n}VV    @V{q_n}VV  \\
F \otimes _kA/\m_A ^{n} @>{1\otimes \varphi _n}>> F \otimes _kB/\m_B^{n},  \\
\end{CD}
$$
is a pull-back diagram of  $R \ctensor A^{op}$-modules.
Denote by  $\alpha _n$  the mapping  $\alpha \otimes _k B/\m_B^n : F \otimes _k B/\m _B^n \to F \otimes _k B/\m _B^n$. 
If we have an $R \otimes  _k A/\m _A ^n$-homomorphism $\beta _n : F \otimes _k A/\m _A ^n \to F \otimes _k A/\m _A ^n$ with 
$(1\otimes\varphi_{n})\cdot\beta _n \cdot p_n = q_n \cdot \alpha _{n+1}\cdot (1\otimes\varphi_{n+1})$, then   it follows that there uniquely exists  $\beta _{n+1} : F \otimes _k A/\m _A ^{n+1} \to F \otimes _k A/\m _A ^{n+1}$ such that 
$(1 \otimes \varphi_{n+1})\cdot \beta  _{n+1}= \alpha _{n+1}\cdot (1\otimes\varphi_{n+1})$  and 
$p_n \cdot \beta _{n+1} = \beta _n\cdot p_n$. 
Therefore, by induction, we have such  $\beta  _n$ for all $n \geq1$. 
Then, setting  $\beta = \varprojlim \beta _n$, we see that 
$\beta$ is a lift of the mapping   $\alpha$.

(b)
In the proof above, if $\alpha$  is an isomorphism, then each $\beta _n$  is also an isomorphism by Lemma \ref{lifting homo}, hence so is $\beta$.
\qed\end{pf}

\begin{lemma}\label{max chain iso}
Let   $(F \ctensor T/I_0, \ \Delta _0)$  be a maximal lift of $(F \ctensor T/I, \ \Delta)$. 
Then, any chain complex  $(F \ctensor T/I_0, \ \Delta _1)$  which is isomorphic to  $(F \ctensor T/I_0, \ \Delta _0)$  as a complex over $R \ctensor (T/I_0)^{op}$  is also a maximal lift of  $(F \ctensor T/I, \ \Delta)$.  
\end{lemma}

\begin{pf}
There is a graded isomorphism $\alpha : F\ctensor T/I_0 \to F \ctensor T/I_0$  with  $\Delta _1 = \alpha \Delta _0 \alpha ^{-1}$. 
If  $(F \ctensor T/I_0, \ \Delta _1)$  is not a maximal lift, then by Lemma \ref{not maximal} there is a nontrivial small extension  $T/I_2$  of  $T/I_0$  such that   $(T/I_2, \Delta _2) > (T/I_0, \Delta _1)$  in  $\I (I, \Delta)$. 
We can lift the isomorphism  $\alpha$  to  $\beta : F \ctensor T/I_2 \to F \ctensor T/ I_2$ by Lemma \ref{lifting homo complete1}.   
Then it is easy to see that  $(T/I_2, \beta ^{-1} \Delta _2 \beta) > (T/I_0, \Delta_0)$ in  $\I (I, \Delta)$, and it contradicts to the assumption. 
\qed\end{pf}

\begin{lemma}\label{automorphism}
Let  $T = \formal$  be a non-commutative formal power series ring.

\begin{itemize}
\item[$(a)$]
For any given  $f_i \in \m_T^2 \ (1 \leq i \leq r)$, we define a $k$-algebra homomorphism  $\varphi : T \to T$  by  $\varphi (t_i ) = t_i + f_i \ (1\leq i \leq r)$. 
Then,  $\varphi$  is an automorphism of  $T$  such that it induces the identity mapping on  $T/\m_T^2$. 

\item[$(b)$]
Any $k$-algebra automorphism of  $T$  which induces the identity on  $T/\m_T^2$  is given as in  $(a)$.

\item[$(c)$]
Let  $I_1 \subseteq I_2 \subseteq \m_T^2$  be closed ideals of  $T$  and 
let  $\psi : T/I_1 \to T/I_2$  be any $k$-algebra homomorphism that induces the identity on  $T/\m_T^2$. 
Then there is a $k$-algebra automorphism  $\varphi : T \to T$  with  $\varphi (I_1) \subseteq I_2$  and  the induced mapping  $\overline{\varphi} : T/I_1 \to T/I_2$  equals  $\psi$. 
\end{itemize}
\end{lemma}

\begin{pf}
(a) 
It is obvious that  $\varphi$  induces the identity on  $\m_T/\m_T^2$. 
Hence it follows from Lemma \ref{surjective} that  $\varphi : T \to T$  is a surjective $k$-algebra homomorphism.
In particular, every induced mapping  $\varphi _n : T/\m_T^n \to T/\m_T^n$  is surjective as well.  
Comparing the lengths we conclude that each  $\varphi _n$ is bijective. 
Hence  $\varphi = \varprojlim _n \varphi _n$ is an automorphism.

(b) 
Trivial. 

(c)
By the assumption, we can choose  $f_i \in \m_T^2 \ (1\leq i \leq r)$  so that 
$\psi (t_i\pmod{I_1}) = t_i+f_i \pmod{I_2} \ (1 \leq i \leq r)$. 
Now define an automorphism  $\varphi : T \to T$  by  $\varphi (t_i) = t_i + f_i \ (1 \leq i \leq r)$, and it is easy to see that $\varphi$ satisfies the desired condition. 
\qed\end{pf}

Now,  as in the beginning of this section, 
we consider the lifting chain complex  
$(F \otimes _k T/\m_T^2, \ \delta)$  with 
$\delta = d \otimes 1 + \sum _{i=1}^r t_i^* \otimes t_i$  as in Equation (\ref{delta}),  where 
 $t_i^* : F \to F[-1] \ ( 1 \leq i \leq r)$  are chain homomorphisms whose equivalence classes $[t_1^*],\ldots , [t_r^*]$  form a $k$-basis of  $\Ext _R^1(\Fx, \Fx)$.

\begin{thm}\label{uniqueness of max lifts}
A maximal lift of  $(F \otimes _k T/\m_T^2, \ \delta)$  is unique up to $k$-algebra automorphisms and chain isomorphisms.
I.e., if  we have two maximal elements  $(T/I_0, \Delta _0)$  and  $(T/I_1, \Delta _1)$  in  $\I (\m_T^2, \delta)$, then there exists a $k$-algebra automorphism  $\varphi : T \to T$ such that  $\varphi$  induces a $k$-algebra isomorphism  $\overline{\varphi} : T/I_0 \to T/I_1$  and  
$(F \ctensor (T/I_0), \Delta _0) \otimes _{T/I_0} \ {}_{\overline{\varphi}}(T/I_1)$  is isomorphic to  $(F \ctensor (T/I_1), \Delta _1)$  as a complex over  $R \ctensor (T/I_1)^{op}$. 
\end{thm}

\begin{pf}
\vspace{2pt}
(1) 
First of all we note from Remark \ref{projlim} that there is an infinite descending  sequence of ideals of  $T$; 
$L_0 = \m_T^2 \supset L_1 \supset L_2 \supset  L_3 \supset  \cdots \supset I_1$ 
 such that 
$\length (L_n/L_{n+1}) =1$  for all $n \geq 1$  and  
$T/I_1 = \varprojlim  T/L_n$. 
Note that  each  $T/L_{n+1} \to T/L_n$  is a small extension in $\A$. 
Note also that 
$F \ctensor T/I_1 = \varprojlim F \otimes _k T/L_n$.

\vspace{4pt}
(2)
By induction on  $n$, we shall construct a $k$-algebra homomorphism 
$$
\overline{\varphi_n} : T/I_0 \to T/L_n, 
$$
 and an automorphism  of a graded  $R\otimes _k (T/L_n)^{op}$-module  
$$
\alpha _n : F \otimes _k (T/L_n) \to F \otimes _k (T/L_n), 
$$
which satisfy the following four conditions.

\begin{itemize}
\item[${\small{(0)}}$]
$\overline{\varphi _0} : T/I_0 \to T/L_0=T/\m_T^2$  is a natural projection and  $\alpha _0 = 1$. 

\item[$(i)$]
The following diagram is commutative: 
$$
\begin{CD}
T/I_0 @.  \\
@VVV  \begin{picture}(-10,-15)
\put(-60,15){\vector(1,-1){25}}
\put(-40,5){\small$\overline{\varphi_{n-1}}$}
\put(-100,5){\small$\overline{\varphi_{n}}$}
\end{picture} \\
T/L_n  @>>> T/L_{n-1},  \\
\end{CD}
$$
where the horizontal map is a natural projection. 

\item[$(ii)$]
$$\alpha _n \otimes _{T/L_n} (T/L_{n-1}) = \alpha _{n-1}$$

\item[$(iii)$]
$$
\Delta _0 \otimes _{T/I_0} \ {}_{\overline{\varphi _n}}(T/L_n) = \alpha _n \cdot \left( \Delta _1 \otimes _{T/I_1} (T/L_n) \right) \cdot \alpha _n ^{-1}
$$
\end{itemize}

\vspace{4pt}
(3)
Suppose we obtain such  $\overline{\varphi _n}$  and $\alpha _n$  for  $n \geq 0$  satisfying the above conditions. 
Then, by $(i)$  and  $(ii)$, we have a $k$-algebra homomorphism 
$$
\overline{\varphi} = \varprojlim \overline{\varphi _n} \ :\  T/I_0 \to T/ I_1
$$
and an automorphism of graded $R \ctensor (T/I_1)^{op}$-modules
$$
\alpha = \varprojlim \alpha _n \ :\  F \ctensor T/I_1 \to F \ctensor T/I_1.
$$ 
And it follows from $(iii)$ that 
$$
\qquad\qquad\qquad\qquad\qquad
\Delta _0 \otimes _{T/I_0} {}_{\overline{\varphi}}(T/I_1) = \alpha  \cdot \Delta _1  \cdot \alpha ^{-1}.  
\qquad\qquad\qquad\qquad\qquad \text{(*)}
$$
Therefore we have the isomorphism 
$$
(F \ctensor (T/I_0), \ \ \Delta _0) \otimes _{T/I_0}  {}_{\overline{\varphi}}(T/I_1)  \cong  (F \ctensor (T/I_1),  \ \ \Delta _1), 
$$
as a chain complex of left $R \ctensor (T/I_1)^{op}$-modules. 

Now we prove that  $\overline{\varphi}$  is an isomorphism. 
Since  $\overline{\varphi}$  induces the identity mapping on  $T/\m_T^2$  by the condition $(0)$, we can apply  Lemma \ref{automorphism} to get a $k$-algebra automorphism $\varphi : T \to T$  with  $\varphi (I_0) \subseteq I_1$ and $\varphi$ induces $\overline{\varphi} : T/I_0 \to T/I_1$. 
Here supposed  $\varphi (I_0) \subsetneqq I_1$. 
Then we would have from (*) that $(F \ctensor (T/I_1), \ \alpha \Delta _1 \alpha ^{-1})$  were liftable to  $(F \ctensor (T/I_0), \ \Delta _0) \otimes _{T/I_0} {}_{\varphi} (T/\varphi(I_0))$. 
This is a contradiction, because  $(F \ctensor (T/I_1), \ \alpha \Delta _1 \alpha ^{-1})$  is a maximal lift by Lemma \ref{max chain iso}. 
Thus we have shown   $\varphi (I_0) = I_1$  and hence 
$\overline{\varphi} : T/I_0 \to T/I_1$  is an isomorphism. 

In such a way, we have verified that the theorem is proved once we have 
  $\overline{\varphi _n}$  and $\alpha _n$  for  $n \geq 0$  satisfying the conditions in  (2).

\vspace{4pt}
(4) 
Now we shall construct   $\overline{\varphi _n}$  and $\alpha _n$  by induction on  $n$. 
To do this, assume we already have   $\overline{\varphi _n}$  and $\alpha _n$ satisfying the conditions in  (2)  for an integer $n \geq 0$. 

We take an  element  $\epsilon \in L_n$ which gives a socle element of  $T/L_{n+1}$  so that 
$L_n = L_{n+1} + (\epsilon)$, 
 and hence we have a small extension 
$$
\begin{CD}
0 @>>> k @>\epsilon>> T/L_{n+1} @>>> T/L_n @>>> 0. 
\end{CD}
$$
By Lemma \ref{automorphism} there is a $k$-algebra automorphism  $\varphi _n : T \to T$  such that  $\varphi _n (I_0) \subseteq L_n$  and  $\overline{\varphi_n}$  is induced from $\varphi _n$.

\vspace{4pt}
(5) 
Under the circumstances as in (4), we claim that
$$
\varphi _n (I_0) \subseteq L_{n+1}.
$$

On the contrary, assume that $\varphi _n (I_0) \not\subseteq L_{n+1}$. 
Then, since  $\varphi _n (I_0) \subseteq L_n$, we have $L_n = L_{n+1} + \varphi _n (I_0)$. 
Therefore, there is a fiber product diagram 
$$
\begin{CD}
T/L_{n+1}\cap \varphi _n (I_0)  @>>> T/\varphi _n (I_0) \\
@VVV  @VVV  \\
T/L_{n+1} @>>> T/L_n. \\
\end{CD}
$$
Now let  $\beta _{n+1}$  be any lift of  $\alpha _n$  to  $T/L_{n+1}$, i.e.
$$
\begin{CD}
F \otimes _k T/L_{n+1}  @>{\beta _{n+1}}>> F \otimes _k T/L _{n+1} \\
@VVV  @VVV  \\
F \otimes _k T/L_{n}  @>{\alpha_n}>> F \otimes _k T/L _{n},  \\
\end{CD}
$$
where  $\beta _{n+1}$  is also an isomorphism of graded left $R \otimes _k (T/L_{n+1})^{op}$-modules by Lemma \ref{lifting homo}. 
Then, the chain complexes   
$$
(F \otimes _k (T/L_{n+1}), \ \ \beta _{n+1}  \cdot \left(\Delta _1 \otimes _{T/I_1} (T/L_{n+1})\right)  \cdot \beta _{n+1}^{-1})
$$  and 
$$
\qquad\qquad\qquad\qquad
(F \ctensor T/ \varphi _n(I_0), \ \  \Delta _0 \otimes _{T/I_0} \ {}_{\varphi_n} (T/\varphi _n (I_0)) 
\qquad\qquad\qquad\qquad \text{(**)}
$$ 
are lifts of 
$(F \otimes _k (T/L_{n}), \ \alpha  _{n}  \cdot \left(\Delta _1 \otimes _{T/I_1} (T/L_{n})\right)  \cdot \alpha _{n}^{-1})$  to  $T/L_{n+1}$  and  $T/\varphi _n (I_0)$ respectively. 
And by Lemma \ref{fiber prod}, they are liftable to  $T/ L_{n+1} \cap \varphi _n (I_0)$. 
Note that $(F \ctensor (T/I_0), \ \Delta _0)$,  as well as the chain complex  (**), is a maximal lift of  $(F, d)$. 
Hence it never be liftable to a nontrivial extension ring. 
Thus it follows that  $L_{n+1} \cap \varphi _n (I_0) = \varphi _n (I_0)$, therefore  $\varphi _n (I_0)\subseteq L_{n+1}$ as claimed above.

\vspace{4pt}
(6) 
By the claim (5), the $k$-algebra automorphism  $\varphi _n$ induces a map $\overline{\varphi _n} : T/I_0 \to T/L_{n+1}$. 
Then, the chain complexes   
$$
(F \otimes _k (T/L_{n+1}), \ \ \beta _{n+1} \cdot  \left(\Delta _1 \otimes _{T/I_1} (T/L_{n+1})\right) \cdot  \beta _{n+1}^{-1})
$$  and 
$$
(F \otimes _k (T/ L _{n+1}),  \ \  \Delta _0 \otimes _{T/I_0} \ {}_{\overline{\varphi_n}} (T/L_{n+1})) 
$$ 
are lifts of 
$(F \otimes _k (T/L_{n}), \ \ \alpha  _{n}  \cdot \left(\Delta _1 \otimes _{T/I_1} (T/L_{n})\right)  \cdot \alpha _{n}^{-1})$. 
Therefore, by Lemma \ref{difference}, there is a chain homomorphism $h : F \to F[-1]$ with 
$$
\qquad
\beta _{n+1} \cdot  \left(\Delta _1 \otimes _{T/I_1} (T/L_{n+1})\right) \cdot  \beta _{n+1}^{-1} 
=  \Delta _0 \otimes _{T/I_0} \ {}_{\overline{\varphi_n}} (T/L_{n+1}) + h \otimes \epsilon. 
\qquad \text{(***)} 
$$
Since the classes of  $t_1^*, \ldots , t_r^*$  form a $k$-basis of  $\Ext _R^1 (\Fx, \Fx)$, 
we may describe  $h$  as 
$$
h = \sum _{i=1} ^r c_i t_i^* + [d, H], 
$$
for some $c_i \in k$  and  a graded homomorphism  $H : F \to F$  of degree $0$. 
Now we define a $k$-algebra automorphism  
$$
\varphi _{n+1} : T \to T \qquad \text{by} \qquad \varphi _{n+1} (t_i) = \varphi _n(t_i) + c_i\epsilon
$$
  for  $1 \leq i \leq r$. Then  $\varphi _{n+1}$  is well-defined, because  $\epsilon \in L_n \subseteq \m_T^2$. 
Note that, for $1 \leq i, j \leq r$, we have 
$
\varphi _{n+1} (t_it_j) = \varphi _n (t_it_j) + c_i (\epsilon \varphi _n (t_j) + \varphi _n (t_i)\epsilon) + c_ic_j \epsilon ^2
\equiv \varphi _n (t_it_j) \pmod{L_{n+1}}. 
$
Thus, we see  
$
\varphi _{n+1} (x) \equiv \varphi _n (x) \pmod{L_{n+1}} 
$   
for all $x \in \m_T^2$. 
Therefore by the claim (5)  and by the fact that $I_0 \subseteq \m_T^2$, we have $\varphi _{n+1} (I_0) \subseteq L_{n+1}$, hence  $\varphi _{n+1}$  induces the $k$-algebra map  $\overline{\varphi _{n+1}} : T/I_0 \to T/L_{n+1}$  and the diagram 
$$
\begin{CD}
T/I_0 @.  \\
@VVV  \begin{picture}(-10,-15)
\put(-55,15){\vector(1,-1){25}}
\put(-35,5){\small$\overline{\varphi_{n}}$}
\put(-110,5){\small$\overline{\varphi_{n+1}}$}
\end{picture} \\
T/L_{n+1}  @>>> T/L_{n},  \\
\end{CD}
$$
is commutative. 

By the definition of  $\varphi _{n+1}$, it follows that 
$$
\Delta _0 \otimes _{T/I_0} \ {}_{\overline{\varphi_{n+1}}}(T/L_{n+1}) 
=   
\Delta _0 \otimes _{T/I_0} \ {}_{\overline{\varphi_{n}}}(T/L_{n+1}) 
+ \sum _{i=1} ^r t_i^* \otimes c_i \epsilon, 
$$
thus we see from  (***)  that 
$$
\beta _{n+1}  \cdot \left( \Delta _1 \otimes _{T/I_1} ( T/L_n ) \right)  \cdot \beta _{n+1} ^{-1}
=
\Delta _0 \otimes _{T/I_0} \ {}_{\overline{\varphi_{n+1}}} ( T/L_{n+1} )+ [d, H].
$$
Now define, as in the proof of Lemma \ref{lifting homo}, an automorpshim 
$\alpha _{n+1}$  by 
$$
\alpha _{n+1} = (1- H \otimes \epsilon ) \cdot \beta _{n+1}. 
$$
Then we have that 
$$
\alpha  _{n+1}  \cdot \left( \Delta _1 \otimes _{T/I_1} ( T/L_n ) \right)  \cdot \alpha  _{n+1} ^{-1}
=
\Delta _0 \otimes _{T/I_0} \ {}_{\overline{\varphi_{n+1}}} ( T/L_{n+1} ), 
$$
and also 
$$
\alpha _{n+1} \otimes _{T/L_{n+1}} T/L_n = \beta _{n+1} \otimes _{T/L_{n+1}} T/L_n = \alpha _n. 
$$
Therefore, we have obtained  $\overline{\varphi_{n+1}}$  and  $\alpha _{n+1}$ satisfing the conditions  $(i)$, $(ii)$ and $(iii)$, and the proof is completed by induction. 
\qed
\end{pf}

\subsection{Universal lifts}

As in the previous section, let  $\Fx = (F, d)$  be a projective complex over $R$   that satisfies the fundamental assumption 
$$
 r = \Kdim _k \ \Ext _R ^1 (\Fx, \Fx) < \infty. 
$$
We define a functor $\F : \A  \to (Sets)$ as follows: 
For any  $A \in \A$, we set  
$$
\F (A)  = \{ (F \otimes _k A, \ \Delta) \ | \ \text{it is a lifting chain complex of  $\Fx$  to  $A$} \}/\cong , 
$$  
where $\cong$ denotes  the isomorphism as chain complexes over $R \otimes _k A^{op}$.
If  $f : A \to B$  is a morphism in $\A$, then we define a mapping 
$\F (f) : \F (A) \to \F (B)$  by 
$$
\F (f) ((F \otimes _k A, \ \Delta)) = (F \otimes _k A, \ \Delta) \otimes _A \ {}_f B. 
$$
Note that  $\F$  is a covariant functor.

\begin{defn}
Let  $P \in \Ahat$ and let  $\Lx = (F \ctensor P, \Delta)$ be a lifting chain complex  of $\Fx$  to  $P$. 

\begin{itemize}
\item[$(a)$]
We define a morphism between functors on  $\A$ ; 
$$
\phi _{\Lx} \ : \ \Hom _{k\text{-alg}}(P, \ \ ) \to  \F,  
$$ 
by  
$$
\phi _{\Lx} (f) = (F \ctensor  P, \ \Delta) \otimes _{P} {}_f A
$$
for  $f \in \Hom _{k\text{-alg}}(P, A)$ with  $A \in \A$. 

\item[$(b)$]
We say that  the chain complex  $\Lx$ is {\bf a universal lift of}  $\Fx$  if the morphism  $\phi _{\Lx}$  is an isomorphism. 
Thus, in this case, the functor  $\F$  on the category  $\A$  is pro-representable by  $P \in \Ahat$.  
If  $\Lx = (F \ctensor P, \Delta)$ is a universal lift of   $\Fx$, then  $P$  is called {\bf a parameter algebra of the universal lift of}  $\Fx$. 
\end{itemize}
\end{defn}

\begin{lemma}\label{unique universal lifting}
If there is a universal lift of $\Fx$, then any parameter algebras  of any universal lifts are isomorphic each other as $k$-algebras. 
\end{lemma}

\begin{pf}
In fact, if  $P_1$  and  $P_2$  are such parameter algebras, then we have an isomorphism  $\Hom _{k\text{-alg}}(P_1, \ \ ) \cong \Hom _{k\text{-alg}}(P_2, \ \ )$   as functors on  $\A$. 
Then, it is easy to see that there are isomorphisms  $P_1/\m_{P_1}^n \cong P_2/\m_{P_2}^n$  as $k$-algebras for any $n \geq 1$, which are compatible with the projections  $P_1/\m_{P_1}^{n+1} \to P_1/\m_{P_1}^{n}$ and  $P_2/\m_{P_2}^{n+1} \to P_2/\m_{P_2}^{n}$. 
Hence  $P_1 \cong P_2$.
\qed\end{pf}

\begin{thm}\label{univ=max}
The following two conditions are equivalent for a lifting chain complex 
$\Lx_0 = (F \ctensor P_0, \Delta_0)$  of  $\Fx$ where  $P_0= T/I_0 \in \Ahat$  with  $I_0 \subseteq \m_T^2$.  

\begin{itemize}
\item[$(a)$]
$\Lx _0$  is a universal lift of  $\Fx = (F, d)$. 
\item[$(b)$]
$\Lx _0$  is a maximal lift of  $(F \otimes _k T/\m_T^2, \delta)$.   
\end{itemize}
In particular, 
there always exists a universal lift of $\Fx$, and 
it is unique up to $k$-algebra automorphisms and chain isomorphisms. 
\end{thm}

\begin{pf}

[(b) $\Rightarrow$ (a)] : 
Let  $\Lx _0$ be a maximal lift of  $(F \otimes _k T/ \m_T^2, \delta)$. 
To simplify the notation we write  $\phi _{\Lx _0}$  as  $\phi$. 
We would like to prove that 
$$
\phi (A) \ :\  \Hom _{k\text{-alg}}(P_0, A) \to \F (A) 
$$   
is a bijection  for any  $A \in \A$. 
We prove this by induction on  the length of $A$. 
If  $\length (A) = 1$, then  $A = k$  and  $\phi (k)$  is clealy  bijective. 
 
[The surjectivity of $\phi (A)$]: 
Take a socle element  $\epsilon \in A$, and we have a small extension 
$$
\begin{CD}
0 @>>> k @>\epsilon>> A  @>\pi>> \overline{A} @>>> 0
\end{CD}
$$
where  $\overline{A} = A/(\epsilon)$. 
By the induction hypothesis,  $\phi (\overline{A})$  is bijective. 

$$
\begin{CD}
\Hom _{k\text{-alg}}(P_0 , A) @>{\phi(A)}>>  \F (A) \\
@V{\pi_*}VV  @V{\F (\pi)}VV  \\
\Hom _{k\text{-alg}}(P_0 , \overline{A}) @>{\phi(\overline{A})}>>  \F (\overline{A}) \\
\end{CD}
$$
To prove the surjectivity of  $\phi (A)$, let  $(F \otimes _k A, \Delta)$  be any element of  $\F (A)$. 
Since  $\phi (\overline{A})$  is surjective, there is $g \in \Hom _{k\text{-alg}}(P_0, \overline{A})$  such that 
$$
(F \otimes _k A, \Delta)\otimes _A \ {}_{\pi}\overline{A} \cong 
(F \ctensor P_0, \Delta _0) \otimes _{P_0} \ {}_g\overline{A}.
$$
Hence, it follows from Lemma \ref{lifting homo} that there is an isomorphism of graded  modules $\alpha : F \otimes _k A \to F \otimes _k A$  such that 
$$
(F \otimes _k A, \alpha \Delta \alpha ^{-1})\otimes _A \ {}_{\pi}\overline{A}  =  (F \ctensor P_0, \Delta _0) \otimes _{P_0} \ {}_g\overline{A}.
$$
Now taking  the fiber product 
$$
\begin{CD}
P @>>> P_0 \\
@VVV @VgVV  \\
A @>\pi>> \overline{A},   \\
\end{CD}
$$
we see from the above equality that the complex  
$(F \ctensor P_0, \Delta _0) \otimes _{P_0} \ {}_g\overline{A}$  is liftable to $A$, and  it follows from Proposition \ref{fiber prod} that 
the chain complex  $(F \ctensor P_0, \Delta _0)$  is liftable to  $P$.
If the extension $P \to P_0$  is  nontrivial, then it contradicts to that   $(F \ctensor P_0, \Delta _0)$  is a maximal lift. 
Hence $P \to P_0$ is a trivial small extension, and  $P \to P_0$  has a right inverse in  $\Ahat$.
In particular, the $k$-algebra map  $g : P_0 \to \overline{A}$  can be lifted to the $k$-algebra map $f : P_0 \to A$,  i.e. $\pi \cdot f = g$. 
Then, note that both $(F \otimes _k A, \Delta _0 \otimes _{P_0} {}_fA)$   and  $(F \otimes _k A, \alpha \Delta \alpha ^{-1})$  are lifts  of  $(F \otimes _k \overline{A}, \Delta _0 \otimes _{P_0} {}_g \overline{A})$. 
Hence, by Lemma \ref{difference}, we have 
$$
\Delta _0 \otimes _{P_0} {}_fA = \alpha \Delta \alpha ^{-1} + h \otimes \epsilon$$
for some chain homomorphism  $h : F \to F[-1]$ of degree $-1$.
Then we may write 
$$
[h] = \sum_{i=1}^r c_i [t_i^*] \qquad (c_i \in k), 
$$
as an element of  $\Ext _R^1(\Fx, \Fx)$. 
Now define a $k$-algebra map $\widetilde{\varphi} : T \to A$  by 
$$
\widetilde{\varphi} (t_i) = f(t_i) - c_i \epsilon \qquad (1 \leq i \leq r). 
$$
It can be easily verified that  $\widetilde{\varphi} (t_it_j) = f(t_it_j)$  for any $i, j$. 
Since  $I_0 \subseteq \m_T ^2$, we have  $\widetilde{\varphi} (I_0) = f(I_0) = 0$, thus  $\widetilde{\varphi}$ induces the $k$-algebra map $\varphi : P_0 \to A$ and $\varphi |_{\m_{P_0}^2} = f |_{\m_{P_0}^2}$. 
Then, by the choice of  $\varphi$, we see that 
$$
\Delta _0 \otimes _{P_0} {}_{\varphi} A = 
\Delta _0 \otimes _{P_0} {}_{f} A - \sum _{i=1} ^r c_i t_i^* \otimes \epsilon
= \alpha \Delta \alpha ^{-1} + ( h - \sum _{i=1} ^r c_i t_i^* ) \otimes \epsilon.
$$
Thus it follows from Lemma \ref{difference}(b) that 
$$
(F \ctensor P_0 , \Delta _0) \otimes _{P_0} {}_{\varphi} A 
\cong 
(F \otimes _k A, \alpha \Delta \alpha ^{-1})
\cong 
(F \otimes _k A, \Delta).  
$$
 This proves the surjectivity of   $\phi (A)$.

\vspace{4pt}
[The injectivity of $\phi (A)$]: 
Let  $\varphi _1 , \varphi _2$  be $k$-algebra homomorphisms $P_0 \to A$  with 
$(F \ctensor P_0, \Delta _0) \otimes _{P_0} {}_{\varphi_1} A 
\cong (F \ctensor P_0, \Delta _0) \otimes _{P_0}  {}_{\varphi_2} A$. 
We want to show  $\varphi _1 = \varphi _2$. 
Take a socle element $\epsilon \in A$ and we consider the small extension 
$$
\begin{CD}
0 @>>> k @>\epsilon>> A @>\pi>> \overline{A} @>>> 0.  \\
\end{CD}
$$
Then, by the induction hypothesis, we have  $\pi \cdot \varphi _1 = \pi \cdot \varphi _2$. 
Now consider the mapping $\psi = \varphi _1 - \varphi _2 : P_0 \to A$, and we see that the image of  $\psi$  is contained in  $k\epsilon$. 
Note that  $\psi (1) = 0$  and that 
$
\psi (xy) = \varphi _1(x)\psi (y) - \psi (x) \varphi _2(y) =0
$
if $x, y \in \m_{P_0}$, since  $\m_{P_0} \epsilon = \epsilon \m_{P_0} = 0$. 
Therefore  $\psi (\m_{P_0}^2) = 0$, and we have  $\varphi _1 |_{\m_{P_0}^2} = \varphi _2 |_{\m_{P_0}^2}$. 
 Since  $\varphi _1 = \varphi _2 + \psi$, it follows 
$$
\Delta _0 \otimes _{P_0} {}_{\varphi _1}A   
= \Delta _0 \otimes _{P_0} {}_{\varphi _2}A + \sum _{i=1} ^r t_i^* \otimes \psi (t_i). 
$$   
Denoting  $\psi (t_i) = c_i \epsilon$  with  $c _i \in k$  for $1 \leq i \leq r$, we have the equality
$$
\Delta _0 \otimes _{P_0} {}_{\varphi _1}A   
= \Delta _0 \otimes _{P_0} {}_{\varphi _2}A + \sum _{i=1} ^r c_i t_i^* \otimes \epsilon. 
$$   
Then it follows from Lemma \ref{difference}(b) that $\sum_{i=1}^r c_i [t_i^*] = 0$  as an element of  $\Ext_R^1(\Fx, \Fx)$. 
Since $\{ [t_1^*], \ldots,  [t_r^*]\}$  is a $k$-basis of $\Ext_R^1(\Fx, \Fx)$, we have  $c_i=0$ for all $i$, hence  $\psi = 0$.

\vspace{4pt}
[(a) $\Rightarrow$ (b)] : 
Suppose that  $\Lx _0 = (F \ctensor P_0, \Delta_0)$  is  a universal lift of  $\Fx$. 
Take any maximal lift $(F \ctensor T/I_1, \Delta _1)$ of  $(F \otimes _k T/\m_T ^2, \delta)$, and since it is a universal lift by the implication (b) $\Rightarrow$ (a),  Lemma \ref{unique universal lifting} forces $P_0$ to be isomorphic to $T/I_1$. 
Thus we may assume that  $P_0 = T/I_1$. 
Lemma \ref{Zorn} implies that we can take a maximal element $(I_2, \Delta _2)$  in  $\I (I_1, \Delta _0)$, which is in fact a maximal lift of  $(F \otimes _k T/\m_T ^2, \delta)$. 
Then, again by the implication  (b) $\Rightarrow$ (a),  $(F \ctensor T/I_2, \Delta _2)$ is also a universal lift of  $\Fx$, and hence  $T/I_2$  is isomorphic to $T/I_1$  by Lemma  \ref{unique universal lifting}. 
Since  $I_2 \subseteq  I_1$, the following lemma forces  $I_2 = I_1$. 
This implies that  $(F \ctensor P_0, \Delta _0 ) = (F \ctensor T/I_2, \Delta _2)$, which is a maximal lift as desired. 
\qed
\end{pf}

\begin{lemma}
Let  $T$  be a non-commutative formal power series ring, and let  $I_2 \subseteq I_1$  be closed ideals of  $T$. 
Suppose  $T/I_1$  is isomorphic to  $T/I_2$  as a $k$-algebra. 
Then  $I_1=I_2$. 
\end{lemma}

\begin{pf}
The isomorphism $T/I_1 \cong T/I_2$  induces isomorphisms 
 $T/I_1 + \m_T^n \cong T/I_2+ \m_T^n$  for any integer $n$. 
Since  $I_2 + \m_T^n \subseteq I_1+ \m_T^n$, comparing the lengths, we have 
the equality $I_2 + \m_T^n = I_1+ \m_T^n$  for each  $n$. 
Thus  $I _2 = \bigcap _n I_2 + \m_T^n = \bigcap _n  I_1+ \m_T^n = I_1$. 
\qed\end{pf}

\subsection{Every complete local algebra is a parameter algebra}

\begin{lemma}\label{lift of resolution}
Let  $A' \to A$  be a surjective morphism in  $\A$. 
Suppose the following two conditions hold. 
\begin{itemize}
\item[$(a)$]
$\Lx = (F\otimes _k A, \Delta)$  is a left $R \otimes _k A^{op}$-free resolution of a left $R \otimes _k A^{op}$-module $M$, and  $\Lx$ is a lift of a free complex  $\Fx = (F, d)$  over $R$. 
\item[$(b)$]
There is a left  $R \otimes _k A'$-module $M'$  such that  $M'$  is flat as a right $A'$-module and  $M' \otimes _{A'} A \cong M$  as left  $R \otimes _k A^{op}$-modules.
\end{itemize}
Then there is  a lifting chain complex  $\Lx ' = (F \otimes _k A', \Delta ')$  of  $\Lx$  that is a left  $R \otimes _k A'^{op}$-free resolution of  $M'$. 
\end{lemma}

\begin{pf}
We may write  $A = A'/I'$  where  $I'$  is an ideal of  $A'$. 
Take a set of generators $\{ x_{\lambda} | \ \lambda \in \Lambda \}$ of  $M$  as a left $R \otimes _k A^{op}$-module. 
 Since  $M \cong  M'/M'I'$  as $R \otimes _k A'^{op}$-modules, 
we can take a subset  $\{ x'_{\lambda}| \ \lambda \in \Lambda\}$ of  $M'$  that is an inverse image of  $\{x_{\lambda}\}$. 
Then the equality 
$M' = R \{ x'_{\lambda}\}A' + M'I'$  holds. 
Since  $I'$  is a nilpotent ideal, we have  $M' = R \{ x'_{\lambda}\}A'$. 

By this argument we can show that every surjective homomorphism  $F_0 \otimes _k A \to M$  of left  $R \otimes _k A^{op}$-modules can be lifted to a surjective homomorphism  $F_0 \otimes _k A' \to M'$  of  left $R \otimes _k A'^{op}$-modules. 
Now take the kernels of these surjective maps and we have exact sequences 
$$
\begin{CD}
0 @>>> M_1 @>>> F_0 \otimes _k A @>>> M @>>> 0, \\
0 @>>> M'_1 @>>> F_0 \otimes _k A' @>>> M' @>>> 0.  \\
\end{CD}
$$
Notice that  $M' _1$  is flat as a right $A'$-module, since  $F_0 \otimes _k A'$  and  $M'$  are flat. 
Thus the isomorphism   $M' \otimes _{A'} A \cong M$  implies  $M'_1 \otimes _{A'} A \cong M_1$. 
Then by the same manner as above, the surjective homomorphism  $F_1 \otimes _k A \to M_1$  is liftable to a surjective homomorphism  $F_1 \otimes _k A' \to M'_1$. 
In such a way, by induction, we can construct a chain complex with the underlying graded module  $F \otimes _k A'$, which is a lift of  $(F \otimes _k A, \Delta )$. \qed\end{pf}

\begin{thm}\label{resol of k}
Let   $R$  be a complete local $k$-algebra, i.e.  $R \in \Ahat$, and let  $\Fx = (F, d)$ be a left $R$-free resolution of the residue field  $k = R/\m_R$. 
Then there is a universal lift of  $\Fx$ that is an acyclic complex  of the form  $(F \ctensor R , \Delta)$  with the homology  $H_0 (F \ctensor R, \Delta) = R$. 
In particular, $R$  is the parameter algebra of the universal lift of  $\Fx$. 
\end{thm}

\begin{pf}
Note that the obvious exact sequence 
$0 \to \m_R \to R \to k \to  0$ implies 
$$
\Ext _R ^1(\Fx, \Fx) = \Ext _R ^1(k, k) \cong \Hom _k ( \m_R/\m_R^2, k).
$$
Thus if we denote $R = T / I$  where  $T = \formal$  and  $I \subseteq \m _T ^2$, we can take the dual bais  $\{ t_1^*, \ldots , t_r^*\}$ as a basis of  $\Ext _R ^1(\Fx, \Fx)$. 
Set  $\delta  =  d \otimes 1 + \sum _{i=1} ^r t_i^* \otimes t_i$, and 
we see that   $(F \otimes _k T/\m_T ^2, \ \delta)$  is a lift of  $\Fx$  as in the beginning of Section 3.2. 
 By the exact sequence of chain complexes 
$$
\begin{CD}
0 @>>> (F \otimes _k \m_T /\m_T^2, \ d \otimes 1) @>>> (F \otimes _k T/\m_T^2, \ \delta) @>>> (F, d) @>>> 0,  \\
\end{CD}
$$
and by the acyclicity of  $(F, d)$, we easily see that  $(F \otimes _k T/\m_T^2, \ \delta)$  is acyclic as well, and  $H_0 (F \otimes _	k T/\m_T^2, \ \delta) \cong  T/\m_T ^2 = R / \m_R ^2$  as an $R \otimes _k  (R/\m_R^2)^{op}$-module. 

Starting from  $\Delta _2 = \delta$, and using Lemma \ref{lift of resolution},  we can inductively construct a sequence of chain complexes $(F \otimes _k R/ \m _R^n, \ \Delta _n )$  for  $n \geq 2$  satisfying the equalities  
 $(F \otimes _k R/ \m _R^{n+1}, \ \Delta _{n+1} ) \otimes _{R/\m_R^{n+1}}R/\m_R^{n} = (F \otimes _k R/ \m _R^n, \ \Delta _n )$  and  
 $H_0 (F \otimes _k R/ \m _R^{n}, \ \Delta _{n} ) = R/\m_R^{n}$. 
Now let $\Delta = \varprojlim \Delta _n$,  and we have a lifting chain complex   $(F \ctensor R, \ \Delta)$  of  $(F\otimes _k T/\m_T^2, \ \delta)$. 
   
First, we claim that  $(F \ctensor R, \ \Delta)$ is acyclic and  $H_0 (F \ctensor R, \ \Delta) \cong R$  as  an $R \ctensor R^{op}$-module. 
For this, let  $\Omega _n^i$  be the kernel of $\Delta _{n, i} : F _i \otimes _k R/\m_R^n \to F_{i-1} \otimes _k R/\m_R^{n}$ for any $n, i \geq 0$  where we understand that  $\Omega _n^0 = R/\m_R^n$. 
  By the proof of Lemma \ref{lift of resolution}, we have a commutative diagram with exact rows
$$
\begin{CD}
0 @>>> \Omega _{n+1}^i @>>> F_i \otimes _k R/\m_R^{n+1} @>>> \Omega _{n+1} ^{i-1} @>>> 0 \\
@. @VVV @VVV @VVV \\
0 @>>> \Omega _n^i @>>> F_i \otimes _k R/\m_R^{n} @>>> \Omega _n ^{i-1} @>>> 0, \\
\end{CD}
$$
where the vertical arrows are surjective. 
This implies the exact sequence 
$$
\begin{CD}
0 @>>> \varprojlim \Omega _{n}^i @>>> F_i \ctensor R @>>> \varprojlim  \Omega _{n} ^{i-1} @>>> 0,  \\
\end{CD}
$$
and hence the complex  $(F \ctensor R , \ \Delta)$  is acyclic and  $H_0 (F \ctensor R , \ \Delta ) = \varprojlim R/\m_R ^n = R$  as desired. 

Now we prove that  $(F \ctensor R, \ \Delta )$ is a maximal lift of  $(F \otimes _k R/\m_R^2, \ \delta)$. 
Suppose that it is not a maximal lift. 
Then there will be a nontrivial small extension  
$$
\begin{CD}
0 @>>> k  @>{\epsilon}>>  R'  @>p>> R @>>> 0,  \\
\end{CD}
$$  
of  $R$  so that  the complex  $(F \ctensor R, \ \Delta )$  is liftable to a chain complex   $(F \ctensor R', \ \Delta ')$. 
The exact sequence 
$$
\begin{CD}
0 @>>> (F ,d ) @>{1 \otimes \epsilon}>>  (F \ctensor R', \ \Delta ') @>>> (F \ctensor R, \ \Delta ) @>>> 0 \\
\end{CD}
$$
forces  $(F \ctensor R', \ \Delta ')$  to be acyclic as well,  
and taking the homologies we have an exact sequence 
$$
\begin{CD}
0 @>>> k  @>{\epsilon}>>  H'_ 0 @>\pi>> R @>>> 0,  \\
\end{CD}
$$
where  $H'_0 = H_0 (F \ctensor R', \ \Delta ')$  and  $\pi$  is a homomorphism of  $R \ctensor {R'}^{op}$-modules. 
Take  an element  $x_0 \in H'_0$  with  $\pi (x_0) = 1$, and obviously we have 
$H'_0 = x_0 R' + H'_0 \epsilon$. 
Since  $\epsilon ^2 =0$, it follows that $H'_0 = x_0 R'$. 

We claim that  $H'_0$  is a free module as a right $R'$-module. 
To prove this, assume  $x_0a' =0$  for $a' \in R'$, and we want to show  $a' =0$. 
Suppose $a'\not= 0$. 
Since  $0 = \pi (x_0 a') = 1_R \cdot a' = p(a')$, we see  $a' \in  \epsilon k$, hence $a' = \epsilon c$  for some $c \not= 0 \in k$. 
Then we have  $x_0 \epsilon = 0$, and the right $R'$-module   $H'_0 = x_0 R'$  is in fact a right $R$-module.
Hence  $0 \to k \to H'_0 \to R \to 0$  is an exact sequence of right $R$-modules. 
Therefore the sequence splits and  $H'_0 \cong k \oplus R$  as a right $R$- (hence $R'$-)module. 
This contradicts that  $H'_0$  is generated by a single element as a right $R'$-module.

Now we have shown  $H'_0$  is a free right $R'$-module. 
Since $H'_0$  is a left $R$-module as well, 
for any  $a \in R$, we find a unique element  $a' \in R'$  with  $a\cdot x_0 = x_0 \cdot a'$. 
Now define a map  $f : R \to R'$  by  $f(a) = a'$. 
Since  $(ab)x_0=a (x_0 f(b))= (ax_0)f(b)= x_0(f(a)f(b))$, we can see that  $f$  is a $k$-algebra homomorphism. 
Since we have an equality 
$a = \pi (ax_0) =\pi (x_0f(a)) = 1_R \cdot f(a) = p(f(a))$ for any $a \in R$, we see  $p \cdot f = 1$. 
This contradicts that  $(R', \epsilon)$  is a nontrivial extension, and the proof is completed. 
\qed\end{pf}

\begin{rem}
If  $R$  is left noetherian, then we can take as each  $F_i$   a finitely generated left  $R$-free module, and  $F_i \ctensor R$  is a left $R \ctensor R^{op}$-free module. 
 In this case  the acyclic complex  $(F \ctensor R, \Delta )$  in the theorem   is a free resolution of  $R$  as a left $R \ctensor R^{op}$-module. 
However, in general, notice that  $(F \ctensor R, \Delta)$  may not be a free complex of  left  $R \ctensor R^{op}$-modules. 
\end{rem}

\subsection{Deformation of modules}

Let  $R$  be a $k$-algebra as before. 
In this section, we consider the case where the complex  $\Fx = (F, d)$  is a    free resolution of a left $R$-module $M$. 
Of course, in this setting, we have the equality 
$\Ext_R^1(M, M) = \Ext _R ^1(\Fx, \Fx)$ and it is assumed to be of finite dimension as before. 

For any  $A \in \A$, a left $R \otimes _k A^{op}$-module  $X$  is said to be a {\bf flat deformation of} $M$ {\bf along} $A$  if  $X$  is a flat module as a right $A$-module, and there is an isomorphism  $X \otimes _A k \cong M$  as left $R$-modules. 
And two flat deformations  $X$  and  $Y$  of  $M$  along $A$  are said to be isomorphic  if they are isomorphic as $R \otimes _kA^{op}$-modules. 
We consider the following two functors $\A \to (Sets)$: 
$$
\begin{array}{ll}
\F (A) &= \text{the set of isomorphism classes of lifting chain complexes of} \  \Fx, \\ 
\F _M (A) &= \text{the set of isomorphism classes of flat deformations of } \ M. \\
\end{array}
$$ 

\begin{thm}\label{deformation}
We have an isomorphism  $\F \cong \F _M$  as functors on  $\A$. 
In particular, the functor  $\F_M$  is pro-representable as is  $\F$, i.e. 
there is an isomorphism  $\F _M \cong \Hom _{k\text{-alg}}(P_0, \ \ )$ of functors on  $\A$,  where  $P_0$  is the parameter algebra of the universal lift of  $\Fx$. 
\end{thm}

\begin{pf}
Let  $A \in \A$  and  let  $(F \otimes _k A, \Delta)$  be a lift of  $\Fx$  to  $A$. 
Notice that, by induction on  the length of  $A$,   one can easily prove that   $(F \otimes _k A, \Delta)$  is acyclic, as  $\Fx$  is acyclic. 
Therefore it gives a left $R \otimes _k A^{op}$-free resolution of  $H_0(F \otimes _k A, \Delta)$. 
Since  $(F \otimes _k A, \Delta) \otimes _A k = \Fx$ is acyclic, we have 
$\Tor _i ^{A^{op}}(H_0(F \otimes _k A, \Delta), k)=0$  for $i >0$. 
Since  $A$  is an artinian local algebra, this implies  that $H_0(F \otimes _k A, \Delta)$  is  flat as a right $A$-module.
We should note that  $H_0(F \otimes _k A, \Delta) \otimes _A k = H_0(\Fx) = M$.  Hence  $H_0(F \otimes _k A, \Delta)$  is a flat deformation of  $M$  along  $A$. 
Thus we obtain a well-defined mapping  $ H_0 : \F (A)  \to \F _M (A)$   by taking the $0$-th homology. 
   It is trivial that the map is injective, and Lemma \ref{lift of resolution}  implies it is surjective.
\qed\end{pf}

\newpage
\section{Properties of parameter algebras}

\subsection{Obstruction maps} 
As before  $\Fx = (F, d)$  is a projective complex over a $k$-algebra  $R$. 
Let  $P \in \Ahat$  and let  $\L = (F \ctensor P, \ \Delta )$  be a lift of   $\Fx$ to  $P$, which we fix in this section. 
We aim at constructing the obstruction map 
$$
\alpha _{\L} : \T (P) \to  \Ext _R ^2(\Fx, \Fx),  
$$
which will enable us to compare the cohomology modules between  $P$  and  $R$.

\begin{onepara}\label{def of alpha}
[To define $\alpha _{\L}$]:

Now, suppose we are given a class of small extension  $[P', \epsilon] \in  \T(P)$. 
Lemma \ref{lifting homo complete1}  forces  $\Delta : F \ctensor P \to F \ctensor P[-1]$ to be lifted to  $\Delta ' : F \ctensor P' \to F \ctensor P'[-1]$. 
Note that  $\Delta'$ is just a lift as a graded homomorphism, and it may not  holds that  ${\Delta '} ^2 =0$. 
Recall from Lemma \ref{exact seq} that we then have a commutative diagram of graded left $R \ctensor P'^{op}$-modules  with exact rows   
$$
\begin{CD}
0 @>>> F @>{1\otimes \epsilon}>> F \ctensor P' @>{1 \otimes \pi}>> F \ctensor P @>>> 0 \\
@. @VdVV  @V{\Delta'}VV  @V{\Delta}VV  \\
0 @>>> F [-1]@>{1\otimes \epsilon}>> F \ctensor P'[-1] @>{1 \otimes \pi}>> F \ctensor P[-1] @>>> 0.  \\
\end{CD}
$$
Since  $d^2 =0$  and  $\Delta ^2 =0$, we have the following commutative diagram.
$$
\begin{CD}
0 @>>> F @>{1\otimes \epsilon}>> F \ctensor P' @>{1 \otimes \pi}>> F \ctensor P @>>> 0 \\
@. @V0VV  @V{{\Delta'}^2}VV  @V{0}VV  \\
0 @>>> F [-2]@>{1\otimes \epsilon}>> F \ctensor P'[-2] @>{1 \otimes \pi}>> F \ctensor P[-2] @>>> 0  \\
\end{CD}
$$
By chasing the diagram, we see that 
there is a graded left $R$-module homomorphism  $\sigma : F \to F [-2]$  with 
 ${\Delta '}^2  = \sigma \otimes \epsilon$, i.e. 
${\Delta '}^2  (\sum x \otimes a) = \sum \sigma (x) \otimes \epsilon a$  
for  any formal infinite sum  $\sum x \otimes a \in F \ctensor P'$.  

First we claim that 

\vspace{3pt}
\noindent (i)  $\sigma : F \to F[-2]$  is a chain map. 
\vspace{3pt}

In fact, it holds that 
$$
\sigma \cdot d \otimes \epsilon = (\sigma \otimes \epsilon) \Delta ' 
= {\Delta '} ^3 = \Delta ' (\sigma \otimes \epsilon) = d\cdot \sigma \otimes \epsilon,  
$$
hence it follows that  $[d, \sigma] = d\cdot \sigma - \sigma \cdot d = 0$. 
Thus the graded homomorphism  $\sigma$  is a chain map of degree $-2$, therefore it defines an element  $[\sigma ] \in \Ext _R ^2(\Fx, \Fx)$. 

Next we claim that 

\vspace{3pt}
\noindent (ii)  the class $[\sigma]$  does not depend on a choice  of a lifting map  $\Delta '$. 
\vspace{3pt}

In fact, if  $\Delta '$ and  $\Delta ''$  are two lifting maps of  $\Delta$, then we have a commutative diagram 
$$
\begin{CD}
0 @>>> F @>{1\otimes \epsilon}>> F \ctensor P' @>{1 \otimes \pi}>> F \ctensor P @>>> 0 \\
@. @V0VV  @V{{\Delta'}- \Delta ''}VV  @V{0}VV  \\
0 @>>> F [-1]@>{1\otimes \epsilon}>> F \ctensor P'[-1] @>{1 \otimes \pi}>> F \ctensor P[-1] @>>> 0,   \\
\end{CD}
$$
  from which we can see the existence of graded homomorphism  $\tau : F \to F [-1]$  with  $\Delta ' - \Delta '' = \tau \otimes \epsilon$. 
Then, we have an equality 
$$
{\Delta '}^2 = (\Delta '' + \tau \otimes \epsilon ) ^2 = {\Delta ''}^2  + 
(d\tau + \tau d ) \otimes \epsilon,  
$$
hence, setting  ${\Delta '}^2= \sigma ' \otimes \epsilon$ and  ${\Delta ''}^2= \sigma '' \otimes \epsilon$, we have  $\sigma ' - \sigma '' = d\tau + \tau d$, i.e. $[\sigma '] = [\sigma '']$.
\end{onepara}

Now we can define a mapping  
$$
\alpha _{\L} : \T (P) \to \Ext _R ^2(\Fx, \Fx)
$$
by  sending  $[P', \epsilon]$  to  the class $[\sigma]$.

By (i) and (ii) above, $\alpha _{\L}$  is a well-defined mapping. 
Furthermore, we can show the following. 

\begin{lemma}\label{k-lin of alpha}
The mapping  $\alpha _{\L} : \T (P) \to \Ext _R ^2 (\Fx,\Fx)$  is $k$-linear. 
\end{lemma}

\begin{pf}
To prove the equality 
$$
\alpha _{\L} (c_1 [P_1, \epsilon _1] + c_2 [P_2, \epsilon _2]) 
= c_1  \alpha _{\L} ([P_1, \epsilon _1]) + c_2  \alpha _{\L} ([P_2, \epsilon _2]) $$
for  $c_i \in k$  and  $[P_i ,\epsilon _i] \in \T (P)\  (i=1,2)$,  
 let us consider the pull-back diagram 
$$
\begin{CD}
F \ctensor (P_1 \times _P P_2)  @>>> F \ctensor P_1 \\
@VVV                @VVV   \\
F \ctensor P_2  @>>> F \ctensor P,   \\
\end{CD}
$$
and take lifting graded homomorphisms $\Delta _i : F \ctensor P_i \to F \ctensor P_i [-1]$  of  $\Delta$  for  $i =1,2$. 
We may assume that  $c _i \not= 0$ for $i =1, 2$. 
Since there is a commutative diagram with exact rows by Lemma \ref{exact seq} ; 
$$
{\small
\begin{CD}
0 @>>> F \ctensor (P_1 \times _P P_2) @>>> (F \ctensor P_1) \oplus  (F \ctensor P_2) @>>> F \ctensor P @>>> 0 \\
@. @. @V{(\Delta _1, \Delta _2)}VV  @V{\Delta}VV \\ 
0 @>>> F \ctensor (P_1 \times _P P_2)[-1] @>>> F \ctensor P_1[-1] \oplus  F \ctensor P_2[-1] @>>> F \ctensor P [-1] @>>> 0,  \\
\end{CD}
}
$$
there is a naturally induced mapping $\widetilde{\Delta} : F \ctensor (P_1 \times _P P_2) \to  F \ctensor (P_1 \times _P P_2) [-1]$  which is a lifting map of both of  $\Delta _1$  and  $\Delta _2$. 
Let us take a chain homomorphism  $\sigma _i$  of  $\Fx$  so that  $\Delta _i ^2 = \sigma _i \otimes \epsilon _i$  for  $i = 1, 2$. 
Then it can be seen that 
$
\widetilde{\Delta}^2 = \sigma _1 \otimes (\epsilon _1, 0) + \sigma _2 \otimes (0, \epsilon _2). 
$ 
Recalling the definition of the sum $[Q, \epsilon] = c_1 [P_1, \epsilon _1] + c_2 [P_2 , \epsilon _2]$, we have  
$Q = P_1 \times _P P_2 /(c_1^{-1} \epsilon _1, -c_2^{-1}\epsilon _2)$  and  $\epsilon$  is the class of  $(c_1^{-1} \epsilon _1, 0)$. 
Thus, setting  $\Delta ' = \widetilde{\Delta}\otimes _{P_1 \times _P P_2} Q$,    we have the mapping  $\Delta ' : F \ctensor Q \to F \ctensor Q [-1]$  which is a lifting map of  $\Delta$  to  $Q$,   
and we easily see that 
$
{\Delta '} ^2 =  ( c_1 \sigma _1 + c_2 \sigma _2 ) \otimes \epsilon.  
$
Consequently, we have 
$
\alpha _{\L} ([Q, \epsilon]) =  c_1 [\sigma _1] + c_2 [\sigma _2]
=  c_1 \alpha _{\L} ([P_1, \epsilon _1]) + c_2 \alpha _{\L} ([P_2, \epsilon _2]). 
$
\qed\end{pf}

\begin{lemma}\label{comm of alpha}
Let  $f : P_1 \to P_2$ be a $k$-algebra map in  $\Ahat$, and let  $\L_2 = (F \ctensor P_2, \Delta _2)$  be a lift of  $\Fx$  to  $P_2$. 
Suppose there exists a lift  $\L _1 = (F \ctensor P_1, \Delta _1)$  of  $\L _2$ to  $P_1$, i.e.  
$\L _1 \otimes _{P_1} {}_{f}P_2 = \L _2$. 
Then there is a commutative diagram 
$$
\begin{CD}
\T(P_2) @. \begin{picture}(0,0)
\put(-35,3){\vector(1,0){60}}
\put(-10,-10){\small${f^*}$}
\end{picture} @.  \T (P_1) \\
\begin{picture}(0,0)
\put(25,15){\vector(1,-1){25}}
\put(5,0){\small${\alpha _{\L_2}}$}
\end{picture} @.      
\begin{picture}(0,0)
\put(30,15){\vector(-1,-1){25}}
\put(30,0){\small${\alpha _{\L_1}}$}
\end{picture} 
\\
@. \hphantom{FF}\Ext _R ^2(\Fx, \Fx).\hphantom{FF} @. \\
\end{CD}
$$
\end{lemma}

\begin{pf}
Set  $[P_1' , \epsilon _1] = f^* ([P_2', \epsilon _2])$  for  $[P_2', \epsilon _2] \in \T (P_2)$. 
From the definition, there is a pull-back diagram
$$
\begin{CD}
0 @>>> k @>{\epsilon _1}>> P_1' @>>> P _1  @>>> 0 \\
@.  @| @V{f'}VV  @V{f}VV  \\
0 @>>> k @>{\epsilon _2}>>  P'_2 @>>> P_2 @>>> 0.  \\
\end{CD}
$$
Let  $\Delta _1'$  be any lift  of  $\Delta _1$  onto  $F \ctensor P_1'$. 
Then  $\Delta _2' := \Delta _1' \otimes _{P_1'} {}_{f'} P_2'$  is a lift of  $\Delta _2 = \Delta _1 \otimes _{P_1} {}_f P_2$  onto  $F \ctensor P_2'$. 
Now write 
${\Delta '}_1^2 = \sigma \otimes \epsilon _1$ so that 
$\alpha _{\L_1} ([P_1', \epsilon _1]) = [\sigma]$. 
Then, 
we have  ${\Delta '}_2^2= {\Delta '}_1^2 \otimes _{P_1'} {}_{f'}P_2' = \sigma \otimes \epsilon _2$, hence 
$\alpha _{\L_2} ([P_2', \epsilon]) = [\sigma]$. 
This shows that  $\alpha _{\L_2} = \alpha _{\L_1}\cdot f^*$. 
\qed\end{pf}

\begin{thm}\label{injection}
Let  $\Lx _0 = (F \ctensor P_0, \Delta _0)$  be the universal lift of  $\Fx = (F, d)$  with parameter algebra  $P_0$. 
Then, the $k$-linear mapping  
$$
\alpha _{\Lx _0} : \T (P_0) \to \Ext _R ^2(\Fx, \Fx)
$$
is an injection. 
\end{thm}

\begin{pf}
Let  $[P_1, \epsilon _1] \in \T (P_0)$  be a nontrivial small extension. 
We only have to show  $\alpha _{\Lx_0} ([P_1, \epsilon _1]) \not= 0$. 
Suppose  $\alpha _{\Lx_0} ([P_1, \epsilon _1]) = 0$. 
Then, for any lifting map  $\Delta _1 : F\ctensor P_1 \to F \ctensor P_1 [-1]$  of  $\Delta _0$  to $P_1$, we have 
 $\Delta _1 ^2 = \sigma \otimes \epsilon _1$, where  $\sigma = [d, h]$  for some graded homomorphism  $h : F \to F [-1]$. 
Now putting  $\Delta _1' = \Delta _1 - h \otimes \epsilon _1$, 
one can see that 
${\Delta _1'}^2 = \Delta _1^2 - [d,h]\otimes \epsilon _1= 0$. 
Therefore, $(F \ctensor P_1, \Delta _1')$  is a lifting chain complex of  
$(F \ctensor P_0, \Delta _0)$. 
This is a contradiction, because $(F \ctensor P_0, \Delta _0)$  is a maximal lift. 
See Theorem \ref{univ=max}. 
\qed\end{pf}

\begin{cor}
Suppose  $\Ext _R ^2(\Fx, \Fx) =0$. 
Then the parameter algebra  $P_0$  of the universal lift of  $\Fx$  is isomorphic to the non-commutative formal power series ring  $\formal$. 
\end{cor}

\begin{pf}
Under the assumption, we have  $\T (P_0) =0$  by Theorem \ref{injection}. 
Therefore if we describe  $P_0 = T/I_0$  where  $T$ is a formal power series ring and  $I \subseteq  \m _T ^2$  is a closed ideal, then Proposition \ref{con iso T} forces  $\Hom _{con} (I, k) =0$.   
Thus we only have to show the following lemma.
\qed\end{pf}

\begin{lemma}
Let  $I$  be a closed ideal of a non-commutative formal power series ring $T$. 
If  $\Hom _{con} (I, k) =0$, then  $I =0$.    
\end{lemma}

\begin{pf}
Suppose  $I \not= 0$. 
Then, by Corollary \ref{cor nakayama}, we have  $I \not= \overline{\m _T I + I \m_T}$. 
Therefore,   $I \not= \m _T I + I \m _T + (\m_T^n \cap I)$  for a large integer $n$. 
Since  $\Hom _{con} (I, k)$  contains every $k$-linear map $I/\m _T I + I \m _T + (\m_T^n \cap I) \to k$, we have  $\Hom _{con}(I, k) \not= 0$. 
\qed
\end{pf}

This corollary can be generalized to the following theorem.

\begin{thm}\label{analytic gen}
Let  $P_0 =T/I_0$  be the parameter algebra of the universal lift of  $\Fx = (F, d)$, where  $T$  is a non-commutative formal power series ring and  $I _0 \subseteq \m _T ^2$  is a closed ideal. 
Suppose  $\ell = \Kdim _k \Ext _R ^2 (\Fx, \Fx)$  is finite.  
Then, the ideal  $I_0$  is analytically generated by at most $\ell$ elements. 
\end{thm}

\begin{pf}
Combining Theorems \ref{con iso T}  and  \ref{injection}, we have an injective $k$-linear map  \\
$\Hom _{con} (I_0, k)     \to   \Ext _R ^2 (\Fx, \Fx)$. 
In particular, $\Hom _{con} (I_0, k)$ is a $k$-vector space of finite dimension. Since we have the equality  
$
\Hom _{con} (I_0, k) = \bigcup _{n=1} ^{\infty} \Hom _k ( I_0/\m_TI_0 + I_0 \m_T + ( \m_T^n \cap I), \ k)    
$
by definition, there is an integer  $n_0$  such that 
$
\Hom _{con} (I_0, k) = \Hom _k ( I_0/\m_TI_0 + I_0 \m_T + ( \m_T^n \cap I), \ k)  $
for $n \geq n_0$. 
Hence we have the equalities 
$$
\m_TI_0 + I_0 \m_T + ( \m_T^n \cap I) = \m_TI_0 + I_0 \m_T + ( \m_T^{n+1}  \cap I) = \cdots = \overline{\m_TI_0 + I_0 \m_T}.
$$ 
Thus it follows that  $\Hom _{con} (I_0, k) = \Hom _k (I_0/\overline{\m_TI_0 + I_0 \m_T}, \ k)$. 
Since this is of dimension at most $\ell$, we have 
$\Kdim _k (I_0/ \overline{\m_TI_0 + I_0 \m_T}) \leq \ell$. 
Therefore, by virtue of Proposition \ref{nakayama},   $I_0$  is analytically generated by at most $\ell$ elements. 
\qed\end{pf}

\subsection{Universal lifts based on commutative algebras}

\begin{rem}\label{commutator}
Let  $T = \formal$  be a non-commutative formal power series ring. 
We denote by $C$  the commutator ideal which is a two-sided ideal generated by the commutators  $t_it_j -t_j t_i \ (1 \leq i < j \leq r)$, i.e. 
$$
C = \left( \{ t_it_j -t_j t_i \ | \ 1 \leq i < j \leq r\} \right). 
$$
Note that  $C$  is not a closed ideal if  $r \geq 2$. 
It is however easy to see that  $T/\overline{C}$  is isomorphic to the commutative formal power series ring  $k[[t_1, \ldots , t_r]]$. 
\end{rem}

\begin{rem}\label{ideal containing C}
Let  $I$  be an ideal of  $T$  that contains  $\overline{C}$. 
Then, $I$ is a closed ideal and there are a finite number of elements  $f_1, \ldots , f_{\ell} \in I$  with the equality   $I = (f_1, \ldots , f_{\ell}) + \overline{C}$. 

In fact, it is well known that any ideal of  $T/\overline{C} = k[[t_1, \ldots , t_r]]$  is closed. 
Since the natural projection  $\pi : T \to T/\overline{C}$  is continuous, 
$I = \pi ^{-1} (I/\overline{C})$ is closed as well. 
Since  $T/\overline{C}$  is noetherian, we can find finite elements  $f_1, \ldots , f_{\ell}$  which generate the ideal $I/\overline{C}$. 
Then we have the equality   $I =  (f_1, \ldots , f_{\ell}) + \overline{C}$. 
\end{rem}

Recall from  Section 3.3  that  $\F : \A \to (Sets)$  is a covariant functor such that  $\F (A)$  is the set of isomorphism classes of lifting chain complexes of  $\Fx$  to A, for any  $A \in \A$. 
We consider here the restriction of  $\F$  to commutative artinian algebras. 
For this end, we denote by  $\C$  the category of commutative artinian local $k$-algebras  $A$  with  $A/\m_A \cong k$  and  $k$-algebra homomorphisms. 
Note that  $\C$  is a full subcategory of  $\A$.

\begin{thm}
Let  $\Lx _0 = (F \ctensor P_0, \Delta _0)$  be the universal lift of  $\Fx$  with parameter algebra  $P_0 = T/I_0$  where  $T$  is a non-commutative formal power series ring and  $I \subseteq \m_T^2$.  
We set  $Q_0 = T/I_0 + \overline{C}$  which is a commutative noetherian complete local $k$-algebra. 
Then, the restricted functor  $\F |_{\C} : \C \to (Sets)$  is pro-represented by  $Q_0$, i.e. 
there is an isomorphism  $\F |_{\C} \cong \Hom _{k\text{-alg}} (Q_0, \ \ \ )$  as functors on  $\C$. 
\end{thm}

\begin{pf}
Note that if $A \in \C$, then $\Hom _{k\text{-alg}} (P_0, A) = \Hom _{k\text{-alg}} (Q_0, A)$. 
The theorem follows from this observation. 
\qed\end{pf}

\begin{defn}
We call  $Q_0$  in the theorem a {\bf commutative parameter algebra of the universal lift of} $\Fx$. 
And we call  $\Lx _0 \otimes _{P_0} Q_0$  the universal lift of  $\Fx$  based on commutative parameter algebra. 
\end{defn}

\begin{rem}
If  $\Fx$  is  a projective resolution of a left $R$-module $M$, 
then the universal lift of  $\Fx$  based on commutative parameter algebra  is nothing but the universal deformation of  $M$  whose existence is mentioned in  Theorem \ref{schlessinger}. 
(See also Proposition \ref{deformation}.)
\end{rem}

The commutative parameter algebra  $Q_0$  is of the form 
$k[[t_1, \ldots , t_r]]/\ga$  where  $\ga = I_0+\overline{C}/\overline{C} \subseteq T/\overline{C}= k[[t_1, \ldots ,t_r]]$.

\begin{prop}
Let  $Q_0 = k[[t_1, \ldots , t_r]]/\ga$  be a commutative parameter algebra of the universal lift of  $\Fx$. 
Then, the minimal number of generators of  $\ga$  is at most 
$\Kdim _k \Ext _R ^2(\Fx, \Fx)$. 
\end{prop}

\begin{pf}
It suffices to argue when  $\ell = \Ext_R^2 (\Fx, \Fx)$ is finite. 
Then, by Theorem \ref{analytic gen}, there are $f_1, \ldots , f_{\ell} \in I_0$ satisfying the equality  $I_0 = \overline{(f_1, \ldots , f_{\ell})}$. 
Thus, by virtue of Remark \ref{ideal containing C},  we have the equalities 
$$
I_0 + \overline{C} = 
\overline{I_0 + C} = 
\overline{C + (f_1, \ldots ,f_{\ell})}= 
\overline{C} +  (f_1, \ldots ,f_{\ell}). 
$$
Therefore  $\ga = \overline{I_0+C}/\overline{C}$  is generated by the images of  $f_1, \ldots , f_{\ell}$  in  $T/\overline{C}=k[[t_1,\ldots , t_r]]$. 
\qed\end{pf}

\subsection{Yoneda products} 

Let  $\Fx = (F,d)$  be a projective complex over $R$ with  $r = \Kdim _k \Ext _R^1(\Fx, \Fx)$  being finite as before. 
Then, as in the beginning of Section 3.2, we may consider the lifting chain complex  $\Lx = (F \ctensor T/\m_T^2 , \delta)$ 
with  $\delta = d \otimes 1 + \sum _{i=1}^r t_i^* \otimes t_i$, 
where  $T = \formal$  is a non-commutative formal power series ring and  $t_1^*,\ldots, t_r^*$  are chain homomorphisms which form a $k$-basis of  $\Ext _R ^1(\Fx, \Fx)$. 
Since  $\Lx$  is a lift of $\Fx$, we have the  $k$-linear map  
$$
\alpha _{\Lx} : \T (T/\m_T^2) \to \Ext _R ^2(\Fx, \Fx).
$$ 
by \ref{def of alpha}. 
See also Lemma \ref{k-lin of alpha}.
 
Note that  $\Ext _R ^{\cdot}(\Fx, \Fx) = \bigoplus _{i=-\infty}^{\infty} \Ext _R^i(\Fx, \Fx)$ is an algebra,  called Yoneda algebra,  whose  multiplication is given by Yoneda product. 
In fact, if  $f: F \to F[-i]$  and  $g: F \to F[-j]$  are chain homomorphisms of degree $-i$ and $-j$ respectively, then the composite  $f\cdot g : F \to F[-i-j]$ is a chain homomorphism of degree $-i-j$, and the product in  $\Ext _R ^{\cdot}(\Fx, \Fx)$ is given by $[f][g] = [f\cdot g]$. 
In the following lemma,  $\Ext _R^1(\Fx, \Fx)^2$ denotes the $k$-subspace of  $\Ext _R^2(\Fx, \Fx)$  generated by all the products of two elements in  $\Ext _R ^1 (\Fx, \Fx)$.

\begin{lemma}\label{image of T/m^2}
Under the circumstances above, 
the image of  $\alpha _{\Lx}$  is exactly  $\Ext _R ^1(\Fx, \Fx)^2$, i.e., 
$\alpha _{\Lx}(\T (T/\m_T^2)) = \Ext _R^1(\Fx, \Fx)^2$. 
\end{lemma}

\begin{pf}
Recall from Proposition \ref{con iso T} that 
there is an isomorphism of $k$-vector spaces 
$\tau : \Hom _k (\m_T^2/\m _T^3, k) = \Hom _{con}(\m_T^2, k) \to \T (T/\m_T^2)$.  Suppose  $\tau (f) = [T/I,  \epsilon]$  for  $f \not= 0 \in \Hom _k (\m_T^2/\m _T^3, k)$. 
Then, by definition of $\tau$, we have  $f (\epsilon) =1$  and $I$ is the kernel of the composition mapping  $\m_T^2 \to \m _T^2/\m_T^3$  with  $f: \m _T^2/\m_T^3 \to k$. 
Note, in this case, that  $T/I$  has $\{ 1, \overline{t_1}, \ldots , \overline{t_r}, \epsilon  \}$ as a $k$-basis, where  $\overline{t_i}$  denotes the image of $t_i$  in  $T/I$. 
Also note that  if  we denotes  $f(\overline{t_i t_j}) = c_{ij} \in k$, then  $f$  is completely determined by these  $c_{ij}$'s. 

We can take the following map  $\Delta$  as a lifting graded homomorphism of $\delta$  on $F \otimes _k T/I$. 
$$
\Delta = d \otimes 1 + \sum _{i=1}^r t_i^* \otimes \overline{t_i} + 0 \otimes \epsilon.   
$$
Then, we have 
$\Delta ^2 = \sum_{i, j =1}^r t_i^* t_j^*  \otimes \overline{t_it_j}$. 
Since  $\overline{t_it_j} = f(\overline{t_it_j})\epsilon = c_{ij}\epsilon$, it follows that 
$\Delta ^2 = \sum _{i,j =1}^r c_{ij} t_i^* t_j ^* \otimes \epsilon$. 
Therefore, from the definition of  $\alpha _{\Lx}$,   we see 
$\alpha _{\Lx} ([T/I, \epsilon]) = \sum _{i,j=1}^r c_{ij} [t_i^*][t_j^*]$, 
which is in fact an element of  $\Ext _R ^1(\Fx, \Fx)^2$. 
Since we can take any elements of $k$ as  $c_{ij}$, the image of  $\alpha _{\Lx}$  is exactly    $\Ext _R ^1(\Fx, \Fx)^2$. 
\qed\end{pf}

\begin{thm}\label{square iso}
Let  $\Lx _0 = ( F \ctensor P_0, \Delta _0)$  be the universal lift of  $\Fx$  with the parameter algebra $P_0 = T/I_0$,  where   $T = \formal$  is a non-commutative formal power series ring and  $I_0 \subseteq \m_T^2$. 
Then the image of the injective map  $\alpha _{\Lx _0} : \T (P_0) \to \Ext _R^2(\Fx, \Fx)$  contains  $\Ext _R^1(\Fx, \Fx)^2$.
And there is an isomorphism of $k$-vector spaces
$$
\Ext _R^1 (\Fx, \Fx)^2 \cong \Hom _k (I_0/I_0\cap \m_T^3, k). 
$$ 
\end{thm}

\begin{pf}
Let $\Lx$  be the lifting chain complex  $(F \otimes _k T/\m_T^2, \delta)$ of  $\Fx$, where  $\delta = d \otimes 1 + \sum _{i=1}^r t_i^*\otimes t_i$ as above. 
And let  $\Lx _0 = (F \ctensor T/I_0, \Delta _0)$  be the universal lift of  $\Fx$. 
We denote by  $q$  the natural injection  $I_0 \to \m_T^2$  and by  $p$  the projection  $T/I_0 \to T/\m_T^2$. 

Combining all the results  in  \ref{con iso T}, \ref{comm of alpha}  and  \ref{image of T/m^2},  we have the following commutative diagram.
\begin{equation}\label{comm diagram}
\begin{CD}
\Hom _{con} (\m_T^2, k) @>{\tau}>{\cong}> \T (T/\m_T^2) @>{\alpha _{\Lx}}>> \Ext _R ^1(\Fx, \Fx)^2 \\
@V{q^*}VV    @V{p^*}VV    @V{\iota}VV   \\
\Hom _{con} (I_0, k) @>{\tau}>{\cong}> \T (T/I_0) @>{\alpha _{\Lx _0}}>> \Ext _R ^2(\Fx, \Fx), \\
\end{CD}
\end{equation}
where  $\iota$  is a natural injection. 
Note from Theorem  \ref{injection} and Lemma  \ref{image of T/m^2}  that $\alpha _{\Lx _0}$ is injective, and  $\alpha _{\Lx}$  is surjective. 
Therefore  $\alpha _{\Lx _0} (\T(P_0))$  contains  $\Ext _R^1(\Fx, \Fx)^2$  and we have an isomorphism of $k$-vector spaces 
$$
q^* (\Hom _{con} (\m _T^2, k)) \cong  \Ext _R^1(\Fx, \Fx)^2. 
$$
Note that  $\Hom _{con} (\m _T^2 , k)= \Hom _k (\m_T^2/\m_T^3, k)$. 
Hence  we may describe as follows:
$$
\begin{array}{rl}
\kernel (q^*) &= \{ f \in \Hom _k(\m_T^2/\m_T^3, \ k) \ | \ f(I_0+\m_T^3/\m_T^3 ) =0 \} \vspace{6pt} \\
&= \Hom _k (\m _T^2/I_0 + \m_T^3, \ k). 
\end{array}
$$
Thus, from the obvious exact sequence 
$$
\begin{CD}
0 @>>> I_0/I_0 \cap \m_T^3 @>>> \m_T ^2 /\m_T^3 @>>> \m_T^2/I_0 + \m_T^3 @>>> 0, \\ 
\end{CD}
$$  
we finally have 
$$
\begin{array}{rlr}
 \Ext _R^1(\Fx, \Fx)^2  &\cong q^*(\Hom _{con} (\m_T^2, k)) & \vspace{6pt}\\
 &\cong \Hom _k (\m_T^2/\m_T^3, k)/\Hom _k (\m _T^2/I_0 + \m_T^3, \ k) & \vspace{6pt}\\ 
&\cong \Hom _k (I_0/I_0 \cap \m _T^3, k ).  &\hphantom{H}\text{\qed}
\end{array}
$$
\end{pf}

Note in the theorem that  $I_0/I_0 \cap \m_T^3$  is a finite dimensional $k$-vector space, since it is a subspace of  $\m_T^2/\m_T^3$. 
As a direct consequence of the theorem we have the following corollary. 

\begin{cor}\label{cor to square iso}
Let  $P_0 = T/I_0$  be the parameter algebra of the universal lift of $\Fx$, where   $T = \formal$  is a non-commutative formal power series ring and  $I_0 \subseteq \m_T^2$. 
Then,  $I_0 \subseteq \m_T^3$  if and only if  $\Ext _R^1(\Fx, \Fx)^2 =0$. 
\end{cor}

\begin{prop}
Let  $\Lx _0= ( F \ctensor P_0, \Delta _0)$  be the universal lift of  $\Fx$  with the parameter algebra $P_0 = T/I_0$,  where   $T = \formal$  is a non-commutative formal power series ring and  $I_0 \subseteq \m_T^2$. 
Then the following two conditions are equivalent. 
\begin{itemize}
\item[$(a)$] 
The image of the mapping  $\alpha _{\Lx _0} : \T( P_0) \to \Ext_R^2(\Fx, \Fx)$  is exactly  $\Ext _R^1(\Fx, \Fx)^2$. 

\item[$(b)$] 
There exist elements  $f_1, \ldots , f_{\ell} \in I_0$  which analytically generate the ideal  $I_0$  such that they give rise to linearly independent elements in  $\m_T^2/\m_T^3$. 
\end{itemize}
\end{prop}

\begin{pf}
By the commutative diagram \ref{comm diagram}  in the proof of Theorem \ref{square iso}, we see that the condition  (a)  is equivalent to that the $k$-linear mapping  $q^* : \Hom _{con}(\m_T^2, k) \to \Hom _{con}(I_0, k)$  is surjective, where  $q : I_0 \to \m_T^2$  is a natural injection. 

To prove the implication  $(a) \Rightarrow (b)$, suppose  $q^*$  is surjective.
 Since  $\Hom _{con}(\m_T^2, k) = \Hom _k (\m_T/\m_T^2, k)$  is a finite dimensional $k$-vector space, so is 
$$
\Hom _{con} (I_0, \ k) = \bigcup _{n \geq 1} \Hom _k (I_0/\m _T I_0 + I_0 \m _T + (I_0 \cap \m _T^n), \ k). 
$$
Hence  there is an integer  $n_0 \geq 1$  such that  
$$
\m _T I_0 + I_0 \m _T + (I_0 \cap \m _T^{n_0}) =
\m _T I_0 + I_0 \m _T + (I_0 \cap \m _T^{n_0+1}) = \cdots = 
\overline{\m _T I_0 + I_0 \m _T}.
$$
Therefore $q^*$ induces a surjective mapping  $\Hom _k (\m_T^2/\m_T^3, k) \to \Hom _k (I_0/\overline{ \m _T I_0 + I_0 \m _T}, \ k)$, hence 
the natural mapping   $I_0/\overline{ \m _T I_0 + I_0 \m _T} \to \m_T^2/\m_T^3$  is injective. 
Now let us take the elements $f_1, \ldots , f_r \in I_0$  whose images in $I_0/\overline{\m_T I_0 + I_0 \m_T}$ form a $k$-basis. 
Then, the images of  $f_1, \ldots, f_r$  in  $\m_T^2/\m_T^3$  are linearly independent and it follows from Proposition \ref{nakayama} that  $I_0$  is analytically generated by  $f_1, \ldots , f_r$.

To prove  $(b) \Rightarrow (a)$, suppose we have elements  $f_1, \ldots , f_{\ell} \in I_0$  such that they give rise to linearly independent elements in  $\m_T^2/\m_T^3$  and  $I_0 = \overline{(f_1,  \ldots , f_r)}$. 
By Proposition \ref{nakayama} we may assume that the images of  $f_1, \ldots , f_r$ in  $I_0/ \overline{ \m _T I_0 + I_0 \m_T}$  form a $k$-basis. 
Note that $I_0 \cap \m _T ^3$  is a closed ideal of  $T$  containing   $\m _T I_0 + I_0 \m_T$, hence we have an inclusion relation 
$\overline{ \m _T I_0 + I_0 \m_T} \subseteq 
I_0 \cap \m _T ^3 \subseteq \m _T^3 $. 
Therefore we obtain the natural map 
$$
I_0 / \overline{ \m _T I_0 + I_0 \m_T} \overset{g}\longrightarrow 
I_0 / I_0 \cap \m _T^3 
\subseteq 
\m_T ^2 / \m _T^3.  
$$
Since $g$  maps the $k$-basis of $I_0 / \overline{ \m _T I_0 + I_0 \m_T}$ to a set of linearly  independent elements in  $\m_T ^2 / \m _T^3$,  we have the injectivity of  $g$. 
In particular, the equality  $\overline{ \m _T I_0 + I_0 \m_T} = I_0 \cap \m _T ^3$  holds. 
Hence it follows that 
$$
\m _T I_0 + I_0 \m_T + (I_0 \cap \m _T^n)  = I_0 \cap \m _T ^3, 
$$
for all  $n \geq 3$. 
Thus we have the equality $\Hom _{con} (I_0, \ k) = \Hom _k (I_0/I_0 \cap \m _T^3, \ k)$. 
Therefore the map  $q^* : \Hom _{con}(\m_T^2, k) \to \Hom _{con}(I_0, k)$ is the same as the $k$-dual of the natural injection  $I_0 / I_0 \cap \m_T^3 \subseteq \m _T^2/ \m _T^3$. 
The surjectivity of  $q^*$  is now obvious. 
\qed\end{pf}

\subsection{Comparison of cohomology} 
As in the previous sections,  $\Fx = (F, d)$  denotes a projective complex over  $R$, where  $R$  is an associative $k$-algebra. 
We assume that  $r = \Kdim _k \Ext _R^1(\Fx, \Fx)$  is finite as before.
Adding to this assumption, we assume in the rest of the paper that the complex  $\Fx$  is right bounded, i.e. there is an integer  $s$  such that  $F_i = 0$  for  $i<s$. 

We also denotes by  $\Lx _0 = (F \ctensor P_0, \Delta _0)$  the universal lift of  $\Fx$  with the parameter algebra  $P_0$. 
Note that $\Lx _0$  may not be projective as a right  $P_0$-module. 
They are even non-flat as seen in Example \ref{not flat}. 

For any integer  $n \geq 1$, we set 
$$
\Lx _0^{(n)} = (F \otimes _k P_0/\m ^n _{P_0}, \ \Delta ^{(n)}) 
= \Lx \otimes _{P_0} P_0/\m ^n _{P_0}. 
$$
In fact, each  $\Lx _0^{(n)}$ is  a right bounded complex of projective left  $R \otimes _k (P_0/\m ^n_{P_0})^{op}$-modules, in particular, it is a free right $P_0/\m^n_{P_0}$-module.

For any associative $k$-algebra $R$, we denote by  $D_{+} (R)$  the derived category consisting of right bounded complexes over $R$. 
Then, tensoring the chain complex  $\Lx _0^{(n)}$  yields the functor $\rho _{n}$  between the derived categories:  
$$
\rho _{n} : D_{+} (P_0/ \m ^n _{P_0}) \to  D_{+} (R),  
$$
which is defined by  $\rho _{n} (X) = \Lx _0 ^{(n)} \otimes _{P_0/\m ^n_{P_0}} X$. 
This is well-defined, since  $\Lx _0 ^{(n)}$  is a right bounded complex of projective left  $R \otimes _k (P_0/\m^n_{P_0})^{op}$-modules. 

Note that the natural projection $P_0/\m^{n+1}_{P_0} \to P_0/\m^n_{P_0}$ induces a natural functor  $D_+ (P_0/\m^n_{P_0}) \to D_+ (P_0/\m^{n+1}_{P_0})$. 
And it is easy to see the diagram 
\begin{equation}\label{derived} 
\begin{CD}
D_+ (P_0/\m^n_{P_0})  @>{\text{\small ${\rho _n}$}}>>  D_+ (R)  \\
@VVV 
\begin{picture}(0,0)
\put(-50,-15){\vector(1,1){25}}
\put(-30,-5){\small {$\rho _{n+1}$}}
\end{picture} @.      
\\
D_+(P_0/\m^{n+1}_{P_0}) \\
\end{CD}
\end{equation}
is commutative. 

Note that  $\Lx _0 ^{(n)} \otimes _{P_0/\m ^n_{P_0}} k = \Fx$, hence we have 
$\rho _n (k) = \Fx$  for each  $n \geq 1$.
It follows that the functor  $\rho _n$  induces the map 
$$
\rho _n ^i : \Ext _{P_0/ \m ^n _{P_0}} ^i (k,k) \to \Ext _R ^i (\Fx, \Fx)
$$
for all integer  $i$. 
Then the commutative diagram (\ref{derived}) forces the commutativity of the following diagram. 
\begin{equation}\label{ext} 
\begin{CD}
\Ext _{P_0/\m^n_{P_0}}^i (k,k) @>{\text{\small ${\rho _n^i}$}}>>  \Ext _R^i (\Fx, \Fx)  \\
@VVV 
\begin{picture}(0,0)
\put(-50,-15){\vector(1,1){25}}
\put(-30,-5){\small {$\rho _{n+1}^i$}}
\end{picture} @.      
\\
\Ext _{P_0/\m^{n+1}_{P_0}}^i (k,k) \\
\end{CD}
\end{equation}

\begin{defn}
From the commutative diagram (\ref{ext}), we can define the inductive limit  $\rho _{\infty} ^i = \varinjlim_n \rho _n ^i $ for each  $i$ ; 
$$
\rho _{\infty} ^i : \varinjlim _n \Ext_{P_0/\m^{n}_{P_0}}^i(k, k) \to \Ext _R ^i (\Fx, \Fx).
$$  
\end{defn}

The aim of this section is to show that  $\rho _{\infty}^i$  is an injective map for  $i = 0, 1, 2$. 

Note that  $\varinjlim _n \Ext_{P_0/\m^{n}_{P_0}}^{\cdot}(k, k)=
\bigoplus _{i\geq 0} \varinjlim _n \Ext_{P_0/\m^{n}_{P_0}}^i(k, k)$,  as well as  $\Ext _R ^{\cdot} (\Fx, \Fx)$,  has a structure of algebra by Yoneda product, and  
$$
\rho _{\infty} ^{\cdot} : \varinjlim _n \Ext_{P_0/\m^{n}_{P_0}}^{\cdot}(k, k) \to \Ext _R ^{\cdot} (\Fx, \Fx)
$$
is an algebra map. 

First, consider the case $i = 0$. 
Since  $\Ext _{P_0/\m^{n}_{P_0}}^0(k, k) = k$  and  $\rho _n ^0 : k \to \Ext _R ^0(\Fx, \Fx)$  is a natural injection for any  $n \geq 1$, we easily see the following lemma holds. 

\begin{lemma}
The mapping  $\rho _{\infty} ^0 : k \to \Ext _R ^0 (\Fx, \Fx)$  is a natural injection. 
\end{lemma}
 
To argue for the case  $i=1$, we should notice that 
$\Ext _{P_0/\m^{n}_{P_0}}^1(k, k) \cong  \Hom _k (\m _{P_0}/ \m ^2_{P_0}, k)$  for all $n \geq 2$  and  the natural maps  $\Ext _{P_0/\m^{n}_{P_0}}^1(k, k) \to \Ext _{P_0/\m^{n+1}_{P_0}}^1(k, k)$ coincide with the identity map on   $\Hom _k (\m _{P_0}/ \m ^2_{P_0}, k)$. 
Hence we have 
$$
\rho _{\infty} ^1 :  \Hom _k (\m _{P_0}/ \m ^2_{P_0}, k) \to \Ext _R ^1(\Fx, \Fx). 
$$   
We can prove this is actually an isomorphism.

\begin{lemma}
The mapping  $\rho _{\infty} ^1 : \Hom _k (\m _{P_0}/ \m ^2_{P_0}, k) \to \Ext _R ^1 (\Fx, \Fx)$  is an isomorphism. 
\end{lemma}

\begin{pf}
By the observation above, we only have to prove that 
$\rho _2^1 : \Hom _k (\m _{P_0}/ \m ^2_{P_0}, k) \to \Ext _R ^1(\Fx, \Fx)$  is an isomorphism.
Let us denote  $P_0 = T/ I_0$  where  $T = \formal$  and  $I_0 \subseteq \m _T^2$. 
Recall that  $\Ext _R ^1(\Fx, \Fx)$  has  a $k$-basis  $\{ [t_1^*], [t_2^*], \ldots , [t_r^*]\}$. 
And, by definition, we have 
$\Lx _0 ^{(2)} = (F \otimes _k T/\m_T^2, \ \delta)$  with 
$\delta = d \otimes 1 + \sum _{i=1}^r t_i^* \otimes \overline{t_i}$, where $\overline{t_i}$  denotes the image of  $t_i$  in  $\m_T/\m_T^2$. 
Therefore, it is easy to see that the mapping  $\rho _2^1$  is defined by 
$$
\rho_2^1  (f )  = \sum _{i=1} ^r f(\overline{t_i}) [t_i^*], 
$$
for  $f \in \Hom _k (\m_T/ \m_T^2, k)$. 
The lemma follows from this.  
\qed\end{pf}

Now we proceed to the case $i=2$. 
The goal here is to prove the following theorem. 

\begin{thm}\label{main comparison} 
There is an isomorphism  $\beta : \T (P_0) \to \varinjlim \Ext _{P_0/\m _{P_0}^n} ^2(k, k)$  which makes the following diagram commutative. 
$$
\begin{CD}
\T (P_0)  @>{\text{\small ${\alpha_{\Lx _0}}$}}>>  \Ext _R^2 (\Fx, \Fx)  \\
@V{\text{\small ${\beta}$}}VV 
\begin{picture}(0,0)
\put(-55,-15){\vector(1,1){25}}
\put(-30,-5){\small {$\rho _{\infty}^2$}}
\end{picture} @.      
\\
\varinjlim _n \Ext _{P_0/\m_{P_0}^n}^2 (k,k)  \\
\end{CD}
$$
In particular,  $\rho _{\infty}^2$  is an injective map as is $\alpha_{\Lx _0}$. \end{thm}

To prove this, let  $A$  be an arbitrary  artinian local $k$-algebra in $\A$, and let  $\Gx ^A = (G, d^A)$  be a free resolution of the left $A$-module $k = A/\m_A$. 
Then, by Theorem \ref{resol of k}, there is a universal lift  of  $\Gx ^A$  of the form  $\Gx _0^A = (G \otimes _k A, \Delta _0^A)$ whose parameter algebra is  $A$. 
Since  $\Gx _0 ^A$  is a lift of  $\Gx ^A$, we have the $k$-linear map 
$$
\alpha _{\Gx _0^A} : \T (A) \to  \Ext _{A} ^2(\Gx ^A, \Gx ^A) = \Ext _A^2(k, k)
$$
which is defined in \ref{def of alpha}.

\begin{lemma}\label{alpha is iso}
Let  $A \in \A$  as above.
Then the map  $\alpha _{\Gx _0^A}$  is an isomorphism. 
\end{lemma}

\begin{pf}
Since $A$  is a parameter algebra of the universal lift  $\Gx _0^A$  of $\Gx$, Theorem \ref{injection} implies that  $\alpha _{\Gx _0^A}$  is injective. 
Thus, to show this is an isomorphism, it is enough to show that 
$\Kdim _k \T (A) =  \Kdim _k \Ext _A^2(k, k)$. 
Let us describe  $A = T/I$  where  $T$  is the non-commutative formal power series ring and  $I \subseteq \m _T ^2$. 
Note that the ideal  $I$ is open and closed in  $T$, since  $A$ is artinian. 
Therefore, by Proposition \ref{con iso T},  we have  $\T(A) \cong \Hom _{con}(I, k) = \Hom _k (I/I \m_A + \m_A I, k)$. 
On the other hand, by the following lemma, 
 we know that   $\Ext _A^2(k ,k) \cong  \Hom _k (I/I\m _T + \m_T I, k) \cong \T (A)$.
Hence  $\Kdim _k \Ext _A^2(k, k) = \Kdim _k \T (A)$.
\qed\end{pf}

\begin{lemma}\label{second syzygy}
Let  $P = T/I \in \Ahat$  where  $T= \formal$  is a non-commutative formal power series ring and  $I \subseteq \m _T^2$. 
Then there is an isomorphism 
$\Ext _P ^2 (k, k) \cong \Hom _k (I/I\m_T+\m_T I, k) \cong \Hom _{\text{$T$-bimod}}(I, k)$. 
\end{lemma}

\begin{pf}
By virtue of Lemma \ref{f g ideals are free},  there is a minimal free resolution of $k$  as a left $T$-module of the form 
$$
\begin{CD}
0 @>>> T^r @>>> T @>>> k @>>> 0.
\end{CD}
$$ 
Therefore, tensoring  $P$  over $T$, we have an exact sequence of left $P$-modules 
$$
\begin{CD}
0 @>>> \Tor _1^T(P, k) @>>> P^r @>>> P @>>> k @>>> 0.
\end{CD}
$$ 
Note that, by the exact sequence of right $T$-modules 
$0 \to I \to T \to  P \to 0$, 
we have  $ \Tor _1^T(P, k) \cong I/I\m_T$ which is an isomorphism of left $P$-modules.  
Therefore we have 
$\Ext _P ^2 (k, k) \cong \Hom _P (I/I \m_T, k) = \Hom _k (I/I\m_T+\m_T I, k)$. \qed\end{pf}

One can show that the isomorphism $\alpha _{\Gx _0^A}$ does not depend on the choice of free resolution $\Gx ^A$ and its lift  $\Gx _0 ^A$. 
This follows from the following more general lemma.

\begin{lemma}\label{independence}
Let  $\Fx^{(1)} =(F^{(1)}, d^{(1)})$ and  $\Fx^{(2)} =(F^{(2)}, d^{(2)})$ be right bounded projective  complexes over $R$. 
For  $A \in \A$, let  $\Gx ^{(i)} = (F ^{(i)}\otimes _k A, \Delta ^{(i)})$  be a lift  of  $\Fx ^{(i)}$ to  $A$  for $i=1, 2$. 
Suppose that there is a quasi-isomorphism  $q : \Gx ^{(1)} \to \Gx ^{(2)}$ of chain complexes over $R \otimes _k A^{op}$. 
Then, there is a commutative diagram:
$$
\begin{CD}
\T (A) @>{\alpha_{\Gx ^{(1)}}}>> \Ext _R^2 (\Fx ^{(1)}, \Fx ^{(1)}) \\
@V{\alpha_{\Gx ^{(2)}}}VV @V{q_*}VV \\
\Ext _R^2 (\Fx ^{(2)}, \Fx ^{(2)}) @>{q^*}>> \Ext _R ^2(\Fx ^{(1)}, \Fx ^{(2)})  \\
\end{CD}
$$
\end{lemma}

\begin{pf}
Let  $[A', \epsilon] \in \T(A)$ and let  $\widetilde{\Delta}^{(i)}$  be a lifting homomorphism  $F^{(i)} \otimes _k A' \to F^{(i)} \otimes _k A'[-1]$ 
of  $\Delta ^{(i)}$  for  $i =1,2$. 
Then, by the definition of  $\alpha _{\Gx ^{(i)}}$, we have the description 
$(\widetilde{\Delta}^{(i)})^2 = h^{(i)} \otimes _k \epsilon$ with $h^{(i)} : F ^{(i)} \to F^{(i)} [-2]$ being a chain homomorphism,  and the equality 
$\alpha _{\Gx ^{(i)}} ([A', \epsilon]) = [h^{(i)}]$ holds. 
Now take a lifting map   $q' : F^{(1)} \otimes _k A' \to F^{(2)} \otimes _k A'$  of a graded homomorphism  $q$  and we have the commutative diagram: 
$$
\begin{CD}
0 @>>> F^{(1)} @>{1 \otimes \epsilon}>> F^{(1)} \otimes _k A' @>>> F^{(1)} \otimes _k A @>>> 0 \\
@. @V{q \otimes _A k}VV @V{q'}VV @V{q}VV \\
0 @>>> F^{(2)} @>{1 \otimes \epsilon}>> F^{(2)} \otimes _k A' @>>> F^{(2)} \otimes _k A @>>> 0 \\
\end{CD}
$$   
Since  $\Delta ^{(2)} q = q \Delta ^{(1)}$, we have the following commutative diagram: 
$$
\begin{CD}
0 @>>> F^{(1)} @>{1 \otimes \epsilon}>> F^{(1)} \otimes _k A' @>>> F^{(1)} \otimes _k A @>>> 0 \\
@. @VV{0}V @VV{\widetilde{\Delta} ^{(2)} q' - q' \widetilde{\Delta} ^{(1)}}V @VV{0}V \\
0 @>>> F^{(2)}[-1] @>{1 \otimes \epsilon}>> F^{(2)} \otimes _k A'[-1] @>>> F^{(2)} \otimes _k A[-1] @>>> 0 \\
\end{CD}
$$   
Hence there is a graded homomorphism  $p : F^{(1)} \to F^{(2)}[-1]$  with the equality 
$$
\widetilde{\Delta} ^{(2)} q' - q' \widetilde{\Delta} ^{(1)}=p \otimes \epsilon.
$$
Multiplying  $\widetilde{\Delta} ^{(2)}$ (resp. $\widetilde{\Delta} ^{(1)}$)  from the left (resp. right),  we have equalities 
$$
d^{(2)}p \otimes \epsilon = \widetilde{\Delta} ^{(2)} (p \otimes \epsilon) = 
(h^{(2)}\otimes \epsilon) q' - \widetilde{\Delta} ^{(2)} q' \widetilde{\Delta} ^{(1)} = h^{(2)} q \otimes \epsilon - \widetilde{\Delta} ^{(2)}q' \widetilde{\Delta} ^{(1)},  
$$
and 
$$
pd^{(1)} \otimes \epsilon = (p \otimes \epsilon)\widetilde{\Delta} ^{(1)} 
= \widetilde{\Delta} ^{(2)} q' \widetilde{\Delta} ^{(1)} - q' (h^{(1)}\otimes \epsilon)  = \widetilde{\Delta} ^{(2)}q' \widetilde{\Delta} ^{(1)} - q h^{(1)}  \otimes \epsilon.  
$$ 
Consequently, the equality  
$$
h^{(2)}q -q h^{(1)} = d^{(2)}p + pd^{(1)}
$$
holds. 
It follows that the equality  $q_*([h^{(1)}]) = [q h^{(1)}] = [h^{(2)}q] = q^* ([h^{(2)}])$ holds as an elements of  $\Ext _R^2(\Fx ^{(1)}, \Fx^{(2)})$. 
\qed\end{pf}

\begin{lemma}\label{comm}
Let  $A \in \A$  and  let  $\Gx ^A = (G, d^A)$  be a free resolution of the left $A$-module $k = A/\m_A$. 
We take a universal lift  of  $\Gx ^A$  of the form  $\Gx _0^A = (G \otimes _k A, \Delta _0^A)$ as above. 
Furthermore, suppose that there is a lifting chain complex  $\Lx = (F \otimes _k A, \Delta )$ of  $\Fx = (F,d)$. 
Then we have the following commutative diagram 
$$
\begin{CD}
\T (A)  @>{\text{\small ${\alpha _{\Lx}}$}}>>  \Ext _R^2 (\Fx, \Fx)  \\
@V{\text{\small ${\alpha _{\Gx _0^A}}$}}VV 
\begin{picture}(0,0)
\put(-50,-15){\vector(1,1){25}}
\put(-30,-5){\small {$\rho _{\Lx}^2$}}
\end{picture} @.      
\\
\Ext _{A}^2 (k,k),  \\
\end{CD}
$$
where  $\rho _{\Lx}^2$  is the map induced by the functor  $\Lx \otimes _A - : D_+(A) \to D_+(R)$. 
\end{lemma}

\begin{pf}
Consider the tensor product of chain complexes 
$$
\Xx : = \Lx \otimes _A \Gx _0^A = ((F \otimes _k A)\otimes _A (G \otimes _k A), \  d_X), 
$$ 
where  $d_X = \Delta \otimes 1 + 1 \otimes \Delta _0^A$. 
Notice from Theorem \ref{resol of k} that $\Gx _0^A$  is a complex of free $A \otimes _k A ^{op}$-modules and it is quasi-isomorphic to  $A$  as a chain complex of  $A\otimes A^{op}$-modules. 
Therefore the chain complex   $\Xx$  is quasi-isomorphic to  $\Lx$  as a chain complex of  $R \otimes _k A^{op}$-modules. 
By virtue of Lemma \ref{independence}, it is sufficient to prove the commutativity of the following diagram. 
$$
\begin{CD}
\T (A)  @>{\text{\small ${\alpha _{\Xx}}$}}>>  \Ext _R^2 (\Xx \otimes _A k , \Xx \otimes _A k)  \\
@V{\text{\small ${\alpha _{\Gx _0^A}}$}}VV 
\begin{picture}(0,0)
\put(-85,-15){\vector(1,1){25}}
\put(-60,-5){\small {$\rho _{\Xx}^2$}}
\end{picture} @.      
\\
\Ext _{A}^2 (k,k) \\
\end{CD}
$$
For this, let  $[A', \epsilon] \in \T (A)$ and take a lifting map $\Gamma : G\otimes _k A' \to G \otimes _k A'[-1]$  of  $\Delta _0 ^A$. 
By definition of  $\alpha _{\Gx _0^A}$, we have 
$\Gamma ^2 = h \otimes \epsilon$  for some  $h : G \to G[-2]$  and  $\alpha _{\Gx _0^A} ([A', \epsilon]) = [h]$. 
Therefore  $\rho _{\Xx}^2 (\alpha _{\Gx _0^A}([A', \epsilon]))$  is represented by a chain map  $1 \otimes h$  on  $\Xx \otimes _A k = \Fx \otimes _A \Gx ^A$. 
On the other hand, $d _X ' = \Delta \otimes 1 + 1 \otimes \Gamma$ is a lifting map of  $d_X$, and we have the equality 
${d_X'}^2 = 1 \otimes h \otimes \epsilon$. 
Hence it follows from the definition that  $\alpha _{\Xx} ([A', \epsilon]) = [1 \otimes h]$ as well.
Hence we have  $\rho _{\Xx}^2\cdot \alpha _{\Gx _0^A} = \alpha _{\Xx}$. 
\qed\end{pf}

Let  $P \in \Ahat$ be any complete local $k$-algebra. 
Then, induced from the natural projections  $p _n : P/ \m_P ^{n+1} \to P/\m _P ^{n}$  and  $\pi _n : P \to P/ \m_P ^n$  for each  $n \geq 1$, there are natural mappings  $p _n^* : \T (P/ \m_P ^n) \to \T (P/\m _P ^{n+1})$ and  $\pi _n^* : \T (P/\m_P^n) \to \T (P)$. 
See Lemma \ref{functoriality of T}. 
By the functorial property of $\T$, it is clear that the diagram 
$$
\begin{CD}
\T (P/\m_P^n)  @>{\text{\small ${\pi ^* _n}$}}>>  \T (P) \\
@V{\text{\small ${p ^* _n}$}}VV 
\begin{picture}(0,0)
\put(-40,-15){\vector(1,1){25}}
\put(-20,-5){\small {$\pi ^* _{n+1}$}}
\end{picture} @.      
\\
\T (P/\m_P^{n+1}) \\
\end{CD}
$$
is commutative for each $n \geq 1$. 
Thus it induces the map 
$$
\gamma _P : = \varinjlim \pi _n^* \ : \ \varinjlim \T (P/\m_P^n) \to \T (P). 
$$

\begin{lemma}\label{gamma is iso}
The map $\gamma _P$ is an isomorphism for any  $P \in \Ahat$. 
\end{lemma}

\begin{pf}
First we show that each  $\pi _n ^* : \T (P/\m_P^n) \to \T (P)$  is injective for $n \geq 2$,  hence so is $\gamma _P$. 
For this, let  $[A, \epsilon] \in \T (P/\m_P^n)$. 
Then take a fiber product and we have the following commutative diagram with exact rows and columns. 
$$
\begin{CD}
@. @. 0 @. 0 \\
@. @. @VVV @VVV \\
@. @. \m_{A'}^n @= \m _P^n \\
@. @. @VVV @VVV \\
0 @>>> k @>{\epsilon'}>>  A'  @>>> P @>>> 0 \\
@. @| @VVV  @V{\pi _n}VV \\
0 @>>> k @>{\epsilon}>>  A  @>>> P/\m_P^n @>>> 0 \\
@. @. @VVV @VVV \\
@. @. 0 @. 0 \\
\end{CD}
$$ 
By definition  $\pi ^* ([A, \epsilon]) = [A', \epsilon ']$. 
Suppose  $[A', \epsilon ']=0$  in  $\T (P)$. 
Then, since the small extension  $(A', \epsilon ')$  is a trivial one, 
 we have  $\epsilon ' \not\in \m_{A'}^2$  by Lemma \ref{trivial or nontrivial small extensions}. 
Then by the diagram above, we see  $\epsilon \not\in \m_A ^2$ as well. 
Hence  $[A, \epsilon] =0$  in  $\T(P/\m_P^n)$ again by Lemma \ref{trivial or nontrivial small extensions}. 

Now we prove  $\gamma _P : \varinjlim \T(P/\m_P ^n) \to \T (P)$  is surjective.  For this, let  $[A', \epsilon ']$  be any element of $\T(P)$. 
Since  $\bigcap _{n=1}^{\infty} \m _{A'}^n =(0)$  and since $(\epsilon ')$ is of finite length, there is an integer  $n_0 \geq 1$ such that $(\epsilon ') \cap \m_{A'} ^n = (0)$  for  $n \geq n_0$. 
For such any $n$, we set  $A _n = A'/\m_{A'}^n$  and  $\epsilon _n = \epsilon ' \mod \m_{A'}^n$. 
And it is easy to see that $[A_n , \epsilon _n] \in \T (P/\m _P ^n)$  and   $\pi _n^* ([A_n , \epsilon _n]) = [A', \epsilon ']$  for  $n \geq n_0$. 
The surjectivity of  $\gamma _P$  follows from this. 
\qed\end{pf}

\vspace{4pt}
\noindent
[Proof of Theorem \ref{main comparison}]

Let  $\Lx _0$  be a universal lift of  $\Fx$  with the parameter algebra  $P _0$ as in the setting of the theorem. 
We denote  $\Lx _0 ^{(n)} = \Lx \otimes _{P_0} P_0/\m_{P_0} ^n$  and  $\Gx _0 ^{(n)} = \Gx _0 \otimes _{P_0} P_0/\m_{P_0} ^n$, where  
$\Gx _0$  is the universal lift of a free left $P_0$-module $k$. 
Then, from Lemma \ref{comm},  we have a commutative diagram 
$$
\begin{CD}
\T (P_0/\m_{P_0} ^n )  @>{\text{\small ${\alpha _{\Lx _0 ^{(n)}}}$}}>>  \Ext _R^2 (\Fx, \Fx)  \\
@V{\text{\small ${\alpha _{\Gx _0^{(n)}}}$}}VV 
\begin{picture}(0,0)
\put(-50,-15){\vector(1,1){25}}
\put(-30,-5){\small {$\rho _{\Lx _0^{(n)}}^2$}}
\end{picture} @.      
\\
\Ext _{P_0/\m_{P_0} ^n }^2 (k,k). \\
\end{CD}
$$ 
Note from Lemma \ref{alpha is iso} that $\alpha _{\Gx _0^{(n)}}$ is an isomorphism. 
Now taking the inductive limit and setting  $\beta = \varinjlim \alpha _{\Gx _0^{(n)}}$, we have a commutative diagram  by Lemma \ref{gamma is iso} 
$$
\begin{CD}
\T (P_0)  @>{\text{\small ${\varinjlim \alpha _{\Lx _0 ^{(n)}}}$}}>>  \Ext _R^2 (\Fx, \Fx)  \\
@V{\text{\small ${\beta}$}}VV 
\begin{picture}(0,0)
\put(-60,-15){\vector(1,1){25}}
\put(-30,-5){\small $\varinjlim \rho ^2_{\Lx _0^{(n)}}$}
\end{picture} @.      
\\
\varinjlim  \Ext _{P_0/\m_{P_0} ^n }^2 (k,k),  \\
\end{CD}
$$
where  $\beta$  is an isomorphism as well. 
It is easy to see from the definition that 
$\varinjlim \alpha _{\Lx _0^{(n)}} = \alpha _{\Lx _0}$ 
and 
$\varinjlim \rho ^2_{\Lx _0^{(n)}} = \rho ^2 _{\infty}$. 
\qed

\vspace{6pt}

We should note that there is  a natural mapping 
$$
\nu : \varinjlim _n \Ext _{P/\m _P^n}^2 (k, k) \to \Ext _P ^2 (k, k). 
$$
However,  the mapping  $\nu$  is not an isomorphism in general. 
In fact, we can show the following proposition.

\begin{prop}
Let  $P = T/I$  be a complete local $k$-algebra where  $T$  is a non-commutative formal power series ring  and  $I \subseteq \m_T^2$. 
Then the natural map  $\nu$  is always injective. 
It is an isomorphism if and only if 
the ideal $\m _T I + I \m_T$  is closed and  
$\Kdim _k I/\m_T I + I \m_T$  is finite. 
\end{prop}

\begin{pf}
By Lemma \ref{second syzygy}, we know that 
$\Ext _P^2(k ,k) \cong \Hom _{\text{$T$-bimod}}(I, k)$. 
On the other hand, it follows from  Theorem \ref{main comparison} and Proposition  \ref{con iso T} that  
\begin{equation}\label{injlim}
\varinjlim _n \Ext _{P/\m_P^n}^2(k ,k) \overset{\beta}\cong \T (P) \cong \Hom _{con}(I, k).  
\end{equation}
Through these isomorphisms, it can be seen that  $\nu$  coincides with the natural map $\Hom _{con}(I, k)  \to \Hom _{\text{$T$-bimod}}(I, k)$, which is of course an injection. 

Suppose  that   $\m _T I + I \m_T$  is closed with $\Kdim _k I/\m_T I + I \m_T < \infty$. 
Then, by Corollary \ref{artin-rees}, the inclusion  $\m _T^n \cap I \subseteq \m_T I + I \m_T$ holds  for large $n \gg 1$. 
Therefore we have  $\Hom _{con}(I, k) = \Hom _{\text{$T$-bimod}}(I, k)$.
See Definition \ref{def of con}. 

On the contrary, assume $\Hom _{con}(I, k) = \Hom _{\text{$T$-bimod}}(I, k)$. 
If  $\Kdim _k I/\m_T I + I \m_T = \infty$, then 
$\Hom_{\text{$T$-bimod}}(I, k) = \Hom _k (I/\m_T I + I \m_T , k)$ has uncountable dimension as a $k$-vector space. 
On th other hand, by the equality (\ref{injlim}), 
$\Hom _{con}(I, k)$  has countable dimension, as it is an inductive limit of finite dimensional $k$-vector spaces. 
By this contradiction, we can conclude that  $\Kdim _k I/\m_T I + I \m_T < \infty$. 
Then, since  $\Kdim _k I/\overline{\m_T I + I \m_T} < \infty$, it follows 
$$
\Hom _{con} (I, k) = \Hom _k ( I/\overline{\m_T I + I \m_T}, k). 
$$  
Since this equals $\Hom _{\text{$T$-bimod}}(I, k) = \Hom _k (I/{\m_T I + I \m_T}, k)$, we see the equality $\overline{\m_T I + I \m_T}= {\m_T I + I \m_T}$.   
\qed\end{pf}


\end{document}